\documentclass[12pt,reqno]{amsart}

\allowdisplaybreaks

\newtheorem{Theorem}{Theorem}

\newtheorem{Lemma}[Theorem]{Lemma}

\theoremstyle{remark}

\newtheorem*{Remark}{Remark}

\numberwithin{equation}{section}

\def\LL{\leavevmode\setbox0=\hbox{L}\hbox to\wd0{\hss\char'40L}}
\def\al{\alpha}

\def\ep{\varepsilon}
\def\ze{\zeta}
\def\et{\eta}
\def\th{\theta}

\def\la{\lambda}

\def\si{\sigma}

\def\om{\omega}

\def\today{\ifcase\month\or
 January\or February\or March\or April\or May\or June\or
 July\or August\or September\or October\or November\or December\fi
 \space\number\day, \number\year}
\def\dfrac#1#2{\frac{\displaystyle#1}{\displaystyle#2}}

\def\({\left(}
\def\){\right)}
\def\[{\left[}
\def\]{\right]}

\def\Re{\operatorname{\mathfrak R\mathfrak e}}

\def\sgn{\operatorname{sgn}}

\def\prodl{\prod\limits}

\def\3{\ss}
\def\sp{\operatorname{\text {\it sp}}}
\def\so{\operatorname{\text {\it so}}}
\def\soe{\operatorname{\widetilde{\text {\it so}}}}
\def\fl#1{\left\lfloor#1\right\rfloor}
\def\cl#1{\left\lceil#1\right\rceil}
\def\coef#1{\left\langle#1\right\rangle}
\def\poq#1#2{(#1;q)_#2}
\def\oddrows{\operatorname{oddrows}}

\begin{document}

\newbox\Adr
\setbox\Adr\vbox{
\centerline{\sc C.~Krattenthaler$^\dagger$}
\vskip18pt
\centerline{Institut Girard Desargues,
Universit\'e Claude Bernard Lyon-I}
\centerline{21, avenue Claude Bernard, F-69622 Villeurbanne Cedex, France}
\centerline{WWW: \footnotesize\tt http://igd.univ-lyon1.fr/\~{}kratt}
}

\title{Asymptotics for random walks in alcoves of affine Weyl groups}
\author[C. Krattenthaler]{\box\Adr}
\address{Institut Girard Desargues,
Universit\'e Claude Bernard Lyon-I,
21, avenue Claude Bernard,
F-69622 Villeurbanne Cedex, France}
\thanks{$^\dagger$ Research partially supported by the Austrian
Science Foundation FWF, grant P13190-MAT, as well as
grant S9607-N13,
in the framework of the National Research Network
``Analytic Combinatorics and Probabilistic Number Theory,"
and by EC's IHRP Programme, grant HPRN-CT-2001-00272, 
``Algebraic Combinatorics in Europe."\newline
\leavevmode\indent
Current address: Fakult\"at f\"ur Mathematik, Universit\"at Wien,
 Nordbergstra{\ss}e~15, A-1090 Vienna, Austria.
}

\dedicatory{D\'edi\'e \`a Alain Lascoux}



\begin{abstract}
Asymptotic results are derived for the number of random walks in
alcoves of affine Weyl groups (which are certain regions in
$n$-dimensional Euclidean space bounded by hyperplanes), 
thus solving problems posed by Grabiner 
[J. Combin.\ Theory Ser.~A {\bf 97} (2002), 285--306]. 
These results include asymptotic
expressions for the number of vicious walkers on a circle, as well as for
the number of vicious walkers in an interval. 
The proofs depart from the exact results of Grabiner [loc.\ cit.], 
and require as diverse means as  
results from symmetric function theory and the saddle point method,
among others.
\end{abstract}
\maketitle

\section{Introduction} \label{sec:1}
The enumeration of
random walks in lattice regions bounded by hyperplanes is a classical and
frequently studied subject in combinatorics and related fields. 
Its attractiveness stems from the fact that this problem has
implications to many other, often seemingly unrelated problems, 
and, thus, to several different fields. 
To mention some
examples, random walk interpretations exist for
ballot problems (see e.g.\ \cite{FilaAA,MohaAE,ZeilAD}), 
standard Young tableaux (see e.g.\ \cite{ZeilAD}), 
semistandard tableaux and plane partitions 
(see e.g.\ \cite{GuOVAA,KarlAB,KrGVAA} and 
\cite[Sec.~8]{ProcAJ}), 
symplectic tableaux (see e.g.\ \cite{FuKrAA,KrGVAA}),
oscillating tableaux (see e.g.\ \cite{BaFoAJ,KratAV}),
cylindric partitions (see \cite{GeKrAA}),
non-intersecting lattice paths and
vicious walkers (see e.g.\ \cite{FIS,GrabAD,GuOVAA,KrGVAA,KrGVAB}), 
and are therefore used for the solution of problems in these areas
(see e.g.\ \cite{GrabAE,GrMrAA,ProcAJ} for applications in
representation theory, and e.g.\ \cite{FIS} 
for applications in statistical physics),
as well as in the analysis of non-parametric statistics
in probability theory (see \cite{MohaAE} for an introduction to that
area).

Clearly, at the very beginning stands the problem of enumerating all
lattice paths in the plane integer lattice 
between two given points, which consist of positive unit steps and do
not cross a given diagonal line. (In fact, the original formulation
is in terms of a two-candidate ballot problem.) As is well-known, a
solution to this problem is provided by the famous {\it reflection
principle}, which is usually attributed to Andr\'e \cite{AndDAA}
(see e.g.\ \cite[p.~22]{ComtAA}). It is more than a hundred years
later, when Gessel and Zeilberger showed in \cite{GeZeAA} how far one
can go by using the reflection principle. Their main result gives 
(under certain conditions) the
number of random walks in regions of $n$-dimensional
Euclidean space which are bounded by hyperplanes. Their formula
involves the elements of the reflection group which is generated by
the reflections with respect to the hyperplanes which bound these regions. (The
same formula for the case of finite Weyl groups has been
independently discovered by Biane \cite{BianAB}. We refer the reader
to \cite{HumpAC} for an introduction to reflection groups and Weyl
groups.) 
It covers numerous formulae that occurred in the literature earlier
(and even afterwards \dots).\footnote{There are in fact only very few
known results on the enumeration of walks in regions bounded by hyperplanes
that are not covered by this result. For earlier results see
\cite[Ch.~1 and 2]{MohaAE}. More recent results include for example
\cite{BaFlAA,HaMoAB,KrNhAA,NiedAC,NiedAJ,NiedAN,SatoAA,SatoAB,SaSaAA},
where, from a conceptual point of view, the papers \cite{BaFlAA} and
\cite{NiedAC,NiedAJ,NiedAN} have to be emphasized most: in
\cite{BaFlAA} the so-called {\it kernel method\/} is exploited (which
seems to be especially well-suited for this type of problem), whereas
in \cite{NiedAC,NiedAJ,NiedAN} it is the {\it umbral calculus} which
is systematically applied to solve lattice path enumeration problems.}

Recently, Grabiner \cite{GrabAD} has revisited the
problem of enumerating random walks in {\it alcoves} of {\it affine}
(i.e., infinite) {\it
Weyl groups}. (See the next section for precise definitions.) 
To be precise, he considered three types of random walks in
these regions: (1) 
lattice walks consisting of {\it positive unit steps} $e_j$
(with $e_j$ denoting the $j$-th standard unit vector), 
(2) lattice walks
consisting of positive {\it and negative} unit steps $\pm e_j$, and
(3) lattice walks consisting of steps of the form $\pm\frac {1}
{2}e_1\pm\frac {1} {2}e_2\pm\dots \pm\frac {1} {2}e_n$ 
(where any sign pattern is
allowed). We will refer to these three types of 
walks as {\it walks with standard steps
in the positive direction} (or {\it walks with positive standard
steps}, for short), {\it walks with standard steps}, and {\it
walks with diagonal steps}, respectively. Starting from the result of
Gessel and Zeilberger, Grabiner derived interesting determinantal formulae for
the enumeration of these three types of walks in alcoves 
of types $\tilde A_{n-1}$, $\tilde B_n$, $\tilde C_n$, and $\tilde D_n$. 
From the $\tilde A_{n-1}$
results he was also able to derive determinantal formulae for the
enumeration of walks {\it on the circle} (see the next section for
the precise definition), which includes the
enumeration of $n$ non-colliding particles on a circle. 

All of Grabiner's formulae are {\it exact\/} results. Hence, as an
afterthought, he posed the problem of determining the asymptotic
behaviour of the number of these walks if the number of steps becomes
large. It happens that this had already been done independently in
\cite{KrGVAB} for walks with diagonal steps in 
the alcove of type $\tilde C_n$, albeit in a completely
different language, the language of vicious walkers.

The purpose of this paper is to carry out the asymptotic analysis of
the number of random walks in alcoves of affine Weyl groups in all
the other cases, and also for the number of random walks on the
circle. To be precise, we determine the asymptotic behaviour of the
number of random walks in an alcove as the number of steps tends to
infinity for the case that starting and end point are held fixed, as
well as for the case where the end point can be arbitrary.
Frequently, the results depend
heavily on the parities of the involved parameters, a phenomenon which
distinguishes the discrete case from the continuous case. (This
phenomenon does also not occur for the corresponding problem for
the walks in chambers of finite Weyl groups; see the last paragraph of
the Introduction.)
While, from an analytic point of view, the order of magnitude
is always rather straight-forward to determine, for which very basic
tools (if at all), such as Stirling's formula, or, in one case, a
rather standard application of the saddle point method, suffice,
the determination of the
multiplicative constant poses quite frequently a substantial
challenge. Interestingly, carrying out the latter task requires quite often
some advanced facts from symmetric function theory (see 
the proofs in Sections~\ref{sec:4}--\ref{sec:7}).
In particular, identities for classical group characters
from \cite{KratBC} come in handy at many places.
It should be observed that the proofs show that the errors are always
exponentially small, with the exception of Theorem~\ref{thm:2-asy},
where the error is dictated by the Stirling approximation of the
binomial coefficient in \eqref{eq:3-1}, and of Theorem~\ref{thm:3-asy},
where the error is dominated by those coming 
from the saddle point approximation given in
Lemma~\ref{lem:1} in Appendix~\ref{app:A}. 

In the next section we provide the basic definitions, in particular,
the definitions of the alcoves to which our walks are confined, and
we summarize all the exact results that exist for the enumeration 
of the three types of walks in these alcoves. These will be the
starting points for our asymptotic calculations, which we carry out
in the subsequent sections. The results for the alcove of type
$\tilde A_{n-1}$ are given in Section~\ref{sec:3}, the results for the
enumeration of walks on the circle are the subject of Section~\ref{sec:4},
we give the results for the alcove of type
$\tilde C_{n}$ in Section~\ref{sec:5}, 
in Section~\ref{sec:6} there follow the results for the alcove of type
$\tilde B_{n}$,
and, finally, we present the results for the alcove of type
$\tilde D_{n}$ in Section~\ref{sec:7}. Auxiliary results that are
needed in the proofs of the theorems are collected in three appendices.

In concluding the introduction, it is probably useful to
review the state of affairs for {\it finite} Weyl groups, i.e., 
the known results on the
asymptotic behaviour of walks in chambers of finite Weyl
groups as the number of steps of the walks becomes large. 
(We refer the reader again to the book \cite{HumpAC} for
definitions and more information on finite Weyl groups.) 
In fact, since, as we already indicated, 
the random walk problems considered in this paper can
be seen from various different angles, numerous results can be found
scattered in the combinatorics, probability, physics, and even
representation theory literature. 
If starting and end point are fixed, the asymptotics 
of walks with diagonal steps in Weyl chambers of types
$A_{n-1}$ and $C_n$ were determined (in the language of vicious walkers) by
Rubey \cite[Ch.~2, Sections~3 and 4]{RubeADD}, \cite[Sections~3 and 4]{RubeAD}, 
with previous results in special cases given in
\cite[Sections~2, 4, 7]{KrGVAA}. For the case that the starting point is the origin
and the end point is fixed,
a result of Biane \cite{BianAE} on the asymptotics of multiplicities
of irreducible representations in tensor powers of irreducible
representations of semisimple Lie groups, combined with an
observation due to Grabiner and Magyar \cite[Sec.~3.3]{GrMrAA} that, under
mild restrictions, the
number of random walks in Weyl chambers is equal
to such multiplicities, implies a uniform asymptotic formula for 
random walks in Weyl chambers of any type, with the exception of
walks with standard steps in a Weyl chamber of type $A_{n-1}$.
For the case that the starting point is arbitrary but fixed and the
end point is not fixed,
Grabiner \cite{GrabAF} has recently shown that, by combining a
result of Kuperberg \cite[Theorem~1.2.1]{KupeAI} on the approximation
of sums of random variables defined on lattices by the corresponding
Brownian motion, and of himself \cite{GrabAB} on Brownian motion in
Weyl chambers, one obtains the dominating term of the asymptotic
behaviour for all types and for all possible step sets at once.
We want to remark that
walks with standard steps in the positive direction in a chamber of
the Weyl group of type $A_{n-1}$ are equivalent to skew standard
Young tableaux with at most $n$ rows, 
the shape depending on starting and end point of the
walk. The asymptotic behaviour of the number of non-skew standard
Young tableaux with at most $n$ rows is covered by the celebrated 
earlier (and more general) result of Regev \cite{RegeAG}. 
More precise results than that of Grabiner's (i.e., with bounds on the
errors also) in
the case of walks with diagonal steps in a Weyl chamber of type
$A_{n-1}$ were found (again in the language of vicious walkers) by
Rubey \cite[Ch.~2, Sec.~3]{RubeADD}, \cite[Sec.~3]{RubeAD}, 
with previous results in special cases given in
\cite[Sections~2, 3, 6]{KrGVAA}. In addition, Rubey 
\cite[Ch.~2, Sec.~4]{RubeADD}, \cite[Sec.~4]{RubeAD} 
also provides more precise results in 
the case of walks with diagonal steps in a Weyl chamber
of type $C_n$. (Again, results in special cases can already be
found in \cite[Sec.~4, 5, 7]{KrGVAA}.)

\section{A summary of exact results of random walks in alcoves
of affine Weyl groups} \label{sec:2}

In this section we summarize the exact results for random walks in 
alcoves of affine
Weyl groups, also including two results 
for random walks on a circle, which
are the starting points for our asymptotic calculations which 
are to follow in the later sections. We also use the opportunity to
point out, in each case, equivalent formulations of the walk
problems (in case they exist).

Before we state the results, let us recall the definitions of these
alcoves. Let $m$ be some given positive integer or
half-integer. (By definition, 
a {\it half-integer} is an odd number divided by 2.)
We define the alcove $\mathcal A^{\tilde A_{n-1}}_m$ of type $\tilde A_{n-1}$ to be
the region
\begin{equation} \label{eq:alcA}
\mathcal A^{\tilde A_{n-1}}_m:=\{(x_1,x_2,\dots,x_n):x_1>x_2>\dots>x_n>x_1-m\}.
\end{equation}
(Strictly speaking, this is a {\it scaled\/} alcove.)
The (scaled) alcove of type $\tilde C_n$ is defined by
\begin{equation} \label{eq:alcC}
\mathcal A^{\tilde C_{n}}_m:=\{(x_1,x_2,\dots,x_n):m>x_1>x_2>\dots>x_n>0\}.
\end{equation}
The (scaled) alcove of type $\tilde B_n$ is defined by
\begin{equation} \label{eq:alcB}
\mathcal A^{\tilde B_{n}}_m:=\{(x_1,x_2,\dots,x_n):x_1>x_2>\dots>x_n>0 \text { and }
x_1+x_2<2m\}.
\end{equation}
Finally, the (scaled) alcove of type $\tilde D_n$ is the region
\begin{equation} \label{eq:alcD}
\mathcal A^{\tilde D_{n}}_m:=\{(x_1,x_2,\dots,x_n):x_1>x_2>\dots>x_{n-1}>\vert
x_n\vert ,\text { and }
x_1+x_2<2m\}.
\end{equation}

\medskip
We begin with results for the enumeration of walks in $\mathcal
A^{\tilde A_{n-1}}_m$. The first result is originally due to Filaseta \cite{FilaAA}.
It is however covered by the general result \cite{GeZeAA}.
In the statement of the theorem, 
and also subsequently, given a vector
$\et=(\et_1,\et_2,\dots,\et_n)$ we use the symbol
$\vert\et\vert$ to denote the sum of its components, i.e.,
$\vert\et\vert:=\et_1+\et_2+\dots+\et_n$.

\begin{Theorem}[\cite{FilaAA}] \label{thm:1}
Let $m$ be a positive integer. Furthermore,
let $\et=(\et_1,\et_2,\dots,\et_n)$ and 
$\la=(\la_1,\la_2,\dots,\la_n)$ be vectors of integers in the
alcove $\mathcal A^{\tilde A_{n-1}}_m$ of type $\tilde A_{n-1}$ {\em(}defined in \eqref{eq:alcA}{\em)}.
Then the number of random walks from $\et$ to $\la$, which consist
entirely of standard steps in the positive direction, and 
which stay in the alcove $\mathcal A^{\tilde A_{n-1}}_m$, is given by
\begin{equation} \label{eq:Ae}
\big(\vert\la\vert-\vert\et\vert\big)!
\sum _{k_1+\dots+k_n=0} ^{}\det_{1\le h,t\le n}\(\frac {1}
{(\la_t-\et_h+mk_h)!}\).
\end{equation}
\end{Theorem}

The corresponding result for positive {\it and\/} negative standard
steps is also a direct consequence of the general result
\cite{GeZeAA}, and is stated explicitly in \cite{GrabAD}.

\begin{Theorem}[{\cite[Eq.~(34)]{GrabAD}}] \label{thm:2}
Let $m$ be a positive integer. Furthermore,
let $\et=(\et_1,\et_2,\dots,\et_n)$ and 
$\la=(\la_1,\la_2,\dots,\la_n)$ be vectors of integers in the
alcove $\mathcal A^{\tilde A_{n-1}}_m$ of type $\tilde A_{n-1}$ {\em(}defined in \eqref{eq:alcA}{\em)}.
Then the number of random walks from $\et$ to $\la$ with exactly $k$
standard steps,
which stay in the alcove $\mathcal A^{\tilde A_{n-1}}_m$, is given by
the coefficient of $x^k/k!$ in
\begin{equation} \label{eq:A+-e}
\sum _{k_1+\dots+k_n=0} ^{}\det_{1\le h,t\le
n}\(I_{\la_t-\et_h+mk_h}(2x)\),
\end{equation}
where $I_\al(x)$ is the modified Bessel function of the first kind
$$I_\al(x)=\sum _{j=0} ^{\infty}\frac {(x/2)^{2j+\al}}
{j!\,(j+\al)!}.$$
\end{Theorem}

This result has also a different interpretation: by
considering each of the $n$ coordinates as a separate walk, the walks
with standard steps in $\mathcal A^{\tilde A_{n-1}}_m$ 
can be seen as $n$ separate particles on the integer line, 
where at each tick of the clock {\it exactly one} particle moves to
the right or to the left by one unit (a move to the right, respectively to
the left, of the $j$-th particle corresponding
to a step $+e_j$, respectively $-e_j$), 
under the constraint that at no time two
particles occupy the same lattice site, and such that in addition a
shift by $m$ of any of the particles never collides with any of the
other particles. Thus we obtain a sub-model of Fisher's \cite{FIS}
{\it random turns vicious walker model}.\footnote{In statistical physics, 
a model of $n$ walkers on the integer line
where at each tick of the clock exactly one walker moves to
the right or to the left, under the constraint that at no time two
walkers occupy the same lattice site, is called the
{\it random turns vicious walker model}.}

Similarly,
random walks with diagonal steps in $\mathcal A^{\tilde A_{n-1}}_m$ with given
starting and end point can be seen in several ways: by
considering each of the $n$ coordinates as a separate walk, such
random walks can be seen as $n$ separate particles on the integer line, 
where at each tick of the clock {\it each} particle moves 
one unit step to the right
(corresponding to a change of $+\frac {1} {2}$ in the corresponding
coordinate) or to
the left (corresponding to a change of $-\frac {1} {2}$ in the corresponding
coordinate), such that they never collide, and such that in addition a
shift by $2m$ of any of the particles never collides with any of the
other particles. Thus we obtain a sub-model of Fisher's \cite{FIS}
{\it lock step vicious walker model}.\footnote{In statistical physics, 
a model of $n$ walkers on the integer line
where at each tick of the clock each walker moves to
the right or to the left, under the constraint that at no time two
walkers occupy the same lattice site, is called the
{\it lock step vicious walker model}.}
An alternative, two-dimensional picture arises if we
convert the movements in each coordinate of
the random walk to a separate path in the plane integer lattice, 
identifying 
a change by $+\frac {1} {2}$ in a coordinate with an up-step $(1,1)$ and 
a change by $-\frac {1} {2}$ with a down-step $(1,-1)$ of the corresponding
path. Thus, such random walks can be seen to be equivalent to families
of {\it non-intersecting}\footnote{A family of paths is called {\it
non-intersecting} if no two paths from the family have any common
points.} {\it lattice paths} in the plane integer lattice with
steps $(1,1)$ and $(1,-1)$
(the starting points of which being
aligned along a vertical line, as well as the end points) 
where a shift of the bottom-most
path dominates the top-most path. The latter objects are in turn in
bijection with (special) {\it cylindric partitions} 
(as defined in \cite{GeKrAA}) of rectangular shape (see \cite[Sec.~3]{GeKrAA}
for that translation; to obtain the presentation of the lattice paths
in \cite{GeKrAA}, the above described picture has to be rotated by 
$45^\circ$). 

The following result is at the same time a
direct consequence of the general result in \cite{GeZeAA} and of
Theorem~3 in \cite{GeKrAA}. It is stated explicitly in
\cite{GrabAD}. It is however important to note that it is only true
for {\it integral\/} $m$ (as well as the ``$m$-circle result"
Theorem~\ref{thm:4} for diagonal steps which it implies, 
as opposed to the companion results Theorems~\ref{thm:7}, \ref{thm:9}, and
\ref{thm:11} for the types
$\tilde C_n$, $\tilde B_n$, and $\tilde D_n$). This is because the reflection
argument from \cite{GeZeAA} (repeated in \cite{GrabAD}, and in an
equivalent form in \cite{GeKrAA}) only guarantees that (using the
vicious walker picture) particles never occupy the same site,
respectively a particle shifted by $2m$ never occupies the same site as
another particle. If $m$ is a half-integer, this does not exclude that
a shifted particle changes sides with another particle, and thus the
formula \eqref{eq:Ad} below would also include walks which violate the
condition $x_n>x_1-m$ which is contained in the definition \eqref{eq:alcA} of
the alcove $\mathcal A^{\tilde A_{n-1}}_m$.

\begin{Theorem}[{\cite[Eq.~(35)]{GrabAD}}] \label{thm:3}
Let $m$ be a positive integer. Furthermore,
let $\et=(\et_1,\et_2,\dots,\et_n)$ and 
$\la=(\la_1,\la_2,\dots,\la_n)$ be vectors of integers or of
half-integers in the
alcove $\mathcal A^{\tilde A_{n-1}}_m$ of type $\tilde A_{n-1}$ {\em(}defined in \eqref{eq:alcA}{\em)}.
Then the number of random walks 
from $\et$ to $\la$ with exactly $k$ diagonal steps,
which stay in the alcove $\mathcal A^{\tilde A_{n-1}}_m$, is given by
\begin{equation} \label{eq:Ad}
\sum _{k_1+\dots+k_n=0} ^{}\det_{1\le h,t\le
n}\(\binom k {\frac {k} {2}+\la_t-\et_h+mk_h}\).
\end{equation}
\end{Theorem}

As observed by Grabiner in \cite{GrabAD}, the above results can be
used to derive results on the enumeration of random walks on the
$m$-circle, where by ``{\it random walks on the $m$-circle}" we mean random
walks in $n$-dimensional Euclidean space, where each coordinate is
reduced modulo $m$ (i.e., a point $(x_1,x_2,\dots,x_n)$ is identified
with $(x_1+k_1m,x_2+k_2m,\dots,x_n+k_nm)$ for any integers
$k_1,k_2,\dots,k_n$).

Whereas in the case of standard steps in the positive direction
this does not define
a different model, it does for standard steps in the
positive {\it and\/}
negative direction, and also for diagonal steps. The result from
\cite{GrabAD} for standard steps reads as follows.

\begin{Theorem}[{\cite[Eq.~(32)]{GrabAD}}] \label{thm:5}
Let $m$ be a positive integer. Furthermore,
let $\et=(\et_1,\et_2,\dots,\et_n)$ be a vector of integers 
with $m>\et_1>\et_2>\dots>\et_n\ge0$, 
and let $\la=(\la_1,\la_2,\dots,\la_n)$ be a vector of integers
with $m>\la_{s+1}>\dots>\la_n>\la_1>\dots>\la_s\ge0$, for some $s$. 
Then the number of random walks on the $m$-circle
from $\et$ to $\la$ with exactly $k$ standard steps,
such that at no time two coordinates of a point on the random
walk are equal, is given by
\begin{equation} \label{eq:circ-e}
\frac {1} {n}\sum _{u=0} ^{n-1}e^{-2\pi ius/n}
\det_{1\le h,t\le n}\bigg(\frac {1} {m}\sum _{r=0} ^{m-1}
e^{-2\pi i(u+nr)(\la_t-\eta_h)/mn}\exp(2x\cos(2\pi(u+nr)/mn))\bigg).
\end{equation}
\end{Theorem}

In the same way as explained above for random walks in $\mathcal
A^{\tilde A_{n-1}}_m$, this result can also seen as counting $n$ non-colliding 
particles moving on
the integer circle of length $m$ (the interval $[0,m]$ with $0$ and
$m$ identified), where at each tick of the clock {\it exactly one} 
particle moves to
the right or to the left by one unit.

Similarly, 
random walks with diagonal steps on the circle, with the property 
that at no time
two coordinates of a point on the walk are equal can be equivalently
seen as the movements of $n$ non-colliding particles on a circle,
where at each tick of the clock each particle moves one unit step to the right
or to the left. This version of the lock step vicious walker model had
been first considered by Forrester \cite{ForrAB}. 
He solved the problem of counting the number of ways $n$ such
particles in this model may move from given starting points to given
end points in the case that $n$ is odd, however, where the particles
may reach the end points in {\it any} (cyclic) order (see  
\cite[Sec.~2.2]{ForrAB}). An
analogous formula for the case that $n$ is even has been recently
found by Fulmek \cite{FulmAF}.
Thus, the result from Grabiner's
paper \cite{GrabAD}, which we state 
below, constitutes a refinement
of Forrester's and Fulmek's formulae, as in Grabiner's formula the
order in which the particles arrive at the end points is fixed. 
In the statement below, a small typo from \cite{GrabAD} has been
corrected (in the determinant in Eq.~(33)
in \cite{GrabAD} the term $\ze^{-(u+nr)(\la_j-\et_i)}$ has to be
replaced by $\ze^{-2(u+nr)(\la_j-\et_i)}$).

\begin{Theorem}[{\cite[Eq.~(33)]{GrabAD}}] \label{thm:4}
Let $m$ be a positive integer or half-integer. Furthermore,
let $\et=(\et_1,\et_2,\dots,\et_n)$ be a vector of integers or of
half-integers with $m>\et_1>\et_2>\dots>\et_n\ge0$, 
and let $\la=(\la_1,\la_2,\dots,\la_n)$ be a vector of integers or of
half-integers
with $m>\la_{s+1}>\dots>\la_n>\la_1>\dots>\la_s\ge0$, for some $s$. 
Then the number of random walks on the $m$-circle
from $\et$ to $\la$ with exactly $k$ diagonal
steps, such that at no time two coordinates of a point on the random
walk are equal, is given by
\begin{equation} \label{eq:circ-d}
\frac {1} {n}\sum _{u=0} ^{n-1}e^{-2\pi ius/n}
\det_{1\le h,t\le n}\bigg(\frac {2^{k-1}} {m}\sum _{r=0} ^{2m-1}
e^{-2\pi i(u+nr)(\la_t-\eta_h)/mn}\cos^k(\pi(u+nr)/mn)\bigg).
\end{equation}
\end{Theorem}

Next we quote the two results from \cite{GrabAD} on the enumeration
of random walks in alcoves of type $\tilde C_n$. In this case, 
there is no separate result for positive standard steps, 
since for such walks the condition $m>x_1$, which appears in the
definition \eqref{eq:alcC} of the alcove $\mathcal A^{\tilde C_{n}}_m$, is without
meaning, as well as the condition $x_n>0$, so that the problem of 
enumerating random walks with positive
standard steps between two
given points which stay in $\mathcal A^{\tilde C_{n}}_m$ is equivalent to counting
random walks with positive standard steps which stay in the {\it Weyl
chamber of type $A_{n-1}$}, the latter being defined by 
\begin{equation} \label{eq:WA} 
\{(x_1,x_2,\dots,x_n):x_1>x_2>\dots>x_n\}.
\end{equation}
As we mentioned in the Introduction, this problem has been dealt with
in \cite{GrabAF,KrGVAA,RubeADD,RubeAD}.
On the other
hand, random walks in $\mathcal A^{\tilde C_{n}}_m$ from $\et$ to $\la$ with  
standard steps in the positive {\it
and\/} negative direction are equivalent to oscillating tableaux
from $(\et_1-n,\et_2-(n-1),\dots,\et_n-1)$ to
$(\la_1-n,\la_{2}-(n-1),\dots,\la_n-1)$ with at most $n$ rows and at
most $m-n$ columns, as is
easily seen by identifying a step $+e_j$ with the augmentation of the
$i$-th row of the Ferrers diagram by a box, respectively identifying
a step $-e_j$ with the deletion of a box from the
$i$-th row of the Ferrers diagram. (The reader should
recall that an oscillating tableau
is a sequence of Ferrers diagrams where successive diagrams in the
sequence differ by exactly one box.) The corresponding result from
\cite{GrabAD} reads as follows.

\begin{Theorem}[{\cite[Eq.~(23)]{GrabAD}}] \label{thm:6}
Let $m$ be a positive integer. Furthermore,
let $\et=(\et_1,\et_2,\dots,\et_n)$ and 
$\la=(\la_1,\la_2,\dots,\la_n)$ be vectors of integers in the
alcove $\mathcal A^{\tilde C_{n}}_m$ of type $\tilde C_n$ {\em(}defined in \eqref{eq:alcC}{\em)}.
Then the number of random walks from $\et$ to $\la$ with exactly $k$
standard steps,
which stay in the alcove $\mathcal A^{\tilde C_{n}}_m$, is given by
the coefficient of $x^k/k!$ in
\begin{equation} \label{eq:C+-e}
\det_{1\le h,t\le n}\(\frac {1} {m}\sum _{r=0} ^{2m-1}
\sin\frac {\pi r\la_t} {m}\cdot\sin\frac {\pi r\et_h} {m}\cdot
\exp\(2x\cos\frac {\pi r} {m}\)\).
\end{equation}
\end{Theorem}

On the other hand, random walks with diagonal steps in $\mathcal A^{\tilde C_{n}}_m$
are equivalent (by means of the translation explained earlier for 
random walks in $\mathcal A^{\tilde A_{n-1}}_m$) to the movements of $n$
non-colliding particles in an interval,
where at each tick of the clock each particle moves 
one unit step to the right
or to the left (see also \cite[Sec.~5]{GrabAD}).
Equivalently, these may seen as families of non-intersecting lattice
paths consisting of up- and down-steps between
two horizontal boundaries (see \cite[Sec.~5]{GrabAD} and
\cite{KrGVAB}). The corresponding 
result from \cite{GrabAD} is the following.

\begin{Theorem}[{\cite[Eq.~(18)]{GrabAD}}] \label{thm:7}
Let $m$ be a positive integer or half-integer. Furthermore,
let $\et=(\et_1,\et_2,\dots,\et_n)$ and 
$\la=(\la_1,\la_2,\dots,\la_n)$ be vectors of integers or of
half-integers in the
alcove $\mathcal A^{\tilde C_{n}}_m$ of type $\tilde C_n$ {\em(}defined in \eqref{eq:alcC}{\em)}.
Then the number of random walks 
from $\et$ to $\la$ with exactly $k$ diagonal steps,
which stay in the alcove $\mathcal A^{\tilde C_{n}}_m$, is given by
\begin{equation} \label{eq:Cd}
\det_{1\le h,t\le n}\(\frac {2^{k-1}} {m}\sum _{r=0} ^{4m-1}
\sin\frac {\pi r\la_t} {m}\cdot\sin\frac {\pi r\et_h} {m}\cdot
\cos^k\frac {\pi r} {2m}\).
\end{equation}
\end{Theorem}

The next two results concern the enumeration of random walks in
alcoves of type $\tilde B_n$. While the first result, the result for standard
steps, is stated explicitly in
\cite{GrabAD}, the second, the result for diagonal steps, 
is not made explicit there, although it is
made clear how to derive it. We state it here for
the sake of completeness. Again,
there is no separate result for positive standard steps, 
since for such walks the condition $2m>x_1+x_2$, which appears in the
definition \eqref{eq:alcB} of the alcove $\mathcal A^{\tilde B_{n}}_m$, is without
meaning, as well as the condition $x_n>0$, so that the problem of 
enumerating random walks with positive
standard steps between two
given points which stay in $\mathcal A^{\tilde B_{n}}_m$ is again equivalent to counting
random walks with positive standard steps which stay in the Weyl
chamber of type $A_{n-1}$, the latter being defined by \eqref{eq:WA}.

\begin{Theorem}[{\cite[Eq.~(43)]{GrabAD}}] \label{thm:8}
Let $m$ be a positive integer or half-integer. Furthermore,
let $\et=(\et_1,\et_2,\dots,\et_n)$ and 
$\la=(\la_1,\la_2,\dots,\la_n)$ be vectors of integers in the
alcove $\mathcal A^{\tilde B_{n}}_m$ of type $\tilde B_n$ 
{\em(}defined in \eqref{eq:alcB}{\em)}.
Then the number of random walks from $\et$ to $\la$ with exactly $k$
standard steps,
which stay in the alcove $\mathcal A^{\tilde B_{n}}_m$, is given by
the coefficient of $x^k/k!$ in
\begin{multline} \label{eq:B+-e}
\frac {1} {2}\det_{1\le h,t\le n}\(\frac {1} {m}\sum _{r=0} ^{2m-1}
\sin\frac {\pi r\la_t} {m}\cdot\sin\frac {\pi r\et_h} {m}\cdot
\exp\(2x\cos\frac {\pi r} {m}\)\)
\\+
\frac {1} {2}\det_{1\le h,t\le n}\(\frac {1} {m}\sum _{r=0} ^{2m-1}
\sin\frac {\pi (2r+1)\la_t} {2m}\cdot\sin\frac {\pi (2r+1)\et_h} {2m}\cdot
\exp\(2x\cos\frac {\pi (2r+1)} {2m}\)\).
\end{multline}
\end{Theorem}

\begin{Theorem}[{\cite{GrabAD}}] \label{thm:9}
Let $m$ be a positive integer or half-integer. Furthermore,
let $\et=(\et_1,\et_2,\dots,\et_n)$ and 
$\la=(\la_1,\la_2,\dots,\la_n)$ be vectors of integers or of
half-integers in the
alcove $\mathcal A^{\tilde B_{n}}_m$ of type $\tilde B_n$ {\em(}defined in \eqref{eq:alcB}{\em)}.
Then the number of random walks 
from $\et$ to $\la$ with exactly $k$ diagonal steps,
which stay in the alcove $\mathcal A^{\tilde B_{n}}_m$, is given by
\begin{multline} \label{eq:Bd}
\frac {1} {2}\det_{1\le h,t\le n}\(\frac {2^{k-1}} {m}\sum _{r=0} ^{4m-1}
\sin\frac {\pi r\la_t} {m}\cdot\sin\frac {\pi r\et_h} {m}\cdot
\cos^k\frac {\pi r}
{2m}\)\\
+\frac {1} {2}\det_{1\le h,t\le n}\(\frac {2^{k-1}} {m}\sum _{r=0} ^{4m-1}
\sin\frac {\pi (2r+1)\la_t} {2m}\cdot\sin\frac {\pi (2r+1)\et_h} {2m}\cdot
\cos^k\frac {\pi (2r+1)}
{4m}\).
\end{multline}
\end{Theorem}

The final two results concern the enumeration of random walks in
alcoves of type $\tilde D_n$. Again, the first result, the result for
standard steps, is stated explicitly in
\cite{GrabAD}, while the second, the result for diagonal steps, 
is not made explicit there, although,
again, it is clearly described how to derive it. 
We state it here for the sake of completeness. Also here,
there is no separate result for positive standard steps, 
since for such walks the condition $2m>x_1+x_2$, which appears in the
definition \eqref{eq:alcD} of the alcove $\mathcal A^{\tilde D_{n}}_m$, is without
meaning, as well as the condition $x_{n-1}>\vert x_n\vert$, 
so that the problem of 
enumerating random walks with positive
standard steps between two
given points which stay in $\mathcal A^{\tilde D_{n}}_m$ is again equivalent to counting
random walks with positive standard steps which stay in the Weyl
chamber of type $A_{n-1}$, the latter being defined by \eqref{eq:WA}.

\begin{Theorem}[{\cite[Eq.~(46)]{GrabAD}}] \label{thm:10}
Let $m$ be a positive integer or half-integer. Furthermore,
let $\et=(\et_1,\et_2,\dots,\et_n)$ and 
$\la=(\la_1,\la_2,\dots,\la_n)$ be vectors of integers in the
alcove $\mathcal A^{\tilde D_{n}}_m$ of type $\tilde D_n$ 
{\em(}defined in \eqref{eq:alcD}{\em)}.
Then the number of random walks from $\et$ to $\la$ with exactly $k$
standard steps,
which stay in the alcove $\mathcal A^{\tilde D_{n}}_m$, is given by
the coefficient of $x^k/k!$ in
\begin{multline} \label{eq:D+-e}
\frac {1} {4}\det_{1\le h,t\le n}\(\frac {1} {m}\sum _{r=0} ^{2m-1}
\sin\frac {\pi r\la_t} {m}\cdot\sin\frac {\pi r\et_h} {m}\cdot
\exp\(2x\cos\frac {\pi r} {m}\)\)
\\+
\frac {1} {4}\det_{1\le h,t\le n}\(\frac {1} {m}\sum _{r=0} ^{2m-1}
\sin\frac {\pi (2r+1)\la_t} {2m}\cdot\sin\frac {\pi (2r+1)\et_h} {2m}\cdot
\exp\(2x\cos\frac {\pi (2r+1)} {2m}\)\)
\\+
\frac {1} {4}\det_{1\le h,t\le n}\(\frac {1} {m}\sum _{r=0} ^{2m-1}
\cos\frac {\pi r\la_t} {m}\cdot\cos\frac {\pi r\et_h} {m}\cdot
\exp\(2x\cos\frac {\pi r} {m}\)\)
\\+
\frac {1} {4}\det_{1\le h,t\le n}\(\frac {1} {m}\sum _{r=0} ^{2m-1}
\cos\frac {\pi (2r+1)\la_t} {2m}\cdot\cos\frac {\pi (2r+1)\et_h} {2m}\cdot
\exp\(2x\cos\frac {\pi (2r+1)} {2m}\)\).
\end{multline}
\end{Theorem}

\begin{Theorem}[{\cite{GrabAD}}] \label{thm:11}
Let $m$ be a positive integer or half-integer. Furthermore,
let $\et=(\et_1,\et_2,\dots,\et_n)$ and 
$\la=(\la_1,\la_2,\dots,\la_n)$ be vectors of integers or of
half-integers in the
alcove $\mathcal A^{\tilde D_{n}}_m$ of type $\tilde D_n$ {\em(}defined in \eqref{eq:alcD}{\em)}.
Then the number of random walks 
from $\et$ to $\la$ with exactly $k$ diagonal steps,
which stay in the alcove $\mathcal A^{\tilde D_{n}}_m$, is given by
\begin{multline} \label{eq:Dd}
\frac {1} {4}\det_{1\le h,t\le n}\(\frac {2^{k-1}} {m}\sum _{r=0} ^{4m-1}
\sin\frac {\pi r\la_t} {m}\cdot\sin\frac {\pi r\et_h} {m}\cdot
\cos^k\frac {\pi r}
{2m}\)\\
+\frac {1} {4}\det_{1\le h,t\le n}\(\frac {2^{k-1}} {m}\sum _{r=0} ^{4m-1}
\sin\frac {\pi (2r+1)\la_t} {2m}\cdot\sin\frac {\pi (2r+1)\et_h} {2m}\cdot
\cos^k\frac {\pi (2r+1)}
{4m}\)\\
+\frac {1} {4}\det_{1\le h,t\le n}\(\frac {2^{k-1}} {m}\sum _{r=0} ^{4m-1}
\cos\frac {\pi r\la_t} {m}\cdot\cos\frac {\pi r\et_h} {m}\cdot
\cos^k\frac {\pi r}
{2m}\)\\
+\frac {1} {4}\det_{1\le h,t\le n}\(\frac {2^{k-1}} {m}\sum _{r=0} ^{4m-1}
\cos\frac {\pi (2r+1)\la_t} {2m}\cdot\cos\frac {\pi (2r+1)\et_h} {2m}\cdot
\cos^k\frac {\pi (2r+1)}
{4m}\).
\end{multline}
\end{Theorem}

\section{Asymptotics for random walks in alcoves of type $\tilde A$}
\label{sec:3}

This section is devoted to finding the asymptotic behaviour of the
number of walks from a given starting point to a given end point
which stay in the alcove $\mathcal A^{\tilde A_{n-1}}_m$ of type $\tilde A_{n-1}$ as the
number of steps becomes large, as well as the asymptotic behaviour
of the number of those walks which start at a given point but may
terminate anywhere. In technical terms, we determine the asymptotic
behaviour of the expressions given by
Theorems~\ref{thm:1}--\ref{thm:3} as $k$ becomes large (in the case of
Theorem~\ref{thm:1} the role of $k$ is played by
$\vert\la\vert-\vert\et\vert$), and as well if these expressions
are summed over all possible end points of the walks.

We begin with the walks with standard steps in the positive
direction. 

\begin{Theorem} \label{thm:1-asy}
Let $m$ be a positive integer. Furthermore,
let $\et=(\et_1,\et_2,\dots,\et_n)$ and 
$\la=(\la_1,\la_2,\dots,\la_n)$ be vectors of integers in the
alcove $\mathcal A^{\tilde A_{n-1}}_m$ of type $\tilde A_{n-1}$ {\em(}defined in \eqref{eq:alcA}{\em)}.
Then, for large $\vert\la\vert$, 
the number of random walks from $\et$ to $\la$, which consist
entirely of standard steps in the positive direction, and 
which stay in the alcove $\mathcal A^{\tilde A_{n-1}}_m$, is asymptotically
\begin{equation} \label{eq:Ae-asy}
\frac {2^{n^2-n}} {m^{n-1}}\(\frac {\sin \frac {n\pi}
{m}} {\sin\frac {\pi} {m}}\)^{\vert\la\vert-\vert\et\vert}
\prod _{1\le h<t\le n} ^{}\(\sin \frac
{\pi(\et_h-\et_t)} {m}\cdot\sin \frac
{\pi(\la_h-\la_t)} {m}\).
\end{equation}
\end{Theorem}
\begin{proof}
We want to estimate the expression \eqref{eq:Ae} for $\vert\la\vert$ large. 
In order to accomplish that, we
write it in the form
\begin{multline*} 
\coef{z^0}
\big(\vert\la\vert-\vert\et\vert\big)!
\sum _{k_1,\dots,k_n=-\infty} ^{\infty} z^{m
\sum _{j=1} ^{n}k_j}\det_{1\le h,t\le n}\(\frac {1}
{(\la_t-\et_h+mk_h)!}\)\\
=\coef{z^0}
\big(\vert\la\vert-\vert\et\vert\big)!
\det_{1\le h,t\le n}\(
\sum _{k_h=-\infty} ^{\infty} \frac {z^{m k_h}}
{(\la_t-\et_h+mk_h)!}\),
\end{multline*}
where here, and in the sequel, the notation $\langle f\rangle F$
stands for the coefficient of $f$ in (an appropriate expansion of) $F$.
(I.e., here, $\coef{z^0}g(z)$ denotes the constant coefficient in the
Laurent series $g(z)$.) 
With $\om=e^{2\pi i/m}$, we can rewrite this expression as
$$
\coef{z^0}
\big(\vert\la\vert-\vert\et\vert\big)!
\det_{1\le h,t\le n}\(
\frac {1} {m}\sum _{r_h=0} ^{m-1}
\sum _{k_h=-\infty} ^{\infty} \frac {\(z\om^{r_h}\)^{ k_h}}
{(\la_t-\et_h+k_h)!}\),
$$
because $\sum _{r=0} ^{m-1}\om^{r k}$ is equal to $m$ if $k$ is divisible by
$m$, and it vanishes otherwise.
Evaluating the sum over $k_h$, we obtain the expression
\begin{align} \notag
\coef{z^0}&
\big(\vert\la\vert-\vert\et\vert\big)!
\det_{1\le h,t\le n}\(
\frac {1} {m}\sum _{r_h=0} ^{m-1}
  {\(z\om^{r_h}\)^{\et_h-\la_t}}\exp\({z\om^{r_h}}\)
\)\\
\notag
&=\coef{z^{\vert\la\vert-\vert\et\vert}}
\frac {\big(\vert\la\vert-\vert\et\vert\big)!} {m^n}
\sum _{r_1,\dots,r_n=0} ^{m-1}
\exp\({z
\sum _{j=1} ^{n}\om^{r_j}}\)
\det_{1\le h,t\le n}\(
  \om^{r_h(\et_h-\la_t)}
\)\\
&=
\frac {1} {m^n}
\sum _{r_1,\dots,r_n=0} ^{m-1}
W(\mathbf r)^{\vert\la\vert-\vert\et\vert}
\det_{1\le h,t\le n}\(
  \om^{r_h(\et_h-\la_t)}
\),\label{eq:3-4}
\end{align}
where in the last line $W(\mathbf r)$ is an abbreviation for $\sum _{j=1}
^{n}\om^{r_j}$.

The sum in \eqref{eq:3-4} is a finite sum of the form
$\sum _{\ell} ^{}c_\ell
b_\ell^k$, with $k=\vert\la\vert-\vert\et\vert$, 
where the $b_\ell$'s are of the form $W(\mathbf r)$, and the
$c_\ell$'s are bounded.
Thus, the asymptotic behaviour of this sum as
$k\to\infty$ is dominated by the terms $c_\ell b_\ell^k$ for which $\vert
b_\ell\vert$ is maximal (and $c_\ell\ne0$ of course). 

Now, if two summation indices $r_h$ and $r_t$, for $h\ne t$, should
be equal, then the determinant in the summand in \eqref{eq:3-4}
vanishes. Therefore we may restrict the sum in \eqref{eq:3-4} to
indices $r_1,r_2,\dots,r_n$ which are pairwise distinct. 
Among the latter, the sets of indices $\{r_1,r_2,\dots,r_n\}$
for which $W(\mathbf r)$ has largest modulus are those for which the
$r_j$'s are as ``close" together as possible, i.e., the sets
\begin{equation} \label{eq:3-2}
\{r_1,r_2,\dots,r_n\}=\{\ell,\ell+1,\dots,\ell+n-1\},
\end{equation}
for some $\ell$ between $0$ and $m-1$.
(On the right-hand side, the elements must be reduced modulo $m$.)
Hence, let $r_j=\ell+\si(j)-1$, $j=1,2,\dots,n$, for some permutation
$\si\in S_n$, with $S_n$ denoting the symmetric group of order $n$. 
For this choice of indices, we have
\begin{align*} \notag
W(\mathbf r)^{\vert\la\vert-\vert\et\vert}&
\det_{1\le h,t\le n}\(
  \om^{r_h(\et_h-\la_t)}\)\\
&=
  \left\vert
\sum _{j=1} ^{n}\om^{j-1}\right\vert^{\vert\la\vert-\vert\et\vert}
\om^{\frac {n-1} {2}
(\vert\la\vert-\vert\et\vert)}
\det_{1\le h,t\le n}\(\om^{(\si(h)-1)(\et_h-\la_t)}\)\\
&=\(\frac {\sin\frac {n\pi} {m}} {\sin\frac {\pi}
{m}}\)^{\vert\la\vert-\vert\et\vert}
\om^{\frac {n-1} {2}(\vert\la\vert-\vert\et\vert)}
\om^{\sum _{j=1} ^{n}(\si(j)-1)\et_j}\det_{1\le h,t\le
n}\(\om^{-(\si(h)-1)\la_t}\).
\end{align*}
Thus, if we combine all our findings, we obtain that, as
$(\vert\la\vert-\vert\et\vert)\to\infty$, 
the expression \eqref{eq:3-4} is asymptotically
\begin{multline*} 
\frac {1} {m^n}
\(\frac {\sin\frac {n\pi} {m}} {\sin\frac {\pi}
{m}}\)^{\vert\la\vert-\vert\et\vert}
\sum _{\ell=0} ^{m-1}
\sum _{\si\in S_n} ^{}
\om^{\frac {n-1} {2}(\vert\la\vert-\vert\et\vert)}
\om^{\sum _{j=1} ^{n}(\si(j)-1)\et_j}(\sgn\si)\det_{1\le h,t\le
n}\(\om^{-(h-1)\la_t}\)\\
=\frac {1} {m^{n-1}}
\(\frac {\sin\frac {n\pi} {m}} {\sin\frac {\pi} {m}}\)^
{\vert\la\vert-\vert\et\vert}
\om^{\frac {n-1} {2}({\vert\la\vert-\vert\et\vert})}
\det_{1\le h,t\le
n}\(\om^{(h-1)\et_t}\)\det_{1\le h,t\le
n}\(\om^{-(h-1)\la_t}\).
\end{multline*}
Both determinants are Vandermonde determinants, and are therefore
easily evaluated. The resulting expression is exactly
\eqref{eq:Ae-asy}.
\end{proof}

Next we address the question of determining the asymptotic behaviour
of the number of {\it all\/} walks which start in a given point and
proceed for $k$ standard steps in the positive direction, 
always staying in the alcove $\mathcal A^{\tilde A_{n-1}}_m$
of type $\tilde A_{n-1}$. Clearly, this amounts to summing the expression
\eqref{eq:Ae-asy} over all $\la$ with $\vert\la\vert=\vert\et\vert+k$.
A moment's reflection shows that this is in fact a finite sum, with
the number of terms bounded by $(2m)^m$ (to give a very crude bound), 
a quantity which is independent of
$k$. Thus, it is obvious that the order of magnitude of the number of
all these walks is $\({\sin \frac {n\pi}
{m}} /{\sin\frac {\pi} {m}}\)^k$. 
However, determining the multiplicative constant poses a formidable
challenge, in
particular if one attempts to do it directly from summing up the expression
\eqref{eq:Ae-asy}. An elegant way to bypass (some of this) difficulty
is to set up a relationship between the enumeration of walks with
positive standard steps and the enumeration of walks with arbitrary
standard steps, which we do in Lemma~\ref{lem:gleich}.

\begin{Theorem} \label{thm:1-asy2}
Let $m$ be a positive integer. Furthermore,
let $\et=(\et_1,\et_2,\dots,\et_n)$ be a vector of integers in the
alcove $\mathcal A^{\tilde A_{n-1}}_m$ of type $\tilde A_{n-1}$ {\em(}defined in \eqref{eq:alcA}{\em)}.
Then, as $k$ tends to infinity, 
the number of random walks which start at $\et$ and proceed for exactly
$k$ standard steps in the positive direction, 
which stay in the alcove $\mathcal A^{\tilde A_{n-1}}_m$, is asymptotically
\begin{multline} \label{eq:Ae-asy-even-even1}
\frac {2^{\binom n2}} {m^{n/2}}
\(\frac {\sin\frac {n\pi} {m}} {\sin\frac {\pi} {m}}\)^k
\prod
_{1\le h<t\le n} ^{}\(\sin\frac {\pi(\eta_h-\eta_t)} {m}\)\\
\times
\(
\prod _{h=1} ^{n/2}\cot\frac {(2h-1)\pi} {2m}
+(-1)^{\vert\et\vert+k+\frac {n} {2}}
\prod _{h=1} ^{n/2}\tan\frac {(2h-1)\pi} {2m}\)
\end{multline}
if both $n$ and $m$ are even,
it is asymptotically
\begin{equation} \label{eq:Ae-asy-odd-even}
\frac {2^{\binom n2}} {m^{n/2}}
\(\frac {\sin\frac {n\pi} {m}} {\sin\frac {\pi} {m}}\)^k
\prod
_{1\le h<t\le n} ^{}\(\sin\frac {\pi(\eta_h-\eta_t)} {m}\)
\prod _{h=1} ^{n/2}\cot\frac {(2h-1)\pi} {2m}
\end{equation}
if $n$ is even and $m$ is odd, and it is asymptotically
\begin{equation} \label{eq:Ae-asy-odd-odd}
\frac {2^{\binom n2}} {m^{(n-1)/2}}
\(\frac {\sin\frac {n\pi} {m}} {\sin\frac {\pi} {m}}\)^k
\prod
_{1\le h<t\le n} ^{}\(\sin\frac {\pi(\eta_h-\eta_t)} {m}\)
\prod _{h=1} ^{(n-1)/2}\cot\frac {h\pi} {m}
\end{equation}
if $n$ is odd (regardless of $m$).
\end{Theorem}

\begin{proof}
Our original proof proceeded as indicated in the paragraph before the
statement of the theorem, namely by summing the expression
\eqref{eq:Ae-asy} over all possible $\la$ with
$\vert\la\vert=\vert\et\vert+k$. Although feasible, this path turned out
to be a thorny one. After the result had been obtained, the surprising
observation was that the asymptotic behaviour of walks {\it in
$\mathcal A^{\tilde A_{n-1}}_m$ with
positive standard steps} is, up to a factor of $2^{k}$, identical with
the asymptotic behaviour of walks {\it on the $m$-circle with
arbitrary} (i.e., positive {\it and\/} negative)
{\it standard steps} (cf.\ Theorem~\ref{thm:5-asy2}). 
As we show in Lemma~\ref{lem:gleich} below,
this is even true {\it non-asymptotically}.
(This is indeed the assertion of Lemma~\ref{lem:gleich} since it does
not matter whether we are in the alcove $\mathcal A^{\tilde A_{n-1}}_m$
or on the $m$-circle if we let the end point of the walks be arbitrary.)
The theorem now follows by the (independent) proof of
Theorem~\ref{thm:5-asy2} given in Section~\ref{sec:4}.
\end{proof}

\begin{Lemma} \label{lem:gleich}
The number of random walks 
which start at $\et$, proceed for exactly
$k$ standard steps in the positive direction, 
and stay in the alcove $\mathcal A^{\tilde A_{n-1}}_m$, is equal to 
$2^{-k}$ times the number of random walks which
start at $\et$, proceed for exactly $k$ standard steps,
and stay in the alcove $\mathcal A^{\tilde A_{n-1}}_m$.
\end{Lemma}

\begin{proof}
We are going to show this inductively. Let
$\et=(\et_1,\et_2,\dots,\et_n)$ be given. We decompose $\et$ into its
maximal {\it circular} subsequences of consecutive elements. The
meaning of ``circular" is that $\et_n$ and $\et_1$ are considered to
be ``consecutive" if $\et_n=\et_1-m+1$. 
In the sequel, we shall use the short hand {\it
maximal circular subsequences} for these subsequences. For
example, if $m=11$ there are 3 maximal circular subsequences in
$(9,8,5,4,3,1)$, namely 
\begin{equation} \label{eq:gaps1} 
(9,8),\ (5,4,3),\ (1).
\end{equation}
On the other
hand, if $m=9$, then the maximal circular subsequences are
\begin{equation} \label{eq:gaps2} 
(1,9,8),\ (5,4,3).
\end{equation}
Let the maximal circular
subsequences in the decomposition of $\et$ have respectively
the lengths $a_1,a_2,\dots,a_\ell$, 
with gaps $b_1,b_2,\dots,b_\ell$. More
precisely, the
gap $b_j$ is the difference between the smallest element in the $j$-th
sub-sequence and the largest element in the $(j+1)$-st sub-sequence,
reduced modulo $m$. Here, $(j+1)$ has to be interpreted as 1 if
$j=\ell$, Thus,
in the example \eqref{eq:gaps1} we have $a_1=2$, $a_2=3$, $a_3=1$, and
$b_1=3$, $b_2=2$, $b_3=3$, while
in the example \eqref{eq:gaps2} we have $a_1=3$, $a_2=3$, and
$b_1=3$, $b_2=2$.

It is obvious that, starting from such an $\et$, there are exactly $\ell$ ways
to move by a positive standard step. To be precise, one would increase
the largest element of a maximal
circular sub-sequence by 1. Similarly,
there are exactly $2\ell$ ways to move by an {\it arbitrary} standard step,
namely the $\ell$ possibilities of 
{\it positive} standard steps described above, together with the
$\ell$ possibilities of decreasing a {\it smallest\/} element of a maximal
circular sub-sequence by 1. Thus, for one step,
i.e., for $k=1$, our claim is true.

The question is whether this persists. As we have seen, the number of
possibilities to walk is (for {\it positive} standard steps, as
well as for {\it arbitrary} standard steps) a multiple of the number of 
maximal circular subsequences. Thus, if we are
able to show that the total number of maximal
circular sub-sequences in the set of all
possible points that we reached from $\et$ by walking one positive
standard step is exactly one half of the corresponding number of
maximal circular sub-sequences in the set of
all possible points that we reached from $\et$ by walking one arbitrary
standard step, then we are sure that the ratio of $1:2$ will continue
to hold for each step. 

The latter claim is easy to establish: if we increase the largest
element of the $j$-th maximal circular sub-sequence of $\et$ by 1, then we
obtain a point whose decomposition has 
$$\ell+1-\chi(a_j=1)-\chi(b_{j-1}=2)$$
maximal circular subsequences, where 
we used the notation $\chi(\mathcal A)$=1 if $\mathcal A$ is
true and $\chi(\mathcal A)$=0 otherwise, and where $b_0$ has to be
interpreted as $b_\ell$. On the other hand,
if we decrease the smallest
element of the $j$-th maximal circular sub-sequence of $\et$ by 1, then we
obtain a point whose decomposition has 
$$\ell+1-\chi(a_j=1)-\chi(b_{j}=2)$$
maximal circular subsequences. Thus, the total number of maximal
circular subsequences in the set of points reached from $\et$ by walking one
{\it positive} standard step is
$$\ell(\ell+1)-
\sum _{j=1} ^{\ell}\big(\chi(a_j=1)+\chi(b_{j}=1)\big),$$
while the total number of maximal
circular subsequences in the set of points reached from $\et$ by walking one
{\it arbitrary} standard step is exactly twice of that.
This finishes the proof of the lemma.
\end{proof}

Now we turn our attention to walks with positive {\it and\/} negative
standard steps. The theorem below gives the asymptotic behaviour of
the walks in $\mathcal A^{\tilde A_{n-1}}_m$ with fixed starting and end point. If
the end point is allowed to be arbitrary, then the enumeration of the
corresponding walks is equivalent to the enumeration of walks with
standard steps on the $m$-circle, and, thus, its asymptotic behaviour
is given by Theorem~\ref{thm:5-asy2} in Section~\ref{sec:4}.

\begin{Theorem} \label{thm:2-asy}
Let $m$ be a positive integer. Furthermore,
let $\et=(\et_1,\et_2,\dots,\et_n)$ and 
$\la=(\la_1,\la_2,\dots,\la_n)$ be vectors of integers in the
alcove $\mathcal A^{\tilde A_{n-1}}_m$ of type $\tilde A_{n-1}$ {\em(}defined in \eqref{eq:alcA}{\em)}.
Then, as $k$ tends to infinity such that $k\equiv
\vert\et\vert+\vert\la\vert$ mod\/~$2$, 
the number of random walks from $\et$ to $\la$ with exactly $k$
standard steps,
which stay in the alcove $\mathcal A^{\tilde A_{n-1}}_m$, is asymptotically
\begin{equation} \label{eq:A+-e-asy}
\frac {2^{n^2-n}} {m^{n-1}}\sqrt{\frac {2} {\pi k}}\(2\frac {\sin \frac {n\pi}
{m}} {\sin\frac {\pi} {m}}\)^k\prod _{1\le h<t\le n} ^{}\(\sin \frac
{\pi(\et_h-\et_t)} {m}\cdot\sin \frac
{\pi(\la_h-\la_t)} {m}\).
\end{equation}
\end{Theorem}

\begin{proof} 
We have to determine the asymptotics of the coefficient of $x^k/k!$ in
\eqref{eq:A+-e} as
$k\to\infty$. To begin with, we rewrite this coefficient as
\begin{multline*} \label{}
\coef{z^0\,\frac {x^k} {k!}}\sum _{k_1,\dots,k_n=-\infty} ^{\infty}
z^{m(k_1+\dots+k_n)}\det_{1\le h,t\le
n}\(I_{\la_t-\et_h+mk_h}(2x)\)\\
=\coef{z^0\,\frac {x^k} {k!}}\det_{1\le h,t\le
n}\(\sum _{k_h=-\infty}
^{\infty}I_{\la_t-\et_h+mk_h}(2x)\,z^{mk_h}\),
\end{multline*}
Here, 
$\langle{z^0\frac {x^k} {k!}}\rangle g(x,z)$ denotes the coefficient of
$z^0 x^k/k!$ in $g(x,z)$, which is
in accordance with our earlier general definition of the
coefficient notation.
As in the proof of Theorem~\ref{thm:1-asy},
with $\om=e^{2\pi i/m}$ we may rewrite this expression as
$$\coef{z^0\,\frac {x^k} {k!}}\det_{1\le h,t\le
n}\(\frac {1} {m}\sum _{r_h=0} ^{m-1}
\sum _{k_h=-\infty} ^{\infty}I_{\la_t-\et_h+k_h}(2x)\(\om^{r_h} z\)^{k_h}\).
$$
Now using the easily verified fact that
$$\sum _{j=-\infty} ^{\infty}I_j(2x)\,z^j=\exp\Big(x\(z+z^{-1}\)\Big),$$
we obtain
\begin{align} \notag
&\coef{z^0\,\frac {x^k} {k!}}\det_{1\le h,t\le
n}\(\frac {1} {m}\(\om^{r_h}z\)^{\et_h-\la_t}\sum _{r_h=0} ^{m-1}
\exp\Big(x\(\om^{r_h}z+(\om^{r_h}z)^{-1}\)\Big)\)\\
\notag
&\quad 
=\coef{z^{\vert\la\vert-\vert\et\vert}\,\frac {x^k} {k!}}
\frac {1} {m^n}\sum _{r_1,\dots,r_n=0} ^{m-1}
\det_{1\le h,t\le
n}\(\om^{r_h(\et_h-\la_t)}\)\\
\notag
&\kern8cm
\cdot
\exp\(x\big({\textstyle z\sum _{j=1} ^{n}\om^{r_j}+z^{-1}\sum _{j=1}
^{n}\om^{-r_j}}\big)\)\\
\notag
&\quad 
=\coef{z^{\vert\la\vert-\vert\et\vert}}
\frac {1} {m^n}\sum _{r_1,\dots,r_n=0} ^{m-1}
\det_{1\le h,t\le
n}\(\om^{r_h(\et_h-\la_t)}\)
\cdot
\big({\textstyle z\sum _{j=1} ^{n}\om^{r_j}+z^{-1}\sum _{j=1}
^{n}\om^{-r_j}}\big)^k\\
\notag
&\quad 
=\frac {1} {m^n}
\binom k{\frac {1} {2}\big(k+\vert\la\vert-\vert\et\vert\big)}
\sum _{r_1,\dots,r_n=0} ^{m-1}
\det_{1\le h,t\le
n}\(\om^{r_h(\et_h-\la_t)}\)
\\
\notag
&\kern7cm
\cdot
\(\sum _{j=1} ^{n}\om^{r_j}\)^{\frac {1}
{2}\(k+\vert\la\vert-\vert\et\vert\)}
\(\sum _{j=1} ^{n}\om^{-r_j}\)^{\frac {1}
{2}\(k-\vert\la\vert+\vert\et\vert\)}
\end{align}
for the coefficient of $x^k/k!$ in \eqref{eq:A+-e}. 
Writing again $W(\mathbf r)$ for $\sum _{j=1}
^{n}\om^{r_j}$, we have obtained that the coefficient of $x^k/k!$ in
\eqref{eq:A+-e} is equal to
\begin{multline} \label{eq:3-1}
\frac {1} {m^n}
\binom k{\frac {1} {2}\big(k+\vert\la\vert-\vert\et\vert\big)}\\
\times
\sum _{r_1,\dots,r_n=0} ^{m-1}
\vert W(\mathbf r)\vert^k
\( {W(\mathbf r)} \Big/{\overline {W(\mathbf r)}}\)^{\frac {1} {2}
(\vert\la\vert-\vert\et\vert)}
\det_{1\le h,t\le
n}\(\om^{r_h(\et_h-\la_t)}\).
\end{multline}

Stirling's formula implies that 
the binomial coefficient in this expression is asymptotically
$2^k\sqrt{\frac 2{k\pi}}$ as $k\to\infty$.

The sum in \eqref{eq:3-1}, on the other hand, is again a finite sum of the form
$\sum _{\ell} ^{}c_\ell
b_\ell^k$, with the $b_\ell$'s of the form $\vert W(\mathbf r)\vert$,
and both the $c_\ell$'s and $b_\ell$'s independent of $k$.
Thus, the asymptotic behaviour of this sum as
$k\to\infty$ is dominated by the terms $c_\ell b_\ell^k$ for which $\vert
b_\ell\vert$ is maximal (and $c_\ell\ne0$ of course). 

If we compare the expression \eqref{eq:3-1} that we have obtained so far with
the expression
\eqref{eq:3-4}, then we see that we are in a very similar situation
here as
at the analogous place in the proof of Theorem~\ref{thm:1-asy}. 
Therefore, if we apply the arguments given there to our situation, we
obtain that, as
$k\to\infty$, the expression \eqref{eq:3-1} is asymptotically
\begin{multline*} 
\frac {1} {m^n}\sqrt{\frac {2} {\pi k}}
\(2\frac {\sin\frac {n\pi} {m}} {\sin\frac {\pi} {m}}\)^k
\sum _{\ell=0} ^{m-1}
\sum _{\si\in S_n} ^{}
\om^{\frac {n-1} {2}(\vert\la\vert-\vert\et\vert)}
\om^{\sum _{j=1} ^{n}(\si(j)-1)\et_j}(\sgn\si)\det_{1\le h,t\le
n}\(\om^{-(\si(h)-1)\la_t)}\)\\
=\frac {1} {m^{n-1}}\sqrt{\frac {2} {\pi k}}
\(2\frac {\sin\frac {n\pi} {m}} {\sin\frac {\pi} {m}}\)^k
\om^{\frac {n-1} {2}(\vert\la\vert-\vert\et\vert)}
\det_{1\le h,t\le
n}\(\om^{(h-1)\et_t}\)\det_{1\le h,t\le
n}\(\om^{-(h-1)\la_t}\).
\end{multline*}
Both determinants are Vandermonde determinants, and are therefore
easily evaluated. The resulting expression is exactly
\eqref{eq:A+-e-asy}.
\end{proof}

As the final issue in this section, we consider the asymptotic
behaviour of walks in the alcove $\mathcal A^{\tilde A_{n-1}}_m$ consisting of
diagonal steps. The theorem below provides the solution of the problem
if the starting and end point are fixed. Exceptionally, to resolve this
problem, we need to apply a more advanced asymptotic method, 
the saddle point method (although a rather basic instance of it). 
Should the end point be
allowed to be arbitrary, then the corresponding enumeration problem is
equivalent to the enumeration of walks with
diagonal steps on the $m$-circle, and, thus, its asymptotic behaviour
is given by Theorem~\ref{thm:4-asy2} in Section~\ref{sec:4}.

\begin{Theorem} \label{thm:3-asy}
Let $m$ be a positive integer or half-integer. Furthermore,
let $\et=(\et_1,\et_2,\dots,\et_n)$ and 
$\la=(\la_1,\la_2,\dots,\la_n)$ be vectors of integers or of
half-integers in the
alcove $\mathcal A^{\tilde A_{n-1}}_m$ of type $\tilde A_{n-1}$ {\em(}defined in \eqref{eq:alcA}{\em)}.
Then, as $k$ tends to infinity such that $k\equiv 2\et_j+2\la_j$ mod\/~$2$, 
the number of random walks 
from $\et$ to $\la$ with exactly $k$ diagonal steps,
which stay in the alcove $\mathcal A^{\tilde A_{n-1}}_m$, is asymptotically
\begin{equation} \label{eq:Ad-asy}
\frac {1} {m^{n-1}}
\frac {2^{n^2-n}} { \sqrt{2\pi c_0 k}}
 \(2^n\prod _{j=1} ^{n} \cos\frac {\pi\(j-\frac {n+1} {2}\)} {m}\)^k
\prod
_{1\le h<t\le n} ^{}\(\sin\frac {\pi(\eta_h-\eta_t)} {m}\cdot
\sin\frac {\pi(\la_h-\la_t)} {m}\),
\end{equation}
where 
$$c_0=\sum _{j=1} ^{n}\(2\cos\frac {\pi\(j-\frac {n+1} {2}\)} {m}\)^{-2}.$$
\end{Theorem}

\begin{proof}
We have to determine the asymptotics of the expression \eqref{eq:Ad} as
$k\to\infty$. To begin with, we write \eqref{eq:Ad} in the form
\begin{multline*} \label{}
\coef{z^0}\sum _{k_1,\dots,k_n=-\infty} ^{\infty}
z^{m(k_1+\dots+k_n)}\det_{1\le h,t\le
n}\(\binom k {\frac {k} {2}+\la_t-\et_h+mk_h}\)\\
=\coef{z^0}\det_{1\le h,t\le
n}\(\sum _{k_h=-\infty} ^{\infty}
\binom k {\frac {k} {2}+\la_t-\et_h+mk_h}z^{mk_h}\),
\end{multline*}
where, again, the notation $\langle z^0\rangle g(z)$
denotes the coefficient of $z^0$ in $g(z)$.
Arguing in the same way as before in the proofs of
Theorems~\ref{thm:1-asy} and \ref{thm:2-asy}, with $\om=e^{2\pi i/m}$
we can rewrite this expression as
\begin{align} \notag
\coef{z^0}&\det_{1\le h,t\le
n}\(\frac {1} {m}\sum _{r_h=0} ^{m-1}
\sum _{k_h=-\infty} ^{\infty}
\binom k {\frac {k} {2}+\la_t-\et_h+k_h}\(\om^{r_h}z\)^{k_h}\)\\
\notag
&=\coef{z^0}\frac {1} {m^n}\det_{1\le h,t\le
n}\(\sum _{r_h=0} ^{m-1}
\(\om^{r_h}z\)^{-\frac {k} {2}-\la_t+\et_h}
(1+\om^{r_h}z)^k\)\\
\label{eq:3-3}
&=
\frac {1} {m^n}\sum _{r_1,\dots,r_n=0} ^{m-1}
\om^{-\frac {k} {2}\vert\mathbf r\vert}
\det_{1\le h,t\le
n}\(\om^{r_h(\et_h-\la_t)}\)
\coef{z^{\frac {nk} {2}+\vert\la\vert-\vert\et\vert}}
\prod _{j=1} ^{n}(1+\om^{r_j}z)^k,
\end{align}
where in the next-to-last line we used the binomial theorem.

We have again obtained a finite sum. Therefore the task now is to isolate the
summands which are asymptotically largest as $k\to\infty$. First of
all, if two summation indices $r_h$ and $r_t$, for $h\ne t$, should
be equal, then the determinant in the summand in \eqref{eq:3-3}
vanishes. Therefore we may restrict the sum in \eqref{eq:3-3} to
the summands corresponding to 
indices $r_1,r_2,\dots,r_n$ which are pairwise distinct. 

According to Lemma~\ref{lem:1} in Appendix~\ref{app:A}, 
among the latter, those will be
asymptotically largest for which
\begin{equation} \label{eq:cosprod} 
\prod _{j=1} ^{n} \left\vert \cos\frac {\pi(\th_0+r_j)} {m}\right\vert
\end{equation}
is largest (where $\th_0$ is a solution of \eqref{eq:th0}). 
As is not difficult to see, these are those sets of
indices $\{r_1,r_2,\dots,r_n\}$ for which the
$r_j$'s are as ``close" together as possible, i.e., again the sets
as given in \eqref{eq:3-2}, 
for some $\ell$ between $0$ and $m-1$. 

Let $\ell$ be fixed. 
Let $\{r_1,r_2,\dots,r_n\}$ be a set of indices 
as in \eqref{eq:3-2}, i.e., $r_j=\ell+\si(j)-1$, $j=1,2,\dots,n$, 
for some permutation $\si\in S_n$. 
For this set, there is a unique $\th_0$ such
that \eqref{eq:cosprod} is maximal, namely 
$\th_0=-\ell-\frac {n-1} {2}$.
Thus, using Lemma~\ref{lem:1} in \eqref{eq:3-3}, we obtain that the
expression \eqref{eq:3-3} is asymptotically
\begin{multline*} \label{}
\frac {1} {m^n}
\sum _{\ell=0} ^{m-1}
\sum _{\si\in S_n} ^{}
\det_{1\le h,t\le n}\(\om^{(\si(h)-1)(\et_h-\la_t)}\)
\frac {\om^{\frac {n-1} {2}
(\vert\la\vert-\vert\et\vert)}} { \sqrt{2\pi c_0 k}}
\prod _{j=1} ^{n} \(2 \cos\frac {\pi\(j-\frac {n+1} {2}\)} {m}\)^k\\
=\frac {1} {m^{n-1}}
\frac {\om^{\frac {n-1} {2}
(\vert\la\vert-\vert\et\vert)}} { \sqrt{2\pi c_0 k}}
\prod _{j=1} ^{n} \(2 \cos\frac {\pi\(j-\frac {n+1} {2}\)} {m}\)^k\\
\times
\det_{1\le h,t\le
n}\(\om^{(h-1)\et_t}\)
\det_{1\le h,t\le
n}\(\om^{-(h-1)\la_t}\),
\end{multline*}
where $c_0$ is given as in the statement of the theorem.
Again, both determinants are Vandermonde determinants, and are therefore
easily evaluated. The resulting expression is exactly
\eqref{eq:Ad-asy}.
\end{proof}

\section{Asymptotics for random walks on the circle}
\label{sec:4}

In this section we find the asymptotic behaviour of the number of
walks from a given starting point to a given end point on the
$m$-circle as the
number of steps becomes large, as well as the asymptotic behaviour
of the number of those walks which start at a given point but may
terminate anywhere. In technical terms, we determine the asymptotic
behaviour of the expressions given by
Theorems~\ref{thm:5} and \ref{thm:4} as $k$ becomes large, 
and as well if these expressions
are summed over all possible end points of the walks.

Before we state the next theorem, which gives the asymptotic
behaviour of walks with standard steps between two fixed points on the
$m$-circle, we need to discuss under which conditions such walks can exist.
In the theorem below we consider walks from 
$\et=(\et_1,\et_2,\dots,\et_n)$ to 
$\la=(\la_1,\la_2,\dots,\la_n)$ where
$\la_{s+1}>\dots>\la_n>\la_1>\dots>\la_s$, which means that, when 
we interpret such walks as the movements of $n$ separate particles,
the first $s$ particles (the particles which start at
$\et_1,\et_2,\dots,\et_s$) wind themselves once more around the
circle than the other particles. Thus, if we want to get from $\et$
to $\la$ in $k$ steps, we must have
\begin{equation} \label{eq:Bed+-e}
k\equiv\vert\la\vert-\vert\et\vert+Nnm+sm
\pmod2,
\end{equation}
where the integer $N$ 
is the number of times the latter particles wind around the circle.
If $m$ is even, then this condition reduces to the familiar
$k\equiv\vert\la\vert+\vert\et\vert$ mod\/~$2$. However, if $m$ is odd,
then there are two possibilities. If in addition $n$ is even, then,
depending on whether $k$ is even or odd, there are walks only for
every second $s$. If both $m$ and $n$ are odd, then there is no
restriction for $k$.

We are now ready to state the theorem. As we discussed in
Section~\ref{sec:2}, this theorem gives at the same time the
asymptotic behaviour of $n$ non-colliding particles 
on the circle in the random turns vicious walker model.

\begin{Theorem} \label{thm:5-asy}
Let $m$ be a positive integer. Furthermore,
let $\et=(\et_1,\et_2,\dots,\et_n)$ be a vector of integers 
with $m>\et_1>\et_2>\dots>\et_n\ge0$, 
and let $\la=(\la_1,\la_2,\dots,\la_n)$ be a vector of integers
with $m>\la_{s+1}>\dots>\la_n>\la_1>\dots>\la_s\ge0$, for some $s$. 
Then, as $k$ tends to infinity such that \eqref{eq:Bed+-e} holds for
some integer $N$, 
the number of random walks on the $m$-circle
from $\et$ to $\la$ with exactly $k$ standard steps,
such that at no time two coordinates of a point on the random
walk are equal, is asymptotically
\begin{equation} \label{eq:circ-e-asy-even}
\frac {2^{n^2-n+1}} {nm^n}
\(2\frac {\sin\frac {n\pi} {m}} {\sin\frac {\pi} {m}}\)^k
\prod
_{1\le h<t\le n} ^{}\(\sin\frac {\pi(\eta_h-\eta_t)} {m}\cdot
\left\vert\sin\frac {\pi(\la_h-\la_t)} {m}\right\vert\)
\end{equation}
if $n$ is even, and as well
if $n$ is odd and $m$ is even, and it is asymptotically
\begin{equation} \label{eq:circ-e-asy-odd2}
\frac {2^{n^2-n}} {nm^n}
\(2\frac {\sin\frac {n\pi} {m}} {\sin\frac {\pi} {m}}\)^k
\prod
_{1\le h<t\le n} ^{}\(\sin\frac {\pi(\eta_h-\eta_t)} {m}\cdot
\left\vert\sin\frac {\pi(\la_h-\la_t)} {m}\right\vert\)
\end{equation}
if both $n$ and $m$ are odd.
\end{Theorem}

\begin{proof} 
We have to determine the asymptotic behaviour of the coefficient of
$x^k/k!$ in \eqref{eq:circ-e}.
We expand the determinant by linearity in the rows and
obtain
\begin{multline} \label{eq:4}
\frac {1} {nm^n}\sum _{u=0} ^{n-1}
\sum _{r_1,\dots,r_n=0} ^{m-1}
\exp\(2x
\sum _{j=1} ^{n}\cos(2\pi
(u+nr_j)/mn)\)\\
\cdot e^{-2\pi ius/n}
\det_{1\le h,t\le n}\bigg(
e^{-2\pi i(u+nr_h)(\la_t-\eta_h)/mn}\bigg).
\end{multline}
In this expression, we have to extract the coefficient of $x^k/k!$ to
obtain the number of walks with exactly $k$ steps. The expression that
we obtain is a finite sum 
of the form $\sum _{\ell} ^{}c_\ell b_\ell^k$, with the
$b_\ell$'s of the form $2\sum _{j=1}
^{n}\cos(2\pi(u+n r_j)/mn)$, and both the $c_\ell$'s and $b_\ell$'s 
independent of $k$. 
Also in \eqref{eq:4} there is the constraint that we must have
$r_h\ne r_t$ for $h\ne t$ in order to obtain a non-vanishing summand,
because otherwise the determinant is zero.
It is not difficult to see that, because of that, the maximal modulus
of such a $b_\ell$ is equal to
$$2\sum _{j=-n/2} ^{n/2-1}\cos\frac {(2j+1)\pi} {m}=
4\sum _{j=0} ^{n/2-1}\cos\frac {(2j+1)\pi} {m}$$
if $n$ is even, and is equal to
$$2\sum _{j=-(n-1)/2} ^{(n-1)/2}\cos\frac {2j\pi} {m}=
2+4\sum _{j=1} ^{(n-1)/2}\cos\frac {2j\pi} {m}$$
if $n$ is odd. Since we have
$$2\sum _{j=0} ^{n/2-1}\cos\frac {(2j+1)\pi} {m}=\frac {\sin\frac
{n\pi} {m}} {\sin\frac {\pi} {m}}$$
and also
$$1+2\sum _{j=1} ^{(n-1)/2}\cos\frac {2j\pi} {m}=\frac {\sin\frac
{n\pi} {m}} {\sin\frac {\pi} {m}},$$
this maximal modulus is $2\sin(n\pi/m)/\sin(\pi/m)$ in both cases.

If $n$ is even, then the maximal modulus is attained by choosing 
\begin{enumerate}
\item[(e1)] $u=n/2$ and $\{r_1,r_2,\dots,r_n\}=
\{0,1,\dots,\tfrac {n} {2}-1,m-\tfrac {n} {2},\dots,m-2,m-1\}$,
regardless of $m$;
\item[(e2)] $u=n/2$ and $\{r_1,r_2,\dots,r_n\}=
\{\frac {m} {2}-\frac {n} {2},\frac {m} {2}-\frac {n} {2}+1,
\dots,\frac {m} {2}+\frac {n} {2}-1\}$ if $m$ is even;
\item[(e3)] $u=0$ and $\{r_1,r_2,\dots,r_n\}=
\{\frac {m+1} {2}-\frac {n} {2},\frac {m+1} {2}-\frac {n} {2}+1,
\dots,\frac {m+1} {2}+\frac {n} {2}-1\}$ if $m$ is odd.
\end{enumerate}

If $n$ is odd, then the maximal modulus is attained by choosing 
\begin{enumerate}
\item[(o1)] $u=0$ and $\{r_1,r_2,\dots,r_n\}=
\{0,1,\dots,\tfrac {n-1} {2},m-\tfrac {n-1} {2},\dots,m-2,m-1\}$,
regardless of $m$;
\item[(o2)] $u=0$ and $\{r_1,r_2,\dots,r_n\}=
\{\frac {m} {2}-\frac {n-1} {2},\frac {m} {2}-\frac {n-1} {2}+1,
\dots,\frac {m} {2}+\frac {n-1} {2}\}$ if $m$ is even.
\end{enumerate}
(If both $n$ and $m$ are odd, then there is no additional choice
beyond (o1).)

Case (e1) yields the contribution
\begin{equation*} 
\frac {(-1)^s} {nm^n}
\(2\frac {\sin\frac {n\pi} {m}} {\sin\frac {\pi} {m}}\)^k
\det_{1\le h,t\le n}\big(e^{-2\pi i((n+1)/2-h)\la_t/m}\big)
\det_{1\le h,t\le n}\big(e^{2\pi i((n+1)/2-t)\eta_h/m}\big)
\end{equation*}
to the asymptotics of \eqref{eq:circ-e}.
The two determinants are again 
Vandermonde-type determinants. They are
therefore easily evaluated. Thus we obtain 
\begin{equation} \label{eq:5}
\frac {2^{n^2-n}} {nm^n}
\(2\frac {\sin\frac {n\pi} {m}} {\sin\frac {\pi} {m}}\)^k
\prod
_{1\le h<t\le n} ^{}\(\sin\frac {\pi(\eta_t-\eta_h)} {m}\cdot
\left\vert\sin\frac {\pi(\la_t-\la_h)} {m}\right\vert\)
\end{equation}
as the contribution of Case (e1) to the asymptotics. As similar
computations show, the contributions of Cases (e2), (e3), (o1) and
(o2) are also
equal to \eqref{eq:5}.
The claims \eqref{eq:circ-e-asy-even} and \eqref{eq:circ-e-asy-odd2}
follow now upon adding up the corresponding terms in each case.
\end{proof}

Now we are in the position to 
address the question of determining the asymptotic behaviour
of the number of {\it all\/} walks on the $m$-circle 
which start in a given point and
proceed for $k$ standard steps (and, thus, at the same time for the
number of $n$ non-colliding particles on the circle in the random
turns vicious walker model).
Clearly, this amounts to summing the expressions
given in Theorem~\ref{thm:5-asy} over all possible $\la$.
Again, this is just a finite sum, with the number of terms bounded by
$nm^{m}$, a quantity which is independent of
$k$. Thus, it is obvious that the order of magnitude of the number of
all these walks is $\(2{\sin \frac {n\pi}
{m}} /{\sin\frac {\pi} {m}}\)^k$. In order to determine the
multiplicative constant, we have to 
make use of identities featuring Schur functions
and odd orthogonal characters.

\begin{Theorem} \label{thm:5-asy2}
Let $m$ be a positive integer. Furthermore,
let $\et=(\et_1,\et_2,\dots,\et_n)$ be a vector of integers 
with $m>\et_1>\et_2>\dots>\et_n\ge0$.
Then, as $k$ tends to infinity, 
the number of random walks on the $m$-circle which
start at $\et$ and proceed for exactly $k$ standard steps,
such that at no time two coordinates of a point on the random
walk are equal, is asymptotically
\begin{multline} \label{eq:circ-e-asy-even-even1}
\frac {2^{\binom n2}} {m^{n/2}}
\(2\frac {\sin\frac {n\pi} {m}} {\sin\frac {\pi} {m}}\)^k
\prod
_{1\le h<t\le n} ^{}\(\sin\frac {\pi(\eta_h-\eta_t)} {m}\)\\
\times
\(
\prod _{h=1} ^{n/2}\cot\frac {(2h-1)\pi} {2m}
+(-1)^{\vert\et\vert+k+\frac {n} {2}}
\prod _{h=1} ^{n/2}\tan\frac {(2h-1)\pi} {2m}\)
\end{multline}
if both $n$ and $m$ are even,
it is asymptotically
\begin{equation} \label{eq:circ-e-asy-odd-even}
\frac {2^{\binom n2}} {m^{n/2}}
\(2\frac {\sin\frac {n\pi} {m}} {\sin\frac {\pi} {m}}\)^k
\prod
_{1\le h<t\le n} ^{}\(\sin\frac {\pi(\eta_h-\eta_t)} {m}\)
\prod _{h=1} ^{n/2}\cot\frac {(2h-1)\pi} {2m}
\end{equation}
if $n$ is even and $m$ is odd, and it is asymptotically
\begin{equation} \label{eq:circ-e-asy-odd-odd}
\frac {2^{\binom n2}} {m^{(n-1)/2}}
\(2\frac {\sin\frac {n\pi} {m}} {\sin\frac {\pi} {m}}\)^k
\prod
_{1\le h<t\le n} ^{}\(\sin\frac {\pi(\eta_h-\eta_t)} {m}\)
\prod _{h=1} ^{(n-1)/2}\cot\frac {h\pi} {m}
\end{equation}
if $n$ is odd (regardless of $m$).
\end{Theorem}

\begin{proof}
As we already explained, in view of
Theorem~\ref{thm:5-asy}, we have to compute the sum of
\eqref{eq:circ-e-asy-even}, respectively of \eqref{eq:circ-e-asy-odd2}, 
over all possible choices
$(\la_1,\la_2,\dots,\la_n)$. That is, we have to sum these
expressions over all integers $\la_1,\la_2,\dots,\la_n$ with
$m>\la_1>\la_2>\dots>\la_n\ge0$, and its cyclic permutations, {\it
such that walks are possible from $\et$ to
$\la=(\la_1,\la_2,\dots,\la_n)$}, respectively to its cyclic
permutations, {\it in $k$ steps}. What the latter means, was
discussed in the paragraph containing \eqref{eq:Bed+-e}. Thus, we 
have to distinguish between several cases.

Before we list these cases, the reader should observe that the
expressions \eqref{eq:circ-e-asy-even} and \eqref{eq:circ-e-asy-odd2}
are invariant under cyclic permutations of $\la_1,\la_2,\dots,\la_n$.
Thus, as the first step, we will multiply them by $n$, respectively
by $n/2$, depending on whether {\it any} cyclic permutation of $\la$
can be reached from $\et$ in $k$ steps, or only {\it every second}.

Now, if both $m$ and $n$ are odd, then {\it any} $\la$, and {\it any}
cyclic permutation of it can be reached from $\et$ if $k$ is large
enough. Thus, in this case, we have to multiply the 
expression \eqref{eq:circ-e-asy-odd2}
by $n$, and subsequently sum it over all integers 
$\la_1,\la_2,\dots,\la_n$ with
$m>\la_1>\la_2>\dots>\la_n\ge0$.

If $m$ is odd but $n$ is even, then only {\it every second\/}
cyclic permutation of a given $\la$ can be reached from $\et$ for a given
$k$ which is large enough. Thus, in this case, we have to multiply the 
expression \eqref{eq:circ-e-asy-even}
by $n/2$, and subsequently sum it over all integers 
$\la_1,\la_2,\dots,\la_n$ with
$m>\la_1>\la_2>\dots>\la_n\ge0$.

On the other hand, if $m$ is even, then \eqref{eq:Bed+-e} implies
that $\vert\la\vert\equiv \vert\et\vert+k$ mod~$2$. In particular, given a $\la$
satisfying this condition, {\it every}
cyclic permutation of $\la$ can be reached from $\et$ for a given
$k$ which is large enough. 
Thus, in this case, we have to multiply the 
expression \eqref{eq:circ-e-asy-even}
by $n$, and subsequently sum it over all integers 
$\la_1,\la_2,\dots,\la_n$ with
$m>\la_1>\la_2>\dots>\la_n\ge0$ and
$\vert\la\vert\equiv\vert\et\vert+k$ mod~$2$.

In summary, what we need is, on the one hand, the sum
\begin{equation} \label{eq:sinsum}
S_1=\sum _{m>\la_1>\dots>\la_n\ge0} ^{}\prod
_{1\le h<t\le n} ^{}\sin\frac {\pi(\la_h-\la_t)} {m},
\end{equation}
and, on the other hand, the same sum \eqref{eq:sinsum}, but restricted to
those $\la$ for which $\vert\la\vert\equiv\vert\et\vert+k$ mod~$2$. 
The latter sum is equal
to $\frac {1} {2}(S_1+(-1)^{\vert\et\vert+k}S_2)$, with
\begin{equation} \label{eq:sinsum-}
S_2=\sum _{m>\la_1>\dots>\la_n\ge0} ^{}(-1)^{\vert\la\vert}\prod
_{1\le h<t\le n} ^{}\sin\frac {\pi(\la_h-\la_t)} {m}.
\end{equation}

We may again write the summand in the sums $S_1$ and $S_2$
as a Vandermonde-type determinant, 
so that \eqref{eq:sinsum} becomes
\begin{equation*} 
(2i)^{-\binom n2}\sum _{m>\la_1>\dots>\la_n\ge0} ^{}
\det_{1\le h,t\le n}\big(e^{2\pi i((n+1)/2-h)\la_t/m}\big),
\end{equation*}
and the expression \eqref{eq:sinsum-} becomes
\begin{equation*} 
(2i)^{-\binom n2}\sum _{m>\la_1>\dots>\la_n\ge0} ^{}
\det_{1\le h,t\le n}\Big(\big(-e^{2\pi i((n+1)/2-h)/m}\big)^{\la_t}\Big).
\end{equation*}
Upon replacing $\la_j$ by $\la_j+n-j$, these expressions are
transformed to
\begin{equation} \label{eq:sumsindet}
(2i)^{-\binom n2}\sum _{m-n\ge \la_1\ge\dots\ge\la_n\ge0} ^{}
\det_{1\le h,t\le n}\Big(\big(e^{2\pi
i((n+1)/2-h)/m}\big)^{\la_t+n-t}\Big),
\end{equation}
respectively
\begin{equation} \label{eq:sumsindet-}
(2i)^{-\binom n2}\sum _{m-n\ge \la_1\ge\dots\ge\la_n\ge0} ^{}
\det_{1\le h,t\le n}\Big(\big(-e^{2\pi
i((n+1)/2-h)/m}\big)^{\la_t+n-t}\Big).
\end{equation}
We may rewrite the latter two determinants using Schur functions
(see \eqref{e4} for the definition).
Using this notation, we may write \eqref{eq:sumsindet} and
\eqref{eq:sumsindet-} in the form
\begin{equation} \label{eq:sumschur}
(2i)^{-\binom n2}\det_{1\le h,t\le n}\big((\ep q^{n+1-2h})^{n-t}\big)
\sum _{m-n\ge \la_1\ge\dots\ge\la_n\ge0} ^{}
s_\la(\ep q^{n-1},\ep q^{n-3},\dots,\ep q^{-n+3},\ep q^{-n+1}),
\end{equation}
where $q=e^{\pi i/m}$, with $\ep=1$ to yield equality with
\eqref{eq:sumsindet}, and $\ep=-1$ to yield equality with 
\eqref{eq:sumsindet-}.
The determinant in this expression is a
Vandermonde determinant, and is therefore easily evaluated. 
The sum
over Schur functions, on the other hand, can be evaluated by means of
(see \cite[proof of Theorem~2]{KrGVAA} for a discussion of this identity, with
references to various proofs)
\begin{equation} \label{eq:3}
\sum _{p\ge \la_1\ge \dots\ge\la_n\ge0} ^{}s_\la(x_1,x_2,\dots,x_n)
=(x_1x_2\cdots
x_n)^{p/2}\so^{odd}_{\big((p/2)^n\big)}
(x_1,x_2,\dots,x_n),
\end{equation}
where $\so^{odd}_\la(x_1,x_2,\dots,x_n)$ is an {\it odd orthogonal
character}. (See \eqref{e12} for the definition.
The notation $\big((p/2)^n\big)$ in \eqref{eq:3} means a vector of
$n$ components, all of them equal to $p/2$.)
Thus, the expression \eqref{eq:sumschur} becomes
\begin{equation} \label{eq:orth}
\bigg(\ep^{\binom n2}\prod _{1\le h<t\le n} ^{}\sin\frac {\pi(t-h)} {m}\bigg)
\ep^{n(m-n)/2}\so^{odd}_{\big(((m-n)/2)^n\big)}(\ep q^{n-1},\ep q^{n-3},\dots,\ep
q^{-n+3},\ep q^{-n+1}).
\end{equation}
If $\ep=1$, the odd orthogonal character, specialized in this manner, is
evaluated in Lemma~\ref{lem:so1}, while for $\ep=-1$ this is done
in Lemma~\ref{lem:so2}. If the results are substituted in
\eqref{eq:orth}, which, as we argued above, is in fact equal to the sum
\eqref{eq:sinsum}, respectively to \eqref{eq:sinsum-}, the sums 
that we wanted to evaluate, then the claims 
\eqref{eq:circ-e-asy-even-even1}--\eqref{eq:circ-e-asy-odd-odd}
follow after some further straight-forward (but tedious) calculations.
\end{proof}

The next two results address the asymptotic behaviour of the number of
walks on the
$m$-circle consisting of diagonal steps. As we discussed in
Section~\ref{sec:2}, these theorems give at the same time the
asymptotic behaviour of $n$ non-colliding particles 
on the circle in the lock-step vicious walker model.
As before, we begin with the result which described the asymptotic
behaviour of the number of walks with fixed starting and end point.

\begin{Theorem} \label{thm:4-asy}
Let $m$ be a positive integer. Furthermore,
let $\et=(\et_1,\et_2,\dots,\et_n)$ be a vector of integers or of
half-integers with $m>\et_1>\et_2>\dots>\et_n\ge0$, 
and let $\la=(\la_1,\la_2,\dots,\la_n)$ be a vector of integers or of
half-integers
with $m>\la_{s+1}>\dots>\la_n>\la_1>\dots>\la_s\ge0$, for some $s$. 
Then, as $k$ tends to infinity such that $k\equiv 2\et_j+2\la_j$ mod\/~$2$, 
the number of random walks on the $m$-circle
from $\et$ to $\la$ with exactly $k$ diagonal
steps, such that at no time two coordinates of a point on the random
walk are equal, is asymptotically
\begin{equation} \label{eq:circ-d-asy-even}
\frac {2^{n^2-n}} {nm^n}
\bigg(2^{n}\prod _{j=1} ^{n}\cos\frac {\pi\(j-\frac {n+1} {2}\)} {m}\bigg)^{k}
\prod
_{1\le h<t\le n} ^{}\(\sin\frac {\pi(\eta_h-\eta_t)} {m}\cdot
\left\vert\sin\frac {\pi(\la_h-\la_t)} {m}\right\vert\).
\end{equation}
\end{Theorem}

\begin{proof} 
Clearly, this time we want to determine the asymptotic behaviour of
the expression \eqref{eq:circ-d}.
First of all, if $k\equiv 2\et_j+\la_j$ mod~$2$, i.e., 
in the case where $k$, the $\eta_j$'s and
the $\la_j$'s are chosen so that walks exist, then we may replace the
sum over $r$ in \eqref{eq:circ-d} by twice the sum of the same summand, but 
where $r$ runs from $0$ to $m-1$
(instead of $2m-1$; here it is important that $m$ is an integer).
Then we expand again the determinant by linearity,
and obtain
\begin{multline} \label{eq:2}
\frac {1} {nm^n}\sum _{u=0} ^{n-1}
\sum _{r_1,\dots,r_n=0} ^{m-1}
\prod _{j=1} ^{n}\(2\cos(\pi
(u+nr_j)/mn)\)^k\\
\cdot e^{-2\pi ius/n}
\det_{1\le h,t\le n}\big(e^{-2\pi i(u+nr_h)(\la_t-\eta_h)/mn}\big).  
\end{multline}
This is again a finite sum of the form $\sum _{\ell} ^{}c_\ell
b_\ell^k$, with the $b_\ell$'s of the form
$2^n\prod _{j=1} ^{n}\cos(\pi
(u+nr_j)/mn)$, and both the $c_\ell$'s and $b_\ell$'s independent of $k$.
Also in \eqref{eq:2} there is the constraint that we must have
$r_h\ne r_t$ for $h\ne t$ in order to obtain a non-vanishing summand,
because otherwise the determinant is zero.
Because of that, one discovers that in order to have $b_\ell$ with maximal
modulus we must have $u=n/2$ and 
$$\{r_1,r_2,\dots,r_n\}=
\{0,1,\dots,\tfrac {n} {2}-1,m-\tfrac {n} {2},\dots,m-2,m-1\}$$
if $n$ is even, and we must have $u=0$ and 
$$\{r_1,r_2,\dots,r_n\}=
\{0,1,\dots,\tfrac {n-1} {2},m-\tfrac {n-1} {2},\dots,m-2,m-1\}$$
if $n$ is odd.

Let first $n$ be even. Then we obtain for the asymptotics of 
\eqref{eq:2} the expression
\begin{multline*} 
\frac {(-1)^s} {nm^n}
\bigg(2^{n}\prod _{j=1} ^{n}\cos\frac {\pi\(j-\frac {n+1} {2}\)} {m}\bigg)^{k}\\
\times
\det_{1\le h,t\le n}\big(e^{-2\pi i((n+1)/2-h)\la_t/m}\big)
\det_{1\le h,t\le n}\big(e^{2\pi i((n+1)/2-t)\eta_h/m}\big).  
\end{multline*}
Both determinants are essentially Vandermonde determinants and are therefore
easily evaluated. The result is exactly \eqref{eq:circ-d-asy-even}.

If $n$ is odd then we obtain for the asymptotics of 
\eqref{eq:2} the expression
\begin{equation*} 
\frac {1} {nm^n}
\bigg(2^{n}\prod _{j=1} ^{n}\cos\frac {\pi\(j-\frac {n+1} {2}\)} {m}\bigg)^{k}
\det_{1\le h,t\le n}\big(e^{-2\pi i((n+1)/2-h)\la_t/m}\big)
\det_{1\le h,t\le n}\big(e^{2\pi i((n+1)/2-t)\eta_h/m}\big),
\end{equation*}
which becomes again \eqref{eq:circ-d-asy-even} if 
the Vandermonde-type determinants are evaluated.
\end{proof}

By summation of the corresponding expressions, the previous result
allows us now to derive the asymptotic behaviour of walks with a fixed
starting point but with arbitrary end point. The sums that need to be
carried out are equivalent to some of those that we already evaluated in the
proof of Theorem~\ref{thm:5-asy2}. 

\begin{Theorem} \label{thm:4-asy2}
Let $m$ be a positive integer. Furthermore,
let $\et=(\et_1,\et_2,\dots,\et_n)$ be a vector of integers or of
half-integers with $m>\et_1>\et_2>\dots>\et_n\ge0$.
Then, as $k$ tends to infinity, the number of random walks on the $m$-circle
which start at $\et$ and proceed for exactly $k$ diagonal
steps, such that at no time two coordinates of a point on the random
walk are equal, is asymptotically
\begin{equation} \label{eq:circ-d-asy-even2}
\frac {2^{\binom n2}} {m^{n/2}}
\bigg(2^{n}\prod _{j=1} ^{n}\cos\frac {\pi\(j-\frac {n+1} {2}\)} {m}\bigg)^{k}
\prod
_{1\le h<t\le n} ^{}\(\sin\frac {\pi(\eta_h-\eta_t)} {m}\)
\prod _{h=1} ^{n/2}\cot\frac {(2h-1)\pi} {2m}
\end{equation}
if $n$ is even, and it is asymptotically
\begin{equation} \label{eq:circ-d-asy-odd2}
\frac {2^{\binom n2}} {m^{(n-1)/2}}
\bigg(2^{n}\prod _{j=1} ^{n}\cos\frac {\pi\(j-\frac {n+1} {2}\)} {m}\bigg)^{k}
\prod
_{1\le h<t\le n} ^{}\(\sin\frac {\pi(\eta_h-\eta_t)} {m}\)
\prod _{h=1} ^{(n-1)/2}\cot\frac {h\pi} {m}
\end{equation}
if $n$ is odd.
\end{Theorem}

\begin{proof} 
Clearly, we have to sum \eqref{eq:circ-d-asy-even} over all possible choices
of $\la$. Depending on the parity of $k+\et_j$,
this means to take the sum over all integers $\la_1,\la_2,\dots,\la_n$ with
$m>\la_1>\la_2>\dots>\la_n\ge0$, or over all half-integers with the
same property, and over all their cyclic permutations. Since the
expression \eqref{eq:circ-d-asy-even} is independent of $s$, every cyclic
permutation yields the same value. Therefore we have to multiply
this expression
by $n$, and subsequently sum it over all integers, respectively
half-integers, $\la_1,\la_2,\dots,\la_n$ with 
$m>\la_1>\la_2>\dots>\la_n\ge0$. 

So, what we need is the sum
\begin{multline} \label{eq:sinsum1}
{\sum _{m>\la_1>\dots>\la_n\ge0} ^{}\kern-.8cm}{}^{\displaystyle\prime}
\kern.8cm
\prod
_{1\le h<t\le n} ^{}\sin\frac {\pi(\la_h-\la_t)} {m}\\
=(2i)^{-\binom n2}{\sum _{m>\la_1>\dots>\la_n\ge0} ^{}
\kern-.8cm}{}^{\displaystyle\prime}
\kern.8cm
\det_{1\le h,t\le n}\big(e^{2\pi i((n+1)/2-h)\la_t/m}\big),
\end{multline}
where the sum $\sum{}^{\textstyle\prime}$ is over all integral
$\la_1,\la_2,\dots,\la_n$, but also the sum \eqref{eq:sinsum1}
where $\sum{}^{\textstyle\prime}$ is restricted to {\it half-integral\/}
$\la_1,\la_2,\dots,\la_n$. (The equality of the two expressions
in \eqref{eq:sinsum1} follows 
again from the Vandermonde determinant evaluation.)

The sum over all integers $\la_1,\la_2,\dots,\la_n$ has already been
evaluated in the proof of Theorem~\ref{thm:5-asy2}, when we evaluated
$S_1$. If we want to
form the sum \eqref{eq:sinsum1} over all half-integers
$\la_1,\la_2,\dots,\la_n$, then we may replace $\la_j$ by $\la_j+n-j+\frac
{1} {2}$, and rewrite it as
\begin{multline*} 
(2i)^{-\binom n2}
\sum _{m-n\ge \la_1\ge \dots\ge \la_n\ge0} ^{}
\det_{1\le h,t\le n}\big(e^{2\pi
i((n+1)/2-h)(\la_t+n-t+1/2)/m}\big)\\
=(2i)^{-\binom n2}
\sum _{m-n\ge \la_1\ge \dots\ge \la_n\ge0} ^{}
\det_{1\le h,t\le n}\big(e^{2\pi
i((n+1)/2-h)(\la_t+n-t)/m}\big),
\end{multline*}
where the sums are now over {\it integral} $\la_1,\la_2,\dots,\la_n$.
The sum in the last line 
is exactly the same sum as \eqref{eq:sumsindet}, which is in turn 
equal to $S_1$ (and, thus, also to
the sum \eqref{eq:sinsum1} when $\sum{}^{\textstyle\prime}$ is
taken over all integral $\la_1,\la_2,\dots,\la_n$). Therefore, regardless of
the parity of $k$,
the result is the same. Since the evaluation of 
$S_1$ in the proof of Theorem~\ref{thm:5-asy2} yielded
two different expressions depending on whether $n$ is even or odd, we
obtain the two cases in the statement of the theorem.
\end{proof}

\section{Asymptotics for random walks in alcoves of type $\tilde C$} 
\label{sec:5}

The subject of this section is the determination of the asymptotic
behaviour of the number of walks from a given starting point to a
given end point 
which stay in the alcove $\mathcal A^{\tilde C_{n}}_m$ of type $\tilde C_n$ as the
number of steps becomes large, as well as the asymptotic behaviour
of the number of those walks which start at a given point but may
terminate anywhere. In technical terms, we determine the asymptotic
behaviour of the expressions given by
Theorems~\ref{thm:6} and \ref{thm:7} as $k$ becomes large,
and as well if these expressions
are summed over all possible end points of the walks.
In fact, for Theorem~\ref{thm:7}, i.e., for the case of diagonal
steps, this had already been carried out in
\cite{KrGVAB}, so that we only copy the corresponding results for the
sake of completeness; see Theorems~\ref{thm:7-asy} and
\ref{thm:7-asy2} below.

Before, however, we address the case of standard steps. As we
discussed in Section~\ref{sec:2}, this case is also equivalent to the
movements of $n$
non-colliding particles in an interval according to the random turns
vicious walker model. We begin,
as usual, with the corresponding results when starting and end point
are fixed.

\begin{Theorem} \label{thm:6-asy}
Let $m$ be a positive integer. Furthermore,
let $\et=(\et_1,\et_2,\dots,\et_n)$ and 
$\la=(\la_1,\la_2,\dots,\la_n)$ be vectors of integers in the
alcove $\mathcal A^{\tilde C_{n}}_m$ of type $\tilde C_n$ {\em(}defined in \eqref{eq:alcC}{\em)}.
Then, as $k$ tends to infinity such that $k\equiv
\vert\et\vert+\vert\la\vert$ mod\/~$2$, 
the number of random walks from $\et$ to $\la$ with exactly $k$
standard steps,
which stay in the alcove $\mathcal A^{\tilde C_{n}}_m$, is asymptotically
\begin{multline} \label{eq:C+-e-asy}
\kern-2pt\frac {2^{2n^2-n+1}} {m^n}
\(\frac {2\sin\frac {n\pi} {2m}\,\cos\frac {(n+1)\pi} {2m}}
{\sin\frac {\pi} {2m}}\) ^k 
\prod
_{1\le h<t\le n} ^{}\(\sin\frac {\pi(\eta_h-\eta_t)} {2m}\cdot
\sin\frac {\pi(\la_h-\la_t)} {2m}\)\\
\times
\prod _{1\le h\le
t\le n} ^{}\(\sin\frac {\pi(\eta_h+\eta_t)} {2m}\cdot
\sin\frac {\pi(\la_h+\la_t)} {2m}\).
\end{multline}
\end{Theorem}

\begin{proof} We have to determine the asymptotic behaviour of the
coefficient of $x^k/k!$ in \eqref{eq:C+-e}.
If we expand expression \eqref{eq:C+-e}, i.e., if we 
use linearity of the determinant in the rows, then we obtain the expression
\begin{multline} \label{eq:1} 
\frac {1} {m^n}\sum _{r_1,\dots,r_n=1} ^{2m-1}\exp\Big(
\hbox{$2x\sum _{j=1}
^{n}\cos(\pi r_j/m)$}\Big)\bigg(\prod _{j=1} ^{n}\sin(\pi r_j\eta_j/m)\bigg)\\
\cdot \det_{1\le h,t\le n}\big(\sin(\pi r_h\la_t/m)\big).
\end{multline}
As we said, 
in this expression we have to extract the coefficient of $x^k/k!$ to
obtain the number of walks with exactly $k$ steps. The expression that
we obtain is a finite sum 
of the form $\sum _{\ell} ^{}c_\ell b_\ell^k$, with the
$b_\ell$'s of the form $2\sum _{j=1}
^{n}\cos(\pi r_j/m)$, and both the $c_\ell$'s and $b_\ell$'s 
independent of $k$. 

Finding the asymptotics of \eqref{eq:1}
means to find the $b_\ell$'s
with largest modulus. This, in turn, means to choose the parameters $r_h$
either close to the lower limit of the summation, 1, respectively 
close to the upper
limit, $2m-1$, or close to $m-1$, respectively close to $m+1$. 
(In the first case, all the cosines $\cos(\pi r_h/m)$ will be close to
1, whereas in the second case all of them will be close to $-1$.
If some $r_j$ is
equal to $m$, then the corresponding term vanishes because of the
expression $\sin(\pi r_j\et_j/m)$ occurring in the summand.)
There are again restrictions however:
if $r_h=r_t$ for
$h\ne t$ then the determinant in \eqref{eq:1} vanishes, as well as if
$r_h=2m-r_t$ for some $h$ and $t$. Therefore we may restrict
ourselves to the cases where $r_h\ne r_t$ and $r_h\ne 2m-r_t$ for all
$h$ and $t$. 

Hence, we will choose the set $\{r_1,r_2,\dots,r_n\}$ either from
$$\{1,2,\dots,n, 2m-n,\dots,2m-2,2m-1\}$$ 
or from
$$\{m-n,\dots,m-2,m-1,m+1,m+2,\dots,m+n\},$$
in such a way that the $r_j$'s are distinct and $r_h\ne 2m-r_t$ for all 
$h$ and $t$. Clearly, there are $2^n$ sets of the first type, and 
$2^n$ sets as well of the second type.
As is not difficult to see, for each fixed set, the sum of the
corresponding terms $c_\ell b_\ell^k$ is equal to
\begin{equation} \label{eq:cont1} 
\frac {1} {m^n}
\bigg(2\sum _{j=1}^{n}\cos(\pi j/m)\bigg)^k
\cdot \det_{1\le h,t\le n}\big(\sin(\pi t\eta_h/m)\big)
\cdot \det_{1\le h,t\le n}\big(\sin(\pi h\la_t/m)\big)
\end{equation}
in both cases.

Now we have
$$\sum _{j=1} ^{n}\cos(\pi j/m)=\frac
{\sin(n\pi/2m)\,\cos((n+1)\pi/2m)} {\sin(\pi/2m)}$$
and
\begin{equation} \label{eq:sindet}
\det_{1\le h,t\le n}\big(\sin(\pi h\la_t/m)\big)
=2^{n^2-n}\prod
_{1\le h<t\le n} ^{}\sin\frac {\pi(\la_h-\la_t)} {2m}\prod _{1\le h\le
t\le n} ^{}\sin\frac {\pi(\la_h+\la_t)} {2m}.
\end{equation}
The latter identity follows from writing the determinant as
$$
\det_{1\le h,t\le n}\big(\sin(\pi h\la_t/m)\big)=
(2i)^{-n}\det_{1\le h,t\le n}\big(e^{\pi ih\la_t/m}
-e^{-\pi ih\la_t/m}\big),
$$
and evaluating it by means of \eqref{eq:sympl}.
Substituting this in \eqref{eq:cont1}, and multiplying the resulting
expression by $2\cdot 2^n=2^{n+1}$ (the number of these sets
$\{r_1,r_2,\dots,r_n\}$), we obtain \eqref{eq:C+-e-asy}.
\end{proof}

Having accomplished the asymptotic analysis of the walks with fixed
starting and end point, we can now turn to the analysis of the walks
with fixed starting point but arbitrary end point. Again, this amounts
to a summation problem, namely summing expression \eqref{eq:C+-e-asy}
over all possible $\la$. To carry out this task, we make use of
identities featuring Schur functions and symplectic characters.

\begin{Theorem} \label{thm:6-asy2}
Let $m$ be a positive integer. Furthermore,
let $\et=(\et_1,\et_2,\dots,\et_n)$ be a vector of integers in the
alcove $\mathcal A^{\tilde C_{n}}_m$ of type $\tilde C_n$ {\em(}defined in \eqref{eq:alcC}{\em)}.
Then, as $k$ tends to infinity, 
the number of random walks which start at $\et$ and proceed for exactly $k$
standard steps,
which stay in the alcove $\mathcal A^{\tilde C_{n}}_m$, is asymptotically
\begin{multline} \label{eq:C+-e-asy2a}
\frac {2^{2n^2-n}} {m^n}
\(\frac {2\sin\frac {n\pi} {2m}\,\cos\frac {(n+1)\pi} {2m}}
{\sin\frac {\pi} {2m}}\) ^k 
 \prod _{1\le h<t\le n} ^{}\(\sin \frac {\pi  (\et_h-\et_t)} {2m} \cdot
    \sin \frac {\pi  (t-h)} {2m}\)\\
\times
 \prod _{1\le h\le t\le n} ^{}\(\sin \frac {\pi  (\et_t+\et_h)} {2m} 
\cdot   \sin \frac {\pi  (t+h)} {2m}\)\\
\times
\(
 \prod _{h=0} ^{n}\prod _{t=1} ^{n}\frac {\sin \frac 
   {\pi  (t-h+{m-1})} 
   {2m}} {\sin \frac {\pi  (t-h+n)} {2m}}
+(-1)^{\vert\et\vert+k+\frac {n} {2}}
2^{-n^2}\frac {\prodl _{h=1} ^{n/2}\tan^2\frac {(2h-1)\pi} {2m}}
{\prodl _{h=1} ^{n+1}\prodl _{t=1} ^{n}\left\vert\sin\frac {(2t-2h+1)\pi}
{2m}\right\vert}
\)
\end{multline}
if both $m$ and $n$ are even, it is asymptotically
\begin{multline} \label{eq:C+-e-asy2b}
\frac {2^{2n^2-n}} {m^n}
\(\frac {2\sin\frac {n\pi} {2m}\,\cos\frac {(n+1)\pi} {2m}}
{\sin\frac {\pi} {2m}}\) ^k 
 \prod _{1\le h<t\le n} ^{}\(\sin \frac {\pi  (\et_h-\et_t)} {2m} \cdot
    \sin \frac {\pi  (t-h)} {2m}\)\\
\times
 \prod _{1\le h\le t\le n} ^{}\(\sin \frac {\pi  (\et_t+\et_h)} {2m} 
\cdot   \sin \frac {\pi  (t+h)} {2m}\)\\
\times
\(
 \prod _{h=0} ^{n}\prod _{t=1} ^{n}\frac {\sin \frac 
   {\pi  (t-h+{m-1})} 
   {2m}} {\sin \frac {\pi  (t-h+n)} {2m}}
+(-1)^{\vert\et\vert+k+\frac {n+1} {2}}
2^{-n^2}\frac {\prodl _{h=1} ^{(n+1)/2}\tan^2\frac {(2h-1)\pi} {2m}}
{\prodl _{h=1} ^{n+1}\prodl _{t=1} ^{n}\left\vert\sin\frac {(2t-2h+1)\pi}
{2m}\right\vert}
\)
\end{multline}
if $m$ is even and $n$ is odd, it is asymptotically
\begin{multline} \label{eq:C+-e-asy2c}
\frac {2^{2n^2-n}} {m^n}
\(\frac {2\sin\frac {n\pi} {2m}\,\cos\frac {(n+1)\pi} {2m}}
{\sin\frac {\pi} {2m}}\) ^k 
 \prod _{1\le h<t\le n} ^{}\(\sin \frac {\pi  (\et_h-\et_t)} {2m} \cdot
    \sin \frac {\pi  (t-h)} {2m}\)\\
\times
 \prod _{1\le h\le t\le n} ^{}\(\sin \frac {\pi  (\et_t+\et_h)} {2m} 
\cdot   \sin \frac {\pi  (t+h)} {2m}\)\\
\times
\(
 \prod _{h=0} ^{n}\prod _{t=1} ^{n}\frac {\sin \frac 
   {\pi  (t-h+{m-1})} 
   {2m}} {\sin \frac {\pi  (t-h+n)} {2m}}
+(-1)^{\vert\et\vert+k+\frac {n} {2}}
2^{-n^2}\prodl _{h=1} ^{n/2}\frac
{\sin^2\frac {(2h-1)\pi} {2m}}
{\cos^2\frac {h\pi} {m}}
\frac {1}
{\prodl _{h=1} ^{n+1}\prodl _{t=1} ^{n}\left\vert\sin\frac {(2t-2h+1)\pi}
{2m}\right\vert}
\)
\end{multline}
if $m$ is odd and $n$ is even, and it is asymptotically
\begin{multline} \label{eq:C+-e-asy2d}
\frac {2^{2n^2-n}} {m^n}
\(\frac {2\sin\frac {n\pi} {2m}\,\cos\frac {(n+1)\pi} {2m}}
{\sin\frac {\pi} {2m}}\) ^k 
 \prod _{1\le h<t\le n} ^{}\(\sin \frac {\pi  (\et_h-\et_t)} {2m} \cdot
    \sin \frac {\pi  (t-h)} {2m}\)\\
\times
 \prod _{1\le h\le t\le n} ^{}\(\sin \frac {\pi  (\et_t+\et_h)} {2m} 
\cdot   \sin \frac {\pi  (t+h)} {2m}\)
\(
 \prod _{h=0} ^{n}\prod _{t=1} ^{n}\frac {\sin \frac 
   {\pi  (t-h+{m-1})} 
   {2m}} {\sin \frac {\pi  (t-h+n)} {2m}}
\)
\end{multline}
if both $m$ and $n$ are odd.
\end{Theorem}

\begin{proof}
As we already observed, 
in view of Theorem~\ref{thm:6-asy}, we have to carry out the
sum of \eqref{eq:C+-e-asy} over all possible $\la$, i.e., over all
$m>\la_1>\la_2>\dots>\la_n>0$, where $\vert\la\vert\equiv
k+\vert\et\vert$ mod~$2$. Leaving away factors which are independent of
$\la$, the problem is to compute the sum
\begin{equation} \label{eq:S}
{\sum _{m>\la_1>\dots>\la_n>0} ^{}\kern-.8cm}{}^{\displaystyle\prime}
\kern.8cm
\prod
_{1\le h<t\le n} ^{}\(
\sin\frac {\pi(\la_h-\la_t)} {2m}\)
\prod _{1\le h\le t\le n} ^{}\(
\sin\frac {\pi(\la_h+\la_t)} {2m}\),
\end{equation}
where the sum $\sum{}^{\textstyle\prime}$ 
is either restricted to those $\la$ for which $\vert\la\vert$
is even, or to those for which $\vert\la\vert$ is odd, depending on
the parity of $k+\vert\et\vert$. This task will be accomplished if we
are able to evaluate the (complete) sum
\begin{equation*} \label{eq:S1}
T_1=\sum _{m>\la_1>\dots>\la_n>0} ^{}
\prod
_{1\le h<t\le n} ^{}\(
\sin\frac {\pi(\la_h-\la_t)} {2m}\)
\prod _{1\le h\le t\le n} ^{}\(
\sin\frac {\pi(\la_h+\la_t)} {2m}\),
\end{equation*}
and its ``signed variant''
\begin{equation*} \label{eq:S2}
T_2=\sum _{m>\la_1>\dots>\la_n>0} ^{}
(-1)^{\vert\la\vert}\prod
_{1\le h<t\le n} ^{}\(
\sin\frac {\pi(\la_h-\la_t)} {2m}\)
\prod _{1\le h\le t\le n} ^{}\(
\sin\frac {\pi(\la_h+\la_t)} {2m}\).
\end{equation*}
The sum \eqref{eq:S} is then equal to $\frac {1} {2}(T_1+T_2)$ if
the sum $\sum{}^{\textstyle\prime}$ is 
restricted to the $\la$'s for which $\vert\la\vert$ is even, and it is
equal to $\frac {1} {2}(T_1-T_2)$ if
the sum $\sum{}^{\textstyle\prime}$  is
restricted to the $\la$'s for which $\vert\la\vert$ is odd.

A sum equivalent to $T_1$ had already been evaluated in 
\cite[first part of the proof of
Theorem~6]{KrGVAB}. The result is
$$
T_1=
 \prod _{1\le h<t\le n} ^{}
    \sin \frac {\pi  (t-h)} {2m}
 \prod _{1\le h\le t\le n} ^{}
   \sin \frac {\pi  (t+h)} {2m}
 \prod _{h=0} ^{n}\prod _{t=1} ^{n}\frac {\sin \frac 
   {\pi  (t-h+{m-1})} 
   {2m}} {\sin \frac {\pi  (t-h+n)} {2m}}.
$$

 In order to evaluate $T_2$, we proceed in a
manner similar to the evaluation in \cite{KrGVAB}. By means
of \eqref{eq:sindet}, we may rewrite $T_2$ as 
\begin{align*} \label{}
\sum _{m>\la_1>\dots>\la_n>0} ^{}&
(-1)^{\vert\la\vert}
2^{-n^2+n}
\det_{1\le h,t\le n}\big(\sin(\pi h\la_t/m)\big)\\
&=\sum _{m>\la_1>\dots>\la_n>0} ^{}
\frac {1} {2^{n^2}i^n}
\det_{1\le h,t\le n}\big((-e^{\pi i h/m})^{\la_t}-(-e^{\pi i
h/m})^{-\la_t}\big).
\end{align*}
Replacing $\la_j$ by $\la_j+n-j+1$, $j=1,2,\dots,n$, we obtain the
expression
\begin{equation} \label{eq:Summe}
\sum _{m-n-1\ge\la_1\ge\dots\ge\la_n\ge0} ^{}
\frac {1} {2^{n^2}i^n}
\det_{1\le h,t\le n}\big((-e^{\pi i h/m})^{\la_t+n-t+1}-(-e^{\pi i
h/m})^{-(\la_t+n-t+1)}\big).
\end{equation}
This determinant can be expressed in terms of a
{\it symplectic character}. Given a partition 
$\la=(\la_1,\la_2,\dots,\la_n)$ 
(i.e., a non-increasing sequence of non-negative integers), 
the symplectic character
$\sp_\la(x_1,x_2,\dots,x_n)$ is
defined by (see \cite[(24.18)]{FuHaAA})
\begin{equation} \label{e24}
\sp_\la(x_1,x_2,\dots,x_n)
=\frac {\det\limits_{1\le h,t\le
n}(x_h^{\la_t+n-t+1}-x_h^{-(\la_t+n-t+1)})} 
{\det\limits_{1\le h,t\le
n}(x_h^{n-t+1}-x_h^{-(n-t+1)})}.
\end{equation}
Therefore, writing again $q$ for $e^{\pi i/m}$, the sum in \eqref{eq:Summe}
equals
\begin{equation} \label{eq:sum-sympl}
\frac {1} {2^{n^2}i^n}
\det_{1\le h,t\le n}
\( (-q^{h})^{n-t+1}-(-q^{h})^{-(n-t+1)}\)
\sum _{m-n-1\ge\la_1\ge\dots\ge\la_n\ge0}
\sp_\la(-q,-q^2,\dots,-q^n).
\end{equation}
Now we appeal to the formula (see
\cite[(3.4)]{KratBC}),
\begin{equation} \label{e31}
s_{(c^r)}(x_1,x_1^{-1},\dots,x_n,x_n^{-1},1)=
\sum _{c\ge\nu_1\ge\nu_2\ge\dots\ge\nu_r\ge0} ^{}\sp_\nu(x_1,\dots,x_n),
\end{equation}
which is valid for $r\le n$,
where on the left-hand side we have again a Schur function (cf.\
\eqref{e4} for the definition). The notation $(c^r)$ is short for the
vector in which the first $r$ components are equal to $c$, followed by
$2n+1-r$ components all of which are 0.
Use of this formula in \eqref{eq:sum-sympl} yields the equivalent expression
\begin{multline*} 
\frac {1} {2^{n^2}i^n}
\det_{1\le h,t\le n}
\( (-q^{h})^{n-t+1}-(-q^{h})^{-(n-t+1)}\)\\
\times
s_{\big((m-n-1)^n\big)}(-q^n,-q^{n-1},\dots,-q,1,-q^{-1},\dots,
-q^{-n+1},-q^n)\\
=(-1)^{(m-n-1)n}\frac {1} {2^{n^2}i^n}
\det_{1\le h,t\le n}
\( (-q^{h})^{n-t+1}-(-q^{h})^{-(n-t+1)}\)\kern3cm\\
\times
s_{\big((m-n-1)^n\big)}(q^n,q^{n-1},\dots,q,-1,q^{-1},\dots,
q^{-n+1},q^n).
\end{multline*}
Clearly, the determinant is easily evaluated by means of
\eqref{eq:sympl}. The specialized Schur function is evaluated in
Lemma~\ref{lem:Schur}. If everything is combined and simplified, the
claimed formulae \eqref{eq:C+-e-asy2a}--\eqref{eq:C+-e-asy2d} are
eventually obtained.
\end{proof}

We conclude this section by reporting the results from \cite{KrGVAB} on
the asymptotic behaviour of walks in the alcove $\mathcal A^{\tilde C_{n}}_m$ which
consist entirely of diagonal steps. 
These results have been stated there in an
equivalent form, namely in the language of walkers in the
lock-step vicious walkers model, which are bounded by two walls.

\begin{Theorem}[{\cite[Theorem~4]{KrGVAB}}] \label{thm:7-asy}
Let $m$ be a positive integer or half-integer. Furthermore,
let $\et=(\et_1,\et_2,\dots,\et_n)$ and 
$\la=(\la_1,\la_2,\dots,\la_n)$ be vectors of integers or of
half-integers in the
alcove $\mathcal A^{\tilde C_{n}}_m$ of type $\tilde C_n$ {\em(}defined in \eqref{eq:alcC}{\em)}.
Then, as $k$ tends to infinity such that $k\equiv 2\et_j+2\la_j$
mod\/~$2$, the number of random walks 
from $\et$ to $\la$ with exactly $k$ diagonal steps,
which stay in the alcove $\mathcal A^{\tilde C_{n}}_m$, is asymptotically
\begin{multline} \label{eq:vicious-with-2-asy}
 \frac {4^{n^2}} {(2m)^n} \bigg(2^n \prod _{j=1}
^{n}\cos\frac {j\pi} {2m}\bigg)^k
 \prod _{1\le h<t\le n} ^{}\(\sin\frac {\pi (\et_h-\et_t)} {2m} \cdot
   \sin\frac {\pi (\la_h-\la_t)} {2m}\)\\
\times
 \prod _{1\le h\le t\le n}\(\sin\frac {\pi (\et_h+\et_t)} {2m} \cdot
   \sin\frac {\pi (\la_h+\la_t)} {2m}\).
\end{multline}
\end{Theorem}

\begin{Theorem}[{\cite[Theorem~6]{KrGVAB}}] \label{thm:7-asy2}
Let $m$ be a positive integer or half-integer. Furthermore,
let $\et=(\et_1,\et_2,\dots,\et_n)$ be a vector of integers or of
half-integers in the
alcove $\mathcal A^{\tilde C_{n}}_m$ of type $\tilde C_n$ {\em(}defined in \eqref{eq:alcC}{\em)}.
Then, as $k$ tends to infinity, the number of random walks which
start at $\et$ and proceed for exactly $k$ diagonal steps,
which stay in the alcove $\mathcal A^{\tilde C_{n}}_m$, is asymptotically
\begin{multline} \label{eq:starvicious-with-2even}
\frac {4^{n^2}} {(2m)^n}  \bigg(2^n \prod _{j=1} ^{n}\cos \frac {j \pi}
{2m}\bigg)^k
 \prod _{1\le h<t\le n} ^{}\(\sin \frac {\pi  (\et_h-\et_t)} {2m} \cdot
    \sin \frac {\pi  (t-h)} {2m}\)\\
\times
 \prod _{1\le h\le t\le n} ^{}\(\sin \frac {\pi  (\et_h+\et_t)} {2m} \cdot
  \sin \frac {\pi  (t+h)} {2m}\)
 \prod _{h=0} ^{n}\prod _{t=1} ^{n}\frac {\sin \frac 
   {\pi  (t-h+\fl{m}\vphantom{\vrule depth3pt})} 
   {2m}} {\sin \frac {\pi  (t-h+n)} {2m}}\\
\times
 \prod _{h=1} ^{n}\frac {\sin \frac 
{\pi(h+\fl{m}-n\vphantom{\vrule depth3pt})} {2m}} 
   {\sin \frac {\pi(2h+\fl{m}-n\vphantom{\vrule depth3pt})} {2m}},
\end{multline}
if $k+2\et_j$ is odd, and
\begin{multline} \label{eq:starvicious-with-2odd}
\frac {4^{n^2}} {(2m)^n}  \bigg(2^n \prod _{j=1} ^{n}\cos \frac {j \pi}
{2m}\bigg)^k
 \prod _{1\le h<t\le n} ^{}\(\sin \frac {\pi  (\et_h-\et_t)} {2m} \cdot
    \sin \frac {\pi  (t-h)} {2m}\)\\
\times
 \prod _{1\le h\le t\le n} ^{}\(\sin \frac {\pi  (\et_h+\et_t)} {2m} 
\cdot   \sin \frac {\pi  (t+h)} {2m}\)
 \prod _{h=0} ^{n}\prod _{t=1} ^{n}\frac {\sin \frac 
   {\pi  (t-h+\cl{m-1}\vphantom{\vrule depth3pt})} 
   {2m}} {\sin \frac {\pi  (t-h+n)} {2m}},
\end{multline}
if $k+2\et_j$ is even.
\end{Theorem}

\section{Asymptotics for random walks in alcoves of type $\tilde B$} 
\label{sec:6}

This section is devoted to finding the asymptotic behaviour of the
number of walks from a given starting point to a given end point
which stay in the alcove $\mathcal A^{\tilde B_{n}}_m$ of type $\tilde B_n$ as the
number of steps becomes large, as well as the asymptotic behaviour
of the number of those walks which start at a given point but may
terminate anywhere.
In technical terms, we determine the asymptotic
behaviour of the expressions given by
Theorems~\ref{thm:8} and \ref{thm:9} as $k$ becomes large,
and as well if these expressions
are summed over all possible end points of the walks.

The following two theorems address the case of walks with standard steps.
The result for fixed starting {\it and\/} end point is the subject
of the first of the two.

\begin{Theorem} \label{thm:8-asy}
Let $m$ be a positive integer or half-integer. Furthermore,
let $\et=(\et_1,\et_2,\dots,\et_n)$ and 
$\la=(\la_1,\la_2,\dots,\la_n)$ be vectors of integers in the
alcove $\mathcal A^{\tilde B_{n}}_m$ of type $\tilde B_n$ 
{\em(}defined in \eqref{eq:alcB}{\em)}.
Then, as $k$ tends to infinity such that $k\equiv
\vert\et\vert+\vert\la\vert$ mod\/~$2$,
the number of random walks from $\et$ to $\la$ with exactly $k$
standard steps,
which stay in the alcove $\mathcal A^{\tilde B_{n}}_m$, is asymptotically
\begin{multline} \label{eq:B+-e-asy}
 \frac {4^{n^2}} {(2m)^n} \bigg(
\frac {\sin\frac {n\pi } {m}} {\sin\frac {\pi} {2m}}\bigg)^k
 \prod _{1\le h<t\le n} ^{}\bigg(\sin\frac {\pi (\et_h-\et_t)} {2m} \cdot
   \sin\frac {\pi (\la_h-\la_t)} {2m}\\
\cdot
\sin\frac {\pi (\et_h+\et_t)} {2m} \cdot
   \sin\frac {\pi (\la_h+\la_t)} {2m}\bigg)
 \prod _{h=1} ^{n}\(\sin\frac {\pi \et_h} {2m} \cdot
   \sin\frac {\pi \la_h} {2m}\).
\end{multline}
\end{Theorem}

\begin{proof} The analysis is analogous to the one in the proof of
Theorem~\ref{thm:6-asy}. Here, we have to estimate the  coefficient of
$x^k/k!$ in \eqref{eq:B+-e}. Expanding the two determinants 
in \eqref{eq:B+-e} by linearity in the rows, we obtain
\begin{multline*} 
\frac {1} {2m^n}\sum _{r_1,\dots,r_n=1} ^{2m-1}
\(2\sum _{j=1} ^{n}\cos\frac {\pi r_j} {m}\)^k
\(\sum _{j=1} ^{n}\sin\frac {\pi r_j\et_j} {m}\)
\det_{1\le h,t\le n}\(
\sin\frac {\pi r_h\la_t} {m}\)
\\+
\frac {1} {2m^n}\sum _{r_1,\dots,r_n=0} ^{2m-1}
\(2\sum _{j=1} ^{n}\cos\frac {\pi (2r_j+1)} {2m}\)^k\kern7cm\\
\cdot
\(
\prod _{j=1} ^{n}\sin\frac {\pi (2r_j+1)\et_j} {2m}\)
\det_{1\le h,t\le n}\(
\sin\frac {\pi (2r_h+1)\la_t} {2m}\).
\end{multline*}
Again, this is a finite sum of the form 
$\sum _{\ell} ^{}c_\ell b_\ell^k$, with both the $c_\ell$'s and $b_\ell$'s 
independent of $k$.

We have to locate the $b_\ell$'s
with largest modulus. As can be seen in a manner similar to the
considerations in the proof of Theorem~\ref{thm:6-asy}, the $b_\ell$'s with
largest modulus come from the expansion of the {\it second\/}
determinant, when we choose distinct $r_1,r_2,\dots,r_n$ from 
$$\{0,1,\dots,n-1,2m-n,2m-n+1,\dots,2m-1\},$$
or from
$$\{m-n,m-n+1,\dots,m-1,m,\dots,m+n-1\},$$
such that $r_h\ne 2m-1-r_t$ for all $h$ and $t$. Clearly, there are
$2\cdot 2^n$ such sets $\{r_1,r_2,\dots,r_n\}$. For each fixed set,
the sum of the corresponding terms $c_\ell b_\ell^k$ is equal to
\begin{equation} \label{eq:6-1}
\frac {1} {2m^n}\(2\sum _{j=1} ^{n}\cos\frac {\pi (2j-1)} {2m}\)^k
\det_{1\le h,t\le n}\(\sin\frac {\pi (2h-1)\et_t} {2m}\)
\det_{1\le h,t\le n}\(\sin\frac {\pi (2h-1)\la_t} {2m}\).
\end{equation}

We have on the one hand
$$\sum _{j=1} ^{n}\cos\frac {\pi (2j-1)} {2m}=\frac {\sin\frac {n\pi}
{m}} {2\sin\frac {\pi} {2m}}.$$
On the other hand, 
the first determinant in \eqref{eq:6-1} can be rewritten in the form
$$
(2i)^{-n}\det_{1\le h,t\le n}\(e^{\frac {\pi i(2h-1)\et_t} {2m}}
-e^{-\frac {\pi i(2h-1)\et_t} {2m}}\),
$$
and can thus be evaluated by means of \eqref{eq:ortho2}. After some
simplification, the result is
\begin{multline} \label{eq:sindet1}
\det_{1\le h,t\le n}\(\sin\frac {\pi (2h-1)\et_t} {2m}\)\\
=
2^{n^2-n} \prod _{1\le h<t\le n} ^{}\bigg(
   \sin\frac {\pi (\et_h-\et_t)} {2m}\cdot
   \sin\frac {\pi (\et_h+\et_t)} {2m}\bigg)
 \prod _{h=1} ^{n}\(
   \sin\frac {\pi \et_h} {2m}\).
\end{multline}
Clearly, the second determinant in \eqref{eq:6-1} is equal to the
same expression with $\et_j$ replaced by $\la_j$, $j=1,2,\dots,n$.
If all this is substituted in \eqref{eq:6-1}, and if the result is
multiplied by $2\cdot 2^n$, then we obtain exactly
\eqref{eq:B+-e-asy}.
\end{proof}

If the starting point is fixed but the end point is not, we have the
following result.

\begin{Theorem} \label{thm:8-asy2}
Let $m$ be a positive integer or half-integer. Furthermore,
let $\et=(\et_1,\et_2,\dots,\et_n)$ be a vector of integers in the
alcove $\mathcal A^{\tilde B_{n}}_m$ of type $\tilde B_n$ 
{\em(}defined in \eqref{eq:alcB}{\em)}.
Then, as $k$ tends to infinity,
the number of random walks which start at $\et$ and proceed for exactly $k$
standard steps,
which stay in the alcove $\mathcal A^{\tilde B_{n}}_m$, is asymptotically
\begin{multline} \label{eq:B+-e-asy2a}
 \frac {4^{n^2}} {2(2m)^n} \bigg(
\frac {\sin\frac {n\pi } {m}} {\sin\frac {\pi} {2m}}\bigg)^k
 \prod _{1\le h<t\le n} ^{}\bigg(\sin\frac {\pi (\et_h-\et_t)} {2m} \cdot
   \sin\frac {\pi (t-h)} {2m}\\
\cdot
\sin\frac {\pi (\et_h+\et_t)} {2m} \cdot
   \sin\frac {\pi (t+h)} {2m}\bigg)
 \prod _{h=1} ^{n}\(\sin\frac {\pi \et_h} {2m} \cdot
   \sin\frac {\pi h} {2m}\)\\
\times
\Bigg(
\prod _{h=1} ^{2n}\frac {\sin\frac {\pi(m-n+h)} {4m}}
{\sin\frac {\pi h} {4m}}
\prod _{h=1} ^{n-1}\prod _{t=1} ^{n}
\frac {\sin\frac {\pi(m+t-h)} {2m}}
{\sin\frac{\pi(n+t-h)} {2m}}
+\prod _{h=1} ^{2n}\frac {\sin\frac {\pi(m-n+h-1)} {4m}}
{\sin\frac {\pi h} {4m}}
\prod _{h=1} ^{n-1}\prod _{t=1} ^{n}
\frac {\sin\frac {\pi(m+t-h-1)} {2m}}
{\sin\frac{\pi(n+t-h)} {2m}}\\
+(-1)^{\vert\et\vert+k+\binom {n+1}2}
\prod _{h=1} ^{n}\frac {\cos\frac {\pi(m-n+2h-1)} {4m}\cdot
   \sin\frac {\pi(m-n+2h)} {4m}}
{\cos\frac {\pi (2h-1)} {4m}\cdot \sin\frac {\pi h} {2m}}
\prod _{h=1} ^{n-1}\prod _{t=1} ^{n}
\frac {\sin\frac {\pi(m+t-h)} {2m}}
{\sin\frac{\pi(n+t-h)} {2m}}\\
+(-1)^{\vert\et\vert+k+\binom {n}2}
\prod _{h=1} ^{n}\frac {\sin\frac {\pi(m-n+2h-2)} {4m}\cdot
   \cos\frac {\pi(m-n+2h-1)} {4m}}
{\cos\frac {\pi (2h-1)} {4m}\cdot \sin\frac {\pi h} {2m}}
\prod _{h=1} ^{n-1}\prod _{t=1} ^{n}
\frac {\sin\frac {\pi(m+t-h-1)} {2m}}
{\sin\frac{\pi(n+t-h)} {2m}}
\Bigg),
\end{multline}
if $m$ is an integer with parity equal to that of $n$,
it is asymptotically
\begin{multline} \label{eq:B+-e-asy2b}
 \frac {4^{n^2}} {2(2m)^n} \bigg(
\frac {\sin\frac {n\pi } {m}} {\sin\frac {\pi} {2m}}\bigg)^k
 \prod _{1\le h<t\le n} ^{}\bigg(\sin\frac {\pi (\et_h-\et_t)} {2m} \cdot
   \sin\frac {\pi (t-h)} {2m}\\
\cdot
\sin\frac {\pi (\et_h+\et_t)} {2m} \cdot
   \sin\frac {\pi (t+h)} {2m}\bigg)
 \prod _{h=1} ^{n}\(\sin\frac {\pi \et_h} {2m} \cdot
   \sin\frac {\pi h} {2m}\)\\
\times
\Bigg(
\prod _{h=1} ^{2n}\frac {\sin\frac {\pi(m-n+h)} {4m}}
{\sin\frac {\pi h} {4m}}
\prod _{h=1} ^{n-1}\prod _{t=1} ^{n}
\frac {\sin\frac {\pi(m+t-h)} {2m}}
{\sin\frac{\pi(n+t-h)} {2m}}
+\prod _{h=1} ^{2n}\frac {\sin\frac {\pi(m-n+h-1)} {4m}}
{\sin\frac {\pi h} {4m}}
\prod _{h=1} ^{n-1}\prod _{t=1} ^{n}
\frac {\sin\frac {\pi(m+t-h-1)} {2m}}
{\sin\frac{\pi(n+t-h)} {2m}}\\
+(-1)^{\vert\et\vert+k+\binom {n}2}
\prod _{h=1} ^{n}\frac {\sin\frac {\pi(m-n+2h-1)} {4m}\cdot
   \cos\frac {\pi(m-n+2h)} {4m}}
{\cos\frac {\pi (2h-1)} {4m}\cdot \sin\frac {\pi h} {2m}}
\prod _{h=1} ^{n-1}\prod _{t=1} ^{n}
\frac {\sin\frac {\pi(m+t-h)} {2m}}
{\sin\frac{\pi(n+t-h)} {2m}}\\
+(-1)^{\vert\et\vert+k+\binom {n+1}2}
\prod _{h=1} ^{n}\frac {\cos\frac {\pi(m-n+2h-2)} {4m}\cdot
   \sin\frac {\pi(m-n+2h-1)} {4m}}
{\cos\frac {\pi (2h-1)} {4m}\cdot \sin\frac {\pi h} {2m}}
\prod _{h=1} ^{n-1}\prod _{t=1} ^{n}
\frac {\sin\frac {\pi(m+t-h-1)} {2m}}
{\sin\frac{\pi(n+t-h)} {2m}}
\Bigg),
\end{multline}
if $m$ is an integer with parity different from that of $n$,
and it is asymptotically 
\begin{multline} \label{eq:B+-e-asy2c}
 \frac {4^{n^2}} {(2m)^n} \bigg(
\frac {\sin\frac {n\pi } {m}} {\sin\frac {\pi} {2m}}\bigg)^k
 \prod _{1\le h<t\le n} ^{}\bigg(\sin\frac {\pi (\et_h-\et_t)} {2m} \cdot
   \sin\frac {\pi (t-h)} {2m}\\
\cdot
\sin\frac {\pi (\et_h+\et_t)} {2m} \cdot
   \sin\frac {\pi (t+h)} {2m}\bigg)
 \prod _{h=1} ^{n}\(\sin\frac {\pi \et_h} {2m} \cdot
   \sin\frac {\pi h} {2m}\)\\
\times
\prod _{h=1} ^{2n}\frac {\sin\frac {\pi(\fl{m}-n+h)} {4m}}
{\sin\frac {\pi h} {4m}}
\prod _{h=1} ^{n-1}\prod _{t=1} ^{n}
\frac {\sin\frac {\pi(\fl{m}+t-h)} {2m}}
{\sin\frac{\pi(n+t-h)} {2m}},
\end{multline}
if $m$ is a half-integer.
\end{Theorem}

\begin{proof}
Here, in view of Theorem~\ref{thm:8-asy}, we have to carry out the
sum of \eqref{eq:B+-e-asy} over all possible $\la$, i.e., over all
$\la_1>\la_2>\dots>\la_n>0$ with $\la_1+\la_2<2m$, where $\vert\la\vert\equiv
k+\vert\et\vert$ mod~$2$. Leaving away factors which are independent of
$\la$, the problem is to compute the sum
\begin{equation} \label{eq:SB1}
{\underset{\la_1+\la_2<2m}
{\sum _{\la_1>\dots>\la_n>0} ^{}}\kern-.5cm}{}^{\displaystyle\prime}
\kern.8cm
 \prod _{1\le h<t\le n} ^{}\bigg( \sin\frac {\pi (\la_h-\la_t)} {2m}
\cdot \sin\frac {\pi (\la_h+\la_t)} {2m}\bigg)
 \prod _{h=1} ^{n}   \sin\frac {\pi \la_h} {2m},
\end{equation}
where the sum $\sum{}^{\textstyle\prime}$ 
is either restricted to those $\la$ for which $\vert\la\vert$
is even, or to those for which $\vert\la\vert$ is odd, depending on
the parity of $k+\vert\et\vert$. 

Let first $m$ be an integer.
In order to get rid of the restriction $\la_1+\la_2<2m$ in the sum in
\eqref{eq:SB1}, we observe that the summand remains unchanged if we
replace $\la_1$ by $2m-\la_1$ and, moreover, that
$(\la_1,\la_2,\dots,\la_n)$ satisfies the conditions on the summation
indices if and only if $(2m-\la_1,\la_2,\dots,\la_n)$ does. Hence, we
may rewrite the sum in \eqref{eq:SB1} as
\begin{multline} \label{eq:SB}
{{\sum _{m\ge\la_1>\dots>\la_n>0} ^{}}\kern-.8cm}{}^{\displaystyle\prime}
\kern.8cm
 \prod _{1\le h<t\le n} ^{}\bigg( \sin\frac {\pi (\la_h-\la_t)} {2m}
\cdot \sin\frac {\pi (\la_h+\la_t)} {2m}\bigg)
 \prod _{h=1} ^{n}   \sin\frac {\pi \la_h} {2m}\\
+{{\sum _{m-1\ge\la_1>\dots>\la_n>0} ^{}}\kern-1cm}{}^{\displaystyle\prime}
\kern1cm
 \prod _{1\le h<t\le n} ^{}\bigg( \sin\frac {\pi (\la_h-\la_t)} {2m}
\cdot \sin\frac {\pi (\la_h+\la_t)} {2m}\bigg)
 \prod _{h=1} ^{n}   \sin\frac {\pi \la_h} {2m},
\end{multline}
where, again, the sums $\sum{}^{\textstyle\prime}$ 
are either restricted to those $\la$ for which $\vert\la\vert$
is even, or to those for which $\vert\la\vert$ is odd, depending on
the parity of $k+\vert\et\vert$.

The task of evaluating the sums in \eqref{eq:SB} will be accomplished if we
are able to evaluate the (complete) sum
\begin{equation*} 
U_1(c)=\sum _{c\ge\la_1>\dots>\la_n>0} ^{}
 \prod _{1\le h<t\le n} ^{}\bigg( \sin\frac {\pi (\la_h-\la_t)} {2m}
\cdot \sin\frac {\pi (\la_h+\la_t)} {2m}\bigg)
 \prod _{h=1} ^{n}   \sin\frac {\pi \la_h} {2m},
\end{equation*}
and its ``signed variant''
\begin{equation*} 
U_2(c)=\sum _{c\ge\la_1>\dots>\la_n>0} ^{}
(-1)^{\vert\la\vert}
 \prod _{1\le h<t\le n} ^{}\bigg( \sin\frac {\pi (\la_h-\la_t)} {2m}
\cdot \sin\frac {\pi (\la_h+\la_t)} {2m}\bigg)
 \prod _{h=1} ^{n}   \sin\frac {\pi \la_h} {2m},
\end{equation*}
for $c=m$ and $c=m-1$.
The expression \eqref{eq:SB} is then equal to 
$$\frac {1} {2}\Big(U_1(m)+U_1(m-1)+U_2(m)+U_2(m-1)\Big)$$ 
if the sum $\sum{}^{\textstyle\prime}$ is 
restricted to the $\la$'s for which $\vert\la\vert$ is even, and it is
equal to 
$$\frac {1} {2}\Big(U_1(m)+U_1(m-1)-U_2(m)-U_2(m-1)\Big)$$ 
if
the sum $\sum{}^{\textstyle\prime}$  is
restricted to the $\la$'s for which $\vert\la\vert$ is odd.

In order to evaluate $U_1(c)$ and $U_2(c)$, 
we use \eqref{eq:sindet1} to rewrite them uniformly as 
\begin{align*} \label{}
\sum _{c\ge\la_1>\dots>\la_n>0} ^{}&
\ep^{\vert\la\vert}
2^{-n^2+n}
\det_{1\le h,t\le n}\big(\sin(\pi (2h-1)\la_t/2m)\big)\\
&=\sum _{c\ge\la_1>\dots>\la_n>0} ^{}
\frac {1} {2^{n^2}i^n}
\det_{1\le h,t\le n}\big((\ep e^{\pi i (2h-1)/2m})^{\la_t}-(\ep e^{\pi i
(2h-1)/2m})^{-\la_t}\big),
\end{align*}
with $\ep=1$ or $\ep=-1$ depending on whether we want to express
$U_1(c)$ or $U_2(c)$.
Replacing $\la_t$ by $\la_t+n-t+1$, $t=1,2,\dots,n$, we obtain the
expression
\begin{equation} \label{eq:SummeB}
\sum _{c-n\ge\la_1\ge\dots\ge\la_n\ge0} ^{}
\frac {1} {2^{n^2}i^n}
\det_{1\le h,t\le n}\big((\ep e^{\pi i (2h-1)/2m})^{\la_t+n-t+1}-(\ep e^{\pi i
(2h-1)/2m})^{-(\la_t+n-t+1)}\big).
\end{equation}
As in the proof of Theorem~\ref{thm:6-asy2},
this determinant can be expressed in terms of a
{\it symplectic character} as defined in \eqref{e24}.
Specifically, writing $q$ for $e^{\pi i/2m}$, the sum in
\eqref{eq:SummeB} equals
\begin{equation} \label{eq:sum-symplB}
\frac {1} {2^{n^2}i^n}
\det_{1\le h,t\le n}
\( (\ep q^{2h-1})^{n-t+1}-(\ep q^{2h-1})^{-(n-t+1)}\)
\sum _{c-n\ge\la_1\ge\dots\ge\la_n\ge0}
\sp_\la(\ep q,\ep q^3,\dots,\ep q^{2n-1}).
\end{equation}
Use of Formula~\eqref{e31} 
in \eqref{eq:sum-symplB} yields the equivalent expression
\begin{multline} \label{eq:Form1}
\frac {1} {2^{n^2}i^n}
\det_{1\le h,t\le n}
\( (\ep q^{2h-1})^{n-t+1}-(\ep q^{2h-1})^{-(n-t+1)}\)\\
\times
s_{\big((c-n)^n\big)}(\ep q^{2n-1},\ep q^{2n-3},\dots,\ep q^3,\ep
q,1,\ep q^{-1},\ep q^{-3},\dots,
\ep q^{-2n+3},\ep q^{-2n+1})\\
=\ep^{(c-n)n}\frac {1} {2^{n^2}i^n}
\det_{1\le h,t\le n}
\( (\ep q^{2h-1})^{n-t+1}-(\ep q^{2h-1})^{-(n-t+1)}\)\kern3cm\\
\times
s_{\big((c-n)^n\big)}(q^{2n-1},q^{2n-3},\dots,q^3,q,\ep,q^{-1},q^{-3},\dots,
q^{-2n+3},q^{-2n+1}).
\end{multline}
Clearly, the determinant is easily evaluated by means of
\eqref{eq:ortho2}. The specialized Schur function is evaluated in
Lemma~\ref{lem:SchurA} in the case that $\ep=1$, and in 
Lemma~\ref{lem:SchurB} in the case that $\ep=-1$. 
If everything is combined and simplified, the
claimed formulae \eqref{eq:B+-e-asy2a} and \eqref{eq:B+-e-asy2b} are
eventually obtained.

In the case that $m$ is a half-integer, an adaption of the above argument
of replacement of $\la_1$ by $2m-\la_1$ in \eqref{eq:SB1} shows that
\eqref{eq:SB1} equals
\begin{multline*} 
{{\sum _{\fl{m}\ge\la_1>\dots>\la_n>0} ^{}}\kern-.9cm}{}^{\displaystyle\prime}
\kern.8cm
 \prod _{1\le h<t\le n} ^{}\bigg( \sin\frac {\pi (\la_h-\la_t)} {2m}
\cdot \sin\frac {\pi (\la_h+\la_t)} {2m}\bigg)
 \prod _{h=1} ^{n}   \sin\frac {\pi \la_h} {2m}\\
+{{\sum _{\fl{m}\ge\la_1>\dots>\la_n>0}
^{}}\kern-.9cm}{}^{\displaystyle\prime\prime}
\kern.8cm
 \prod _{1\le h<t\le n} ^{}\bigg( \sin\frac {\pi (\la_h-\la_t)} {2m}
\cdot \sin\frac {\pi (\la_h+\la_t)} {2m}\bigg)
 \prod _{h=1} ^{n}   \sin\frac {\pi \la_h} {2m},
\end{multline*}
where one of the sums $\sum{}^{\textstyle\prime}$ and
$\sum{}^{\textstyle\prime\prime}$
is restricted to those $\la$ for which $\vert\la\vert$ is even
while the other to those for which $\vert\la\vert$ is odd, depending on
the parity of $k+\vert\et\vert$.
Hence, it equals $U_1(\fl{m})$, regardless of the parity of
$k+\vert\et\vert$. On
the other hand, we have seen above that $U_1(\fl{m})$ can be
rewritten as \eqref{eq:Form1} with $\ep=1$ and $c=\fl{m}$.
Application of Lemma~\ref{lem:SchurA} with $c=\fl{m}-n$ and some
simplification then lead to \eqref{eq:B+-e-asy2c}.
\end{proof}

We proceed with the corresponding theorems for the case of walks with
diagonal steps.

\begin{Theorem} \label{thm:9-asy}
Let $m$ be a positive integer or half-integer. Furthermore,
let $\et=(\et_1,\et_2,\break\dots,\et_n)$ and 
$\la=(\la_1,\la_2,\dots,\la_n)$ be vectors of integers or of
half-integers in the
alcove $\mathcal A^{\tilde B_{n}}_m$ of type $\tilde B_n$ 
{\em(}defined in \eqref{eq:alcB}{\em)}.
Then, as $k$ tends to infinity such that $k\equiv 2\et_j+2\la_j$ mod\/~$2$, 
the number of random walks 
from $\et$ to $\la$ with exactly $k$ diagonal steps,
which stay in the alcove $\mathcal A^{\tilde B_{n}}_m$, is asymptotically
\begin{multline} \label{eq:Bd-asy}
 \frac {4^{n^2}} {2(2m)^n} \bigg(2^n \prod _{j=1}
^{n}\cos\frac {(2j-1)\pi} {4m}\bigg)^k
 \prod _{1\le h<t\le n} ^{}\bigg(\sin\frac {\pi (\et_h-\et_t)} {2m} \cdot
   \sin\frac {\pi (\la_h-\la_t)} {2m}\\
\cdot
\sin\frac {\pi (\et_h+\et_t)} {2m} \cdot
   \sin\frac {\pi (\la_h+\la_t)} {2m}\bigg)
 \prod _{h=1} ^{n}\(\sin\frac {\pi \et_h} {2m} \cdot
   \sin\frac {\pi \la_h} {2m}\).
\end{multline}
\end{Theorem}

\begin{proof}We proceed in the same way as in the proof of
Theorem~\ref{thm:8-asy}.
This time we have to estimate the two determinants 
in \eqref{eq:Bd}. As a preparatory step, we replace the two sums over
$r$ in the determinants by 2 times the same sums, but restricted to $r$ from
$0$ to $2m-1$, which can be safely done when $\et$, $\la$ and $k$ are
chosen so that paths from $\et$ to $\la$ in $k$ steps exist. (We
already did an analogous transformation in the proof of
Theorem~\ref{thm:4-asy}).
Then, by expanding the two determinants, we obtain 
\begin{multline*} 
\frac {1} {2m^n}\sum _{r_1,\dots,r_n=1} ^{2m-1}
\(2^n\prod _{j=1} ^{n}\cos\frac {\pi r_j}
{2m}\)^k
\(
\prod _{j=1} ^{n}\sin\frac {\pi r_j\et_j} {m}\)\det_{1\le h,t\le n}\(
\sin\frac {\pi r_h\la_t} {m}\)\\
+\frac {1} {2m^n}\sum _{r_1,\dots,r_n=0} ^{2m-1}
\(2^n\prod _{j=1} ^{n}\cos\frac {\pi (2r_j+1)} {4m}\)^k\kern7cm\\
\cdot
\(
\prod _{j=1} ^{n}\sin\frac {\pi (2r_j+1)\et_j} {2m}\)
\det_{1\le h,t\le n}\(\sin\frac {\pi (2r_h+1)\la_t} {2m}
\),
\end{multline*}
which is again a finite sum of the form $\sum _{\ell} ^{}c_\ell
b_\ell^k$, with both the $c_\ell$'s and $b_\ell$'s 
independent of $k$. 

As it turns out, also here the dominating terms 
(the terms for which $b_\ell$ has
largest modulus) come from the second determinant. 
More precisely,
these are the terms corresponding to the subsets $\{r_1,r_2,\dots,r_n\}$ 
of 
$$\{0,1,\dots,n,2m-n,2m-n+1,\dots,2m-1\},$$
with the property that all $r_h$'s are distinct and $r_h\ne 2m-1-r_t$
for all $h$ and $t$. Clearly, there are
$2^n$ such sets $\{r_1,r_2,\dots,r_n\}$. For each fixed set,
the sum of the corresponding terms $c_\ell b_\ell^k$ is equal to
\begin{equation} \label{eq:6-2}
\frac {1} {2m^n}\(2^n\prod _{j=1} ^{n}\cos\frac {\pi (2j-1)} {4m}\)^k
\det_{1\le h,t\le n}\(\sin\frac {\pi (2h-1)\la_t} {2m}\)
\det_{1\le h,t\le n}\(\sin\frac {\pi (2h-1)\et_t} {2m}\).
\end{equation}
The two determinants in this expression have already been evaluated
in \eqref{eq:sindet1}. If this is substituted in \eqref{eq:6-2}, and
if the result is multiplied by $2^n$, then we obtain exactly
\eqref{eq:Bd-asy}.
\end{proof}

\begin{Theorem} \label{thm:9-asy2}
Let $m$ be a positive integer or half-integer. Furthermore,
let $\et=(\et_1,\et_2,\break\dots,\et_n)$ be a vector of integers or of
half-integers in the
alcove $\mathcal A^{\tilde B_{n}}_m$ of type $\tilde B_n$ 
{\em(}defined in \eqref{eq:alcB}{\em)}.
Then, as $k$ tends to infinity 
the number of random walks which start at $\et$ and proceed for exactly $k$
diagonal steps,
which stay in the alcove $\mathcal A^{\tilde B_{n}}_m$, is asymptotically
\begin{multline} \label{eq:Bd-asy2a}
 \frac {4^{n^2}} {2(2m)^n} \bigg(2^n \prod _{j=1}
^{n}\cos\frac {(2j-1)\pi} {4m}\bigg)^k
 \prod _{1\le h<t\le n} ^{}\bigg(\sin\frac {\pi (\et_h-\et_t)} {2m} \cdot
   \sin\frac {\pi (t-h)} {2m}\\
\cdot
\sin\frac {\pi (\et_h+\et_t)} {2m} \cdot
   \sin\frac {\pi (t+h)} {2m}\bigg)
 \prod _{h=1} ^{n}\(\sin\frac {\pi \et_h} {2m} \cdot
   \sin\frac {\pi h} {2m}\)\\
\times
\Bigg(
\prod _{h=1} ^{2n}\frac {\sin\frac {\pi(m-n+h)} {4m}}
{\sin\frac {\pi h} {4m}}
\prod _{h=1} ^{n-1}\prod _{t=1} ^{n}
\frac {\sin\frac {\pi(m+t-h)} {2m}}
{\sin\frac{\pi(n+t-h)} {2m}}
+\prod _{h=1} ^{2n}\frac {\sin\frac {\pi(m-n+h-1)} {4m}}
{\sin\frac {\pi h} {4m}}
\prod _{h=1} ^{n-1}\prod _{t=1} ^{n}
\frac {\sin\frac {\pi(m+t-h-1)} {2m}}
{\sin\frac{\pi(n+t-h)} {2m}}
\Bigg),
\end{multline}
if $m$ is an integer and $k+2\et_j$ is even, it is asymptotically
\begin{multline} \label{eq:Bd-asy2b}
 \frac {4^{n^2}} {(2m)^n} \bigg(2^n \prod _{j=1}
^{n}\cos\frac {(2j-1)\pi} {4m}\bigg)^k
 \prod _{1\le h<t\le n} ^{}\bigg(\sin\frac {\pi (\et_h-\et_t)} {2m} \cdot
   \sin\frac {\pi (t-h)} {2m}\\
\cdot
\sin\frac {\pi (\et_h+\et_t)} {2m} \cdot
   \sin\frac {\pi (t+h)} {2m}\bigg)
 \prod _{h=1} ^{n}\(\sin\frac {\pi \et_h} {2m} \cdot
   \sin\frac {\pi h} {2m}\)\\
\times
\prod _{h=1} ^{2n}\frac {\sin\frac {\pi(\fl m-n+h)} {4m}}
{\sin\frac {\pi h} {4m}}
\prod _{h=1} ^{n-1}\prod _{t=1} ^{n}
\frac {\sin\frac {\pi(\fl m+t-h)} {2m}}
{\sin\frac{\pi(n+t-h)} {2m}},
\end{multline}
if $m$ is a half-integer and $k+2\et_j$ is even, it is asymptotically
\begin{multline} \label{eq:Bd-asy2c}
 \frac {4^{n^2}} {(2m)^n} \bigg(2^n \prod _{j=1}
^{n}\cos\frac {(2j-1)\pi} {4m}\bigg)^k
 \prod _{1\le h<t\le n} ^{}\bigg(\sin\frac {\pi (\et_h-\et_t)} {2m} \cdot
   \sin\frac {\pi (t-h)} {2m}\\
\cdot
\sin\frac {\pi (\et_h+\et_t)} {2m} \cdot
   \sin\frac {\pi (t+h-1)} {2m}\bigg)
 \prod _{h=1} ^{n}\(\sin\frac {\pi \et_h} {2m} \cdot
   \sin\frac {\pi (2h-1)} {4m}\)\\
\times
\prod _{h=1} ^{2n}\frac {\sin{\frac {\pi(m-n+h)} {4m}}} 
{\sin{\frac {\pi h} {4m}}}
\prod _{h=1} ^{n}\prod _{t=1} ^{n}
\frac {\sin\frac {\pi(m+t-h)} {2m}} 
{\sin\frac {\pi(n+t-h)} {2m}}
\prod _{h=1} ^{n}
\frac {\sin\frac {\pi h} {2m}}
{\sin\frac {\pi(m-n+2h)} {4m}}
\frac {\cos\frac {\pi(2h-1)} {4m}}
{\cos\frac {\pi(m-n+2h-1)} {4m}},
\end{multline}
if $m$ is an integer and $k+2\et_j$ is odd, and it is asymptotically
\begin{multline} \label{eq:Bd-asy2d}
 \frac {4^{n^2}} {2(2m)^n} \bigg(2^n \prod _{j=1}
^{n}\cos\frac {(2j-1)\pi} {4m}\bigg)^k
 \prod _{1\le h<t\le n} ^{}\bigg(\sin\frac {\pi (\et_h-\et_t)} {2m} \cdot
   \sin\frac {\pi (t-h)} {2m}\\
\cdot
\sin\frac {\pi (\et_h+\et_t)} {2m} \cdot
   \sin\frac {\pi (t+h-1)} {2m}\bigg)
 \prod _{h=1} ^{n}\(\sin\frac {\pi \et_h} {2m} \cdot
   \sin\frac {\pi (2h-1)} {4m}\)\\
\times
\Bigg(
\prod _{h=1} ^{2n}\frac {\sin{\frac {\pi(\cl{m}-n+h)} {4m}}} 
{\sin{\frac {\pi h} {4m}}}
\prod _{h=1} ^{n}\prod _{t=1} ^{n}
\frac {\sin\frac {\pi(\cl{m}+t-h)} {2m}} 
{\sin\frac {\pi(n+t-h)} {2m}}
\prod _{h=1} ^{n}
\frac {\sin\frac {\pi h} {2m}}
{\sin\frac {\pi(\cl{m}-n+2h)} {4m}}
\frac {\cos\frac {\pi(2h-1)} {4m}}
{\cos\frac {\pi(\cl{m}-n+2h-1)} {4m}}\\
+
\prod _{h=1} ^{2n}\frac {\sin{\frac {\pi(\fl{m}-n+h)} {4m}}} 
{\sin{\frac {\pi h} {4m}}}
\prod _{h=1} ^{n}\prod _{t=1} ^{n}
\frac {\sin\frac {\pi(\fl{m}+t-h)} {2m}} 
{\sin\frac {\pi(n+t-h)} {2m}}
\prod _{h=1} ^{n}
\frac {\sin\frac {\pi h} {2m}}
{\sin\frac {\pi(\fl{m}-n+2h)} {4m}}
\frac {\cos\frac {\pi(2h-1)} {4m}}
{\cos\frac {\pi(\fl{m}-n+2h-1)} {4m}}
\Bigg),
\end{multline}
if $m$ is a half-integer and $k+2\et_j$ is odd.
\end{Theorem}
\begin{proof}
Here, in view of Theorem~\ref{thm:9-asy}, we have to carry out the
sum of \eqref{eq:Bd-asy} over all possible $\la$, i.e., over all
$\la_1>\la_2>\dots>\la_n>0$ with $\la_1+\la_2<2m$, where 
$k\equiv 2\et_j+2\la_j$ mod\/~$2$. 
Leaving away factors which are independent of
$\la$, the problem is to compute the sum
\begin{equation} \label{eq:SB2}
{\underset{\la_1+\la_2<2m}
{\sum _{\la_1>\dots>\la_n>0} ^{}}\kern-.5cm}{}^{\displaystyle\prime}
\kern.8cm
 \prod _{1\le h<t\le n} ^{}\bigg( \sin\frac {\pi (\la_h-\la_t)} {2m}
\cdot \sin\frac {\pi (\la_h+\la_t)} {2m}\bigg)
 \prod _{h=1} ^{n}   \sin\frac {\pi \la_h} {2m},
\end{equation}
where the sum $\sum{}^{\textstyle\prime}$ 
is restricted to integral, respectively to half-integral
$\la_1,\la_2,\dots,\la_n$, depending on whether $k+2\et_j$ is even or odd. 

As in previous proofs, we have to distinguish between several cases.
Let us first restrict the sum \eqref{eq:SB2} to integral
$\la_1,\la_2,\dots,\la_n$. Then, by repeating the argument from the
proof of Theorem~\ref{thm:8-asy2} of replacement of $\la_1$ by
$2m-\la_1$, and using notation from that proof, 
we obtain that the sum \eqref{eq:SB2} equals
$U_1(m)+U_1(m-1)$ if $m$ is an integer, and it equals
$2U_1(\fl m)$ if $m$ is a half-integer. 
In the proof of Theorem~\ref{thm:8-asy2}
it was shown how to evaluate the sum $U_1(c)$ by means of
Lemma~\ref{lem:SchurA}. If the result is substituted, we arrive at
\eqref{eq:Bd-asy2a} and \eqref{eq:Bd-asy2b}.

If we restrict the sum \eqref{eq:SB2} to half-integral
$\la_1,\la_2,\dots,\la_n$, then we are faced with the problem of
evaluating the sum
\begin{equation*} 
V(c)=\sum _{c\ge\la_1>\dots>\la_n>0} ^{}
\kern-.8cm{}^{\displaystyle\prime}
\kern.8cm
 \prod _{1\le h<t\le n} ^{}\bigg( \sin\frac {\pi (\la_h-\la_t)} {2m}
\cdot \sin\frac {\pi (\la_h+\la_t)} {2m}\bigg)
 \prod _{h=1} ^{n}   \sin\frac {\pi \la_h} {2m},
\end{equation*}
where the sum $\sum{}^{\textstyle\prime}$ is
over all {\it half-integral\/} $\la_1,\la_2,\dots,\la_n$.
More specifically, using the replacement of $\la_1$ by $2m-\la_1$
another time, one sees readily that
the sum \eqref{eq:SB2} equals $2V(m-\frac {1} {2})$
if $m$ is an integer, and it equals $V(m)+V(m-1)$ if $m$ is a
half-integer.

In order to evaluate $V(c)$, where $c$ is a half-integer, 
we use \eqref{eq:sindet1} to rewrite it as
\begin{align*} 
\sum _{c\ge\la_1>\dots>\la_n>0} ^{}
\kern-.8cm{}^{\displaystyle\prime}
\kern.8cm
&
2^{-n^2+n}
\det_{1\le h,t\le n}\big(\sin(\pi (2h-1)\la_t/2m)\big)\\
&=\sum _{c\ge\la_1>\dots>\la_n>0} ^{}
\kern-.8cm{}^{\displaystyle\prime}
\kern.8cm
\frac {1} {2^{n^2}i^n}
\det_{1\le h,t\le n}\big((e^{\pi i (2h-1)/2m})^{\la_t}-(e^{\pi i
(2h-1)/2m})^{-\la_t}\big).
\end{align*}
Replacing $\la_t$ by $\la_t+n-t+\frac {1} {2}$, $t=1,2,\dots,n$, we obtain the
expression
\begin{equation} \label{eq:SummeB2}
\sum _{\cl{c-n}\ge\la_1\ge\dots\ge\la_n\ge0} ^{}
\frac {1} {2^{n^2}i^n}
\det_{1\le h,t\le n}\big((e^{\pi i (2h-1)/2m})^{\la_t+n-t+1/2}-(e^{\pi i
(2h-1)/2m})^{-(\la_t+n-t+1/2)}\big).
\end{equation}
(Note that the last sum is over {\it integral} 
$\la_1,\la_2,\dots,\la_n$.)
This determinant can be expressed in terms of an {\it odd orthogonal
character}. (See \eqref{e12} for the definition.)
Writing again $q$ for $e^{\pi i/2m}$, the sum in \eqref{eq:SummeB2}
equals
\begin{equation} \label{eq:sum-ortho}
\frac {1} {2^{n^2}i^n}
\det_{1\le h,t\le n}\big((q^{2h-1})^{n-t+\frac
{1} {2}}-(q^{2h-1})^{-(n-t+\frac {1} {2})}\big)
\sum _{\cl{c-n}\ge\la_1\ge\dots\ge\la_n\ge0}
\so^{odd}_\la(q^{2n-1},q^{2n-3},\dots,q).
\end{equation}
Again there is a formula which allows us to evaluate the sum in the last
line (see \cite[(3.2)]{KratBC}),
\begin{equation*} 
s_{((a^{n-p},(a-1)^p)}(x_1,x_1^{-1},\dots,x_p,x_p^{-1},1)=\underset
{\oddrows\!\big((a^r)/\nu\big)=p}
{\sum _{a\ge\nu_1\ge\dots\ge\nu_n\ge0} ^{}}\so^{odd}_\nu(x_1,\dots,x_p),
\end{equation*}
where the notation $(a^{n-p},(a-1)^p\big)$ is a short
notation for the vector in which the first $n-p$ components are
$a$, the next $p$ components are $a-1$, 
followed by $n+1$ components all of which are $0$, and
where $\oddrows\!\big((a^n)/\nu\big)=p$ means that the number of
rows of odd length in the skew shape $(a^n)/\nu$ equals exactly
$p$.
Use of this formula in \eqref{eq:sum-ortho} gives
\begin{multline*} 
\frac {1} {2^{n^2}i^n}
\det_{1\le h,t\le n}\big((q^{2h-1})^{n-t+\frac
{1} {2}}-(q^{2h-1})^{-(n-t+\frac {1} {2})}\big)\\
\times
\sum _{p=0} ^{n}
s_{\big(\cl{c-n}^{n-p},\cl{c-n-1}^p\big)}
(q^{2n-1},q^{2n-3},\dots,q^3,q,1,q^{-1},q^{-3},
\dots,q^{-2n+3},q^{-2n+1}).
\end{multline*}
We now use
Lemma~\ref{lem:SchurA} to evaluate the Schur function in the last
line to obtain the expression
\begin{multline*} 
\frac {1} {2^{n^2}i^n}
\det_{1\le h,t\le n}\big((q^{2h-1})^{n-t+\frac
{1} {2}}-(q^{2h-1})^{-(n-t+\frac {1} {2})}\big)\\
\times
\prod _{h=1} ^{2n}\frac {\(q^{\frac {\cl{c}-n+h} {2}}-q^{-\frac 
{\cl{c}-n+h} {2}}\)} 
{\(q^{\frac {h} {2}}-q^{-\frac {h} {2}}\)}
\prod _{h=1} ^{n}\prod _{t=1} ^{n}
\frac {\(q^{\cl{c}+t-h}-q^{-\cl{c}-t+h}\)} 
{\(q^{n+t-h}-q^{-n-t+h}\)}
\prod _{h=1} ^{n}
\(q^{h}-q^{-h}\)^2
\\
\times
\sum _{p=0} ^{n}
\prod _{h=1} ^{n}
\frac {1} 
{\(q^{\cl{c}-n+p+h}-q^{-\cl{c}+n-p-h}\)}
\frac {1} {\displaystyle\prod _{h=1} ^{p}\(q^h-q^{-h}\)
\prod _{h=1} ^{n-p}\(q^h-q^{-h}\)}\\
\cdot
\frac {(q^{\frac {\cl{c}-n} {2}}-q^{-\frac {\cl{c}-n} {2}})
(q^{\frac {\cl{c}-n} {2}+p}+q^{-\frac {\cl{c}-n} {2}-p})} 
{(q^{\cl{c}-n+p}-q^{-\cl{c}+n-p})}.
\end{multline*}
In terms of the standard basic hypergeometric notation
\begin{equation} \label{eq:qhyp}
{}_r\phi_s\!\left[\begin{matrix} a_1,\dots,a_r\\ b_1,\dots,b_s\end{matrix}; q,
z\right]=\sum _{\ell=0} ^{\infty}\frac {\poq{a_1}{\ell}\cdots\poq{a_r}{\ell}}
{\poq{q}{\ell}\poq{b_1}{\ell}\cdots\poq{b_s}{\ell}}\left((-1)^\ell q^{\binom
\ell2}\right)^{s-r+1}z^\ell\ ,
\end{equation}
where the shifted $q$-factorials $(a;q)_\ell$ are defined by
$(a;q)_\ell:=(1-a)(1-aq)\cdots(1-aq^{\ell-1})$, $\ell\ge1$,
$(a;q)_0:=1$, this can be written in the form
\begin{multline} \label{eq:sum-SchurA}
\frac {1} {2^{n^2}i^n}
\det_{1\le h,t\le n}\big((q^{2h-1})^{n-t+\frac
{1} {2}}-(q^{2h-1})^{-(n-t+\frac {1} {2})}\big)\\
\times
\prod _{h=1} ^{2n}\frac {\(q^{\frac {\cl{c}-n+h} {2}}-q^{-\frac 
{\cl{c}-n+h} {2}}\)} 
{\(q^{\frac {h} {2}}-q^{-\frac {h} {2}}\)}
\prod _{h=1} ^{n}\prod _{t=1} ^{n}
\frac {\(q^{\cl{c}+t-h}-q^{-\cl{c}-t+h}\)} 
{\(q^{n+t-h}-q^{-n-t+h}\)}
\prod _{h=1} ^{n}
\(q^{h}-q^{-h}\)^2
\\
\times
{\frac {{q^{\cl{c} n + n }} }
   {({\let \over / \def\frac#1#2{#1 / #2} {q^2}}; {q^2}) _{n} \,
     ({\let \over / \def\frac#1#2{#1 / #2} {q^{2\cl{c}-2n+2}}}; {q^2}) _{n}
}}
{} _{3} \phi _{2} \! \left [             \begin{matrix}
      \let \over / \def\frac#1#2{#1 / #2} {q^{2 \cl{c}-2n}},
-{q^{\cl{c}-n+2}},
      {q^{-2 n}}\\ \let \over / \def\frac#1#2{#1 / #2} -{q^{\cl{c}-n}}, {q^{2
\cl{c} +
      2 }}\end{matrix} ;{q^2}, {\displaystyle -{q^{2n+1}}} \right ].
\end{multline}
The determinant is again easily evaluated by means of
\eqref{eq:ortho2}. On the other hand, 
the $_3\phi_2$-series in the last line can be evaluated by a limit
case of Jackson's very-well-poised $_8\phi_7$-summation
(see \cite[(2.6.2); Appendix (II.22)]{GaRaAA})
\begin{multline} \label{eq:Jackson}
{}_8\phi _7\!\left [ \begin{matrix} \let\over/ A,{\sqrt{A}} q,- {\sqrt{A}} q 
   ,B,C,D,{A^2 q^{1+N} \over B C D },{q^{-N}}\\ \let\over/  {\sqrt{A}},-{\sqrt{A}},
   {{A q}\over B},{{A q}\over
   C},{{A q}\over D},{{B C D}\over {A {q^N}}},A {q^{1 + N}}\end{matrix}
   ;q,q\right ] \\= {\frac {{(\let\over/ A q;q)}_{N} {(\let\over/ {{A q}\over {B C}};q)}_{N} 
      {(\let\over/ {{A q}\over {B D}};q)}_{N} {(\let\over/ {{A q}\over {C D}};q)}_{N}} 
    {{(\let\over/ {{A q}\over B};q)}_{N} {(\let\over/ {{A q}\over C};q)}_{N} {(\let\over/ {{A q}\over D};q)}_{N} 
      {(\let\over/ {{A q}\over {B C D}};q)}_{N}}},
\end{multline}
where $N$ is a nonnegative integer. Namely, if in \eqref{eq:Jackson} we
let $N\to\infty$, put $D=\sqrt{A}$ and $C=-\sqrt{A}$, and finally
replace $q$ by $q^2$, we are left
with
\begin{equation} \label{eq:Jackson2}
{} _{3} \phi _{2} \! \left [             \begin{matrix} \let \over / \def\frac#1#2{#1
   / #2} A, - {\sqrt{A}} {q^2} , B\\ \let \over /
   \def\frac#1#2{#1 / #2} -{\sqrt{A}}, {{A {q^2}}\over B}\end{matrix} ;{q^2},
   {\displaystyle -{\frac q B}} \right ] = 
  {\frac {({\let \over / \def\frac#1#2{#1 / #2} -q}; {q^2}) _{\infty} \,
      ({\let \over / \def\frac#1#2{#1 / #2} -{{{\sqrt{A}} q}\over B}}; {q^2})
       _{\infty} \,({\let \over / \def\frac#1#2{#1 / #2} A {q^2}}; {q^2})
       _{\infty} \,({\let \over / \def\frac#1#2{#1 / #2}
       {{{\sqrt{A}} {q^2}}\over B}}; {q^2}) _{\infty} } 
    {({\let \over / \def\frac#1#2{#1 / #2} - {\sqrt{A}} q  };
       {q^2}) _{\infty} \,({\let \over / \def\frac#1#2{#1 / #2} -{q\over B}};
       {q^2}) _{\infty} \,({\let \over / \def\frac#1#2{#1 / #2}
       {\sqrt{A}} {q^2}}; {q^2}) _{\infty} \,
      ({\let \over / \def\frac#1#2{#1 / #2} {{A {q^2}}\over B}}; {q^2})
       _{\infty} }}.
\end{equation}
The $_3\phi_2$-series in \eqref{eq:sum-SchurA} is a special case of
the above $_3\phi_2$-series in which $A={q^{2 \cl{c}-2n}}$ and
$B=q^{-2n}$. If we substitute the corresponding right-hand side of
\eqref{eq:Jackson2} in \eqref{eq:sum-SchurA}, we have evaluated the
sum $V(c)$. This, in its turn, leads to the expressions
\eqref{eq:Bd-asy2c} and \eqref{eq:Bd-asy2d}.
\end{proof}

\section{Asymptotics for random walks in alcoves of type $\tilde D$} 
\label{sec:7}

In this final section we turn our attention to the walks in the 
alcove $\mathcal
A^{\tilde D_{n}}_m$ of type $\tilde D_n$. We
determine the asymptotic behaviour of the
number of walks from a given starting point to a given end point
which stay in $\mathcal A^{\tilde D_{n}}_m$ as the
number of steps becomes large, as well as the asymptotic behaviour
of the number of those walks which start at a given point but may
terminate anywhere.
In technical terms, we determine the asymptotic
behaviour of the expressions given by
Theorems~\ref{thm:10} and \ref{thm:11} as $k$ becomes large,
and as well if these expressions
are summed over all possible end points of the walks.

The first two theorems in this section give our results for the 
case of walks with standard steps.

\begin{Theorem} \label{thm:10-asy}
Let $m$ be a positive integer. Furthermore,
let $\et=(\et_1,\et_2,\dots,\et_n)$ and 
$\la=(\la_1,\la_2,\dots,\la_n)$ be vectors of integers in the
alcove $\mathcal A^{\tilde D_{n}}_m$ of type $\tilde D_n$ 
{\em(}defined in \eqref{eq:alcD}{\em)}.
Then, as $k$ tends to infinity such that $k\equiv
\vert\et\vert+\vert\la\vert$ mod\/~$2$,
the number of random walks from $\et$ to $\la$ with exactly $k$
standard steps,
which stay in the alcove $\mathcal A^{\tilde D_{n}}_m$, is asymptotically
\begin{multline} \label{eq:D+-e-asy}
 \frac {4^{n^2}} {(8m)^n} \bigg(2
\frac {\sin\frac {n\pi } {2m}\cos\frac {(n-1)\pi} {2m}} 
{\sin\frac {\pi} {2m}}\bigg)^k\\
\times
 \prod _{1\le h<t\le n} ^{}\bigg(\sin\frac {\pi (\et_h-\et_t)} {2m} \cdot
   \sin\frac {\pi (\la_h-\la_t)} {2m}\cdot
\sin\frac {\pi (\et_h+\et_t)} {2m} \cdot
   \sin\frac {\pi (\la_h+\la_t)} {2m}\bigg).
\end{multline}
\end{Theorem}

\begin{proof}We proceed as in the proof of Theorem~\ref{thm:8-asy}.
Since we have done similar calculations already several times, we
shall be brief here.

We have to find the asymptotic behaviour of \eqref{eq:D+-e} (instead
of the rather similar \eqref{eq:B+-e}, with which we dealt in
the proof of Theorem~\ref{thm:8-asy}). By applying arguments very
similar to those in the proof of Theorem~\ref{thm:8-asy}, we infer
that the dominating terms in the expansion of \eqref{eq:D+-e} come
from the expansion of the third determinant. To be precise, these
dominating terms add up to $2\cdot 2^{n-1}$ times
\begin{equation} \label{eq:6-3}
\frac {1} {4m^n}
\(2\sum _{j=0} ^{n-1}\cos\frac {\pi j} {m}\)^k
\det_{1\le h,t\le n}\(
\cos\frac {\pi (h-1)\et_t} {m}\)
\det_{1\le h,t\le n}\(
\cos\frac {\pi (h-1)\la_t} {m}\).
\end{equation}
We have
$$\sum _{j=0} ^{n-1}\cos\frac {\pi j} {m}=
\frac {\sin\frac {n\pi } {2m}\cdot\cos\frac {(n-1)\pi} {2m}} 
{\sin\frac {\pi} {2m}}.
$$
On the other hand, the first determinant in \eqref{eq:6-3} can be
rewritten in the form
$$2^{-n}\det_{1\le h,t\le n}\(
e^{\frac {\pi i (h-1)\et_t} {m}}+e^{-\frac {\pi i (h-1)\et_t} {m}}\),
$$
and can thus be evaluated by means of \eqref{eq:ortho3}. The result
is that 
\begin{equation} \label{eq:cosdet}
\det_{1\le h,t\le n}\(
\cos\frac {\pi (h-1)\et_t} {m}\)
=2^{n^2-2n+1}
 \prod _{1\le h<t\le n} ^{}\bigg(
   \sin\frac {\pi (\et_h-\et_t)} {2m}\cdot
   \sin\frac {\pi (\et_h+\et_t)} {2m}\bigg).
\end{equation}
If this is substituted in \eqref{eq:6-3} (and as well the analogous
evaluation of the second determinant in \eqref{eq:6-3}), and
if the result is multiplied by $2\cdot 2^{n-1}=2^n$, then we obtain exactly
\eqref{eq:D+-e-asy}.
\end{proof}

If the starting point is fixed but the end point is not, we have the
following result.

{\allowdisplaybreaks
\begin{Theorem} \label{thm:10-asy2}
Let $m$ be a positive integer or half-integer. Furthermore,
let $\et=(\et_1,\et_2,\dots,\et_n)$ be a vector of integers in the
alcove $\mathcal A^{\tilde D_{n}}_m$ of type $\tilde D_n$ 
{\em(}defined in \eqref{eq:alcD}{\em)}.
Then, as $k$ tends to infinity,
the number of random walks which start at $\et$ and proceed for exactly $k$
standard steps,
which stay in the alcove $\mathcal A^{\tilde D_{n}}_m$, is asymptotically
\begin{multline} \label{eq:D+-e-asy2a}
 \frac {4^{n^2}} {2^n(8m)^n} \bigg(2
\frac {\sin\frac {n\pi } {2m}\cos\frac {(n-1)\pi} {2m}} 
{\sin\frac {\pi} {2m}}\bigg)^k\\
\times
 \prod _{1\le h<t\le n} ^{}\bigg(
\sin\frac {\pi (\et_h-\et_t)} {2m} \cdot
\sin\frac {\pi (t-h)} {2m} \cdot
\sin\frac {\pi (\et_h+\et_t)} {2m} \cdot
\sin\frac {\pi (t+h-2)} {2m} \bigg)\\
\times
\left(\frac {1} {2^{n-1}}
     \frac{\prod_{h=1}^{n}\sin\frac {\pi (m-n+2h)} {2m}}
     {\prod_{h=1}^{n-1}\sin\frac {\pi h} {2m}}
 \prod _{1\le h<t\le n} ^{}
\frac {\sin\frac {\pi (m-n+h+t-1)} {2m}\cdot \sin\frac {\pi (m-n+t+h)} {2m}} 
{\sin^2\frac {\pi (t+h-2)} {2m}}
\right.
\\
\cdot
     \sum_{k=1}^{n}\frac {(-1)^{n-k}(m-n+2k)}
        {\prod_{h=1}^{n} 
\sin\frac {\pi (m-n+k+h)} {2m}
         \prod_{h=1}^{k-1}\sin\frac {\pi h} {2m}
         \prod_{h=1}^{n-k}\sin\frac {\pi h} {2m}}\\
+\frac {1} {2^{n-1}}
     \frac{\prod_{h=1}^{n}\sin\frac {\pi (m-n+2h-1)} {2m}}
     {\prod_{h=1}^{n-1}\sin\frac {\pi h} {2m}}
 \prod _{1\le h<t\le n} ^{}
\frac {\sin\frac {\pi (m-n+h+t-2)} {2m}\cdot \sin\frac {\pi (m-n+t+h-1)} {2m}} 
{\sin^2\frac {\pi (t+h-2)} {2m}}
\\
\cdot
     \sum_{k=1}^{n}\frac {(-1)^{n-k}(m-n+2k-2)}
        {\prod_{h=1}^{n} 
\sin\frac {\pi (m-n+k+h-1)} {2m}
         \prod_{h=1}^{k-1}\sin\frac {\pi h} {2m}
         \prod_{h=1}^{n-k}\sin\frac {\pi h} {2m}}\\
+\frac {1} {2^{n-1}}
     \frac{\prod_{h=1}^{n}\sin\frac {\pi (m-n+2h-1)} {2m}}
     {\prod_{h=1}^{n-1}\sin\frac {\pi h} {2m}}
 \prod _{1\le h<t\le n} ^{}
\frac {\sin\frac {\pi (m-n+h+t)} {2m}\cdot \sin\frac {\pi (m-n+t+h-1)} {2m}} 
{\sin^2\frac {\pi (t+h-2)} {2m}}
\\
\cdot
     \sum_{k=1}^{n}\frac {(-1)^{n-k}(m-n+2k-1)}
        {\prod_{h=1}^{n} 
\sin\frac {\pi (m-n+k+h-1)} {2m}
         \prod_{h=1}^{k-1}\sin\frac {\pi h} {2m}
         \prod_{h=1}^{n-k}\sin\frac {\pi h} {2m}}\\
+\frac {1} {2^{n-1}}
     \frac{\prod_{h=1}^{n}\sin\frac {\pi (m-n+2h-2)} {2m}}
     {\prod_{h=1}^{n-1}\sin\frac {\pi h} {2m}}
 \prod _{1\le h<t\le n} ^{}
\frac {\sin\frac {\pi (m-n+h+t-1)} {2m}\cdot \sin\frac {\pi (m-n+t+h-2)} {2m}} 
{\sin^2\frac {\pi (t+h-2)} {2m}}
\\
\cdot
     \sum_{k=1}^{n}\frac {(-1)^{n-k}(m-n+2k-3)}
        {\prod_{h=1}^{n} 
\sin\frac {\pi (m-n+k+h-2)} {2m}
         \prod_{h=1}^{k-1}\sin\frac {\pi h} {2m}
         \prod_{h=1}^{n-k}\sin\frac {\pi h} {2m}}\\
+(-1)^{mn+k+\vert\et\vert}
\prod _{h=1} ^{n-1}\frac {1} {\cos\frac {\pi (h-1)} {2m}}
 \prod _{1\le h<t\le n} ^{}
\frac {\sin\frac {\pi (m-n+h+t-1)} {2m}\cdot \sin\frac {\pi (m-n+t+h)} {2m}} 
{\sin^2\frac {\pi (t+h-2)} {2m}}\\
+(-1)^{mn+k+\vert\et\vert}
\prod _{h=1} ^{n-1}\frac {1} {\cos\frac {\pi (h-1)} {2m}}
 \prod _{1\le h<t\le n} ^{}
\frac {\sin\frac {\pi (m-n+h+t-2)} {2m}\cdot \sin\frac {\pi (m-n+t+h-1)} {2m}} 
{\sin^2\frac {\pi (t+h-2)} {2m}}\\
+(-1)^{mn+k+\vert\et\vert}
\prod _{h=1} ^{n-1}\frac {1} {\cos\frac {\pi (h-1)} {2m}}
 \prod _{1\le h<t\le n} ^{}
\frac {\sin\frac {\pi (m-n+h+t)} {2m}\cdot \sin\frac {\pi (m-n+t+h-1)} {2m}} 
{\sin^2\frac {\pi (t+h-2)} {2m}}\\
\left.
+(-1)^{mn+k+\vert\et\vert}
\prod _{h=1} ^{n-1}\frac {1} {\cos\frac {\pi (h-1)} {2m}}
 \prod _{1\le h<t\le n} ^{}
\frac {\sin\frac {\pi (m-n+h+t-1)} {2m}\cdot \sin\frac {\pi (m-n+t+h-2)} {2m}} 
{\sin^2\frac {\pi (t+h-2)} {2m}}
\right)
\end{multline}
if $m$ is an integer, and it is asymptotically
\begin{multline} \label{eq:D+-e-asy2b}
 \frac {4^{n^2-n+1}} {(8m)^n} \bigg(2
\frac {\sin\frac {n\pi } {2m}\cos\frac {(n-1)\pi} {2m}} 
{\sin\frac {\pi} {2m}}\bigg)^k\\
\times
 \prod _{1\le h<t\le n} ^{}\bigg(
\sin\frac {\pi (\et_h-\et_t)} {2m} \cdot
\sin\frac {\pi (t-h)} {2m} \cdot
\sin\frac {\pi (\et_h+\et_t)} {2m} \cdot
\sin\frac {\pi (t+h-2)} {2m} \bigg)\\
\left(
     \frac{\prod_{h=1}^{n}\sin\frac {\pi (\fl{m}-n+2h)} {2m}}
     {\prod_{h=1}^{n-1}\sin\frac {\pi h} {2m}}
 \prod _{1\le h<t\le n} ^{}
\frac {\sin\frac {\pi (\fl{m}-n+h+t-1)} {2m}\cdot \sin\frac {\pi
(\fl{m}-n+t+h)} {2m}} 
{\sin^2\frac {\pi (t+h-2)} {2m}}
\right.
\\
\cdot
     \sum_{k=1}^{n}\frac {(-1)^{n-k}(\fl{m}-n+2k)}
        {\prod_{h=1}^{n} 
\sin\frac {\pi (\fl{m}-n+k+h)} {2m}
         \prod_{h=1}^{k-1}\sin\frac {\pi h} {2m}
         \prod_{h=1}^{n-k}\sin\frac {\pi h} {2m}}\\
+
     \frac{\prod_{h=1}^{n}\sin\frac {\pi (\fl{m}-n+2h-1)} {2m}}
     {\prod_{h=1}^{n-1}\sin\frac {\pi h} {2m}}
 \prod _{1\le h<t\le n} ^{}
\frac {\sin\frac {\pi (\fl{m}-n+h+t)} {2m}\cdot \sin\frac {\pi
(\fl{m}-n+t+h-1)} {2m}} 
{\sin^2\frac {\pi (t+h-2)} {2m}}
\\
\left.\cdot
     \sum_{k=1}^{n}\frac {(-1)^{n-k}(\fl{m}-n+2k-1)}
        {\prod_{h=1}^{n} 
\sin\frac {\pi (\fl{m}-n+k+h-1)} {2m}
         \prod_{h=1}^{k-1}\sin\frac {\pi h} {2m}
         \prod_{h=1}^{n-k}\sin\frac {\pi h} {2m}}
\right)
\end{multline}
if $m$ is a half-integer.
\end{Theorem}
}

\begin{proof}
Here, in view of Theorem~\ref{thm:10-asy}, we have to carry out the
sum of \eqref{eq:D+-e-asy} over all possible $\la$, i.e., over all
$\la_1>\la_2>\dots>\la_{n-1}>\vert\la_n\vert$ with $\la_1+\la_2<2m$, 
where $\vert\la\vert\equiv
k+\vert\et\vert$ mod~$2$. Leaving away factors which are independent of
$\la$, the problem is to compute the sum
\begin{equation} \label{eq:SD1}
{\underset{\la_1+\la_2<2m}
{\sum _{\la_1>\la_2>\dots>\la_{n-1}>\vert\la_n\vert} 
^{}}\kern-1.2cm}{}^{\displaystyle\prime}
\kern1.3cm
 \prod _{1\le h<t\le n} ^{}\bigg(
   \sin\frac {\pi (\la_h-\la_t)} {2m}\cdot
   \sin\frac {\pi (\la_h+\la_t)} {2m}\bigg),
\end{equation}
where the sum $\sum{}^{\textstyle\prime}$ 
is either restricted to those $\la$ for which $\vert\la\vert$
is even, or to those for which $\vert\la\vert$ is odd, depending on
the parity of $k+\vert\et\vert$. 

Let first $m$ be an integer.
As in the proof of Theorem~\ref{thm:8-asy2}, instead of a sum
$\sum{}^{\textstyle\prime}$ where we have to deal with the
unwieldy constraints $\la_1+\la_2<2m$ and $\la_{n-1}>\vert\la_n\vert$, 
we would rather prefer sums where the only constraint is of the form 
$\la_1>\la_2>\dots>\la_n$. In order to get rid of the constraint
$\la_1+\la_2<2m$, we apply again the argument 
of replacement of $\la_1$ by $2m-\la_1$ that we already used in the
proofs of Theorems~\ref{thm:8-asy2} and \ref{thm:9-asy2}.
Similarly, in order to get rid of the constraint
$\la_{n-1}>\vert\la_n\vert$, we observe that the summand in
\eqref{eq:SD1} remains invariant under the replacement of $\la_n$ by
$-\la_n$, and that $(\la_1,\la_2,\dots,\la_n)$ satisfies the conditions on the summation
indices if and only if $(\la_1,\la_2,\dots,-\la_n)$ does.
Hence, we may rewrite the sum in \eqref{eq:SD1} as
\begin{multline} \label{eq:SD}
{{\sum _{m\ge\la_1>\dots>\la_n\ge0} ^{}}\kern-.8cm}{}^{\displaystyle\prime}
\kern.8cm
 \prod _{1\le h<t\le n} ^{}\bigg(
   \sin\frac {\pi (\la_h-\la_t)} {2m}\cdot
   \sin\frac {\pi (\la_h+\la_t)} {2m}\bigg)\\
+{{\sum _{m-1\ge\la_1>\dots>\la_n\ge0} ^{}}\kern-.99cm}{}^{\displaystyle\prime}
\kern1cm
 \prod _{1\le h<t\le n} ^{}\bigg(
   \sin\frac {\pi (\la_h-\la_t)} {2m}\cdot
   \sin\frac {\pi (\la_h+\la_t)} {2m}\bigg)\\
+{{\sum _{m\ge\la_1>\dots>\la_n>0} ^{}}\kern-.8cm}{}^{\displaystyle\prime}
\kern.8cm
 \prod _{1\le h<t\le n} ^{}\bigg(
   \sin\frac {\pi (\la_h-\la_t)} {2m}\cdot
   \sin\frac {\pi (\la_h+\la_t)} {2m}\bigg)\\
+{{\sum _{m-1\ge\la_1>\dots>\la_n>0} ^{}}\kern-1cm}{}^{\displaystyle\prime}
\kern1cm
 \prod _{1\le h<t\le n} ^{}\bigg(
   \sin\frac {\pi (\la_h-\la_t)} {2m}\cdot
   \sin\frac {\pi (\la_h+\la_t)} {2m}\bigg),
\end{multline}
where, again, the sums $\sum{}^{\textstyle\prime}$ 
are either restricted to those $\la$ for which $\vert\la\vert$
is even, or to those for which $\vert\la\vert$ is odd, depending on
the parity of $k+\vert\et\vert$.

The task of evaluating the sums in \eqref{eq:SD} will be accomplished if we
are able to evaluate the (complete) sums
\begin{equation*} 
W_1(c)=\sum _{c\ge\la_1>\dots>\la_n\ge0} ^{}
 \prod _{1\le h<t\le n} ^{}\bigg(
   \sin\frac {\pi (\la_h-\la_t)} {2m}\cdot
   \sin\frac {\pi (\la_h+\la_t)} {2m}\bigg)
\end{equation*}
and
\begin{equation*} 
W_2(c)=\sum _{c\ge\la_1>\dots>\la_n>0} ^{}
 \prod _{1\le h<t\le n} ^{}\bigg(
   \sin\frac {\pi (\la_h-\la_t)} {2m}\cdot
   \sin\frac {\pi (\la_h+\la_t)} {2m}\bigg),
\end{equation*}
and their ``signed variants''
\begin{equation*} 
W_3(c)=\sum _{c\ge\la_1>\dots>\la_n\ge0} ^{}
(-1)^{\vert\la\vert}
 \prod _{1\le h<t\le n} ^{}\bigg(
   \sin\frac {\pi (\la_h-\la_t)} {2m}\cdot
   \sin\frac {\pi (\la_h+\la_t)} {2m}\bigg)
\end{equation*}
and
\begin{equation*} 
W_4(c)=\sum _{c\ge\la_1>\dots>\la_n>0} ^{}
(-1)^{\vert\la\vert}
 \prod _{1\le h<t\le n} ^{}\bigg(
   \sin\frac {\pi (\la_h-\la_t)} {2m}\cdot
   \sin\frac {\pi (\la_h+\la_t)} {2m}\bigg)
\end{equation*}
for $c=m$ and $c=m-1$.
The expression \eqref{eq:SD} is then equal to 
$$\frac {1} {2}\Big(W_1(m)+W_1(m-1)+W_2(m)+W_2(m-1)+
W_3(m)+W_3(m-1)+W_4(m)+W_4(m-1)\Big)$$ 
if the sum $\sum{}^{\textstyle\prime}$ is 
restricted to the $\la$'s for which $\vert\la\vert$ is even, and it is
equal to 
$$\frac {1} {2}\Big(W_1(m)+W_1(m-1)+W_2(m)+W_2(m-1)-
W_3(m)-W_3(m-1)-W_4(m)-W_4(m-1)\Big)$$ 
if the sum $\sum{}^{\textstyle\prime}$  is
restricted to the $\la$'s for which $\vert\la\vert$ is odd.

We first show how to evaluate $W_1(c)$ and $W_3(c)$. Using
\eqref{eq:cosdet}, we rewrite them uniformly as 
\begin{align*} \label{}
\sum _{c\ge\la_1>\dots>\la_n\ge0} ^{}&
\ep^{\vert\la\vert}
2^{-n^2+2n-1}
\det_{1\le h,t\le n}\big(\cos(\pi (h-1)\la_t/m)\big)\\
&=\sum _{c\ge\la_1>\dots>\la_n\ge0} ^{}
\frac {1} {2^{n^2-n+1}}
\det_{1\le h,t\le n}\big((\ep e^{\pi i (h-1)/m})^{\la_t}+(\ep e^{\pi i
(h-1)/m})^{-\la_t}\big),
\end{align*}
with $\ep=1$ or $\ep=-1$ depending on whether we want to express
$W_1(c)$ or $W_3(c)$.
Replacing $\la_t$ by $\la_t+n-t$, $t=1,2,\dots,n$, we obtain the
expression
\begin{equation} \label{eq:SummeD}
\sum _{c-n+1\ge\la_1\ge\dots\ge\la_n\ge0} ^{}
\frac {1} {2^{n^2-n+1}}
\det_{1\le h,t\le n}\big((\ep e^{\pi i (h-1)/m})^{\la_t+n-t}+(\ep e^{\pi i
(h-1)/m})^{-(\la_t+n-t)}\big).
\end{equation}
The determinant in the summand appears as a part in the definition of an
{\it even orthogonal character} given in \eqref{e25}. 
Specifically, if we write $q$ for $e^{\pi i/m}$, the sum in
\eqref{eq:SummeD} equals
\begin{multline} \label{eq:sum-symplD}
\frac {1} {2^{n^2-n+1}}
\det_{1\le h,t\le n}\big((\ep q^{h-1})^{n-t}+(\ep
q^{h-1})^{-(n-t)}\big)\\
\times
\sum _{c-n+1\ge\la_1\ge\dots\ge\la_n\ge0} ^{}
\so^{even}_\la(\ep q^{n-1},\ep q^{n-2},\dots,\ep q,\ep).
\end{multline}
(The reader should note that the second determinant in \eqref{e25} vanishes
if one of the variables $x_t$, $t=1,2,\dots,n$, is $1$ or $-1$, 
which is the case in our situation.) 
The sum in \eqref{eq:sum-symplD} can be simplified by means of
another character identity from \cite{KratBC}: given non-negative integers or
half-integers $a,b$ with $a\ge b$, the identity \cite[(3.15)]{KratBC} 
implies
\begin{multline} \label{eq:soeven}
\so^{even}_{(a^n)}(x_1,x_2,\dots,x_n)\cdot
\sum _{p=-b} ^{b}\so^{even}_{(b^{n-1},p)}(x_1,x_2,\dots,x_n)\\=
{\sum _{a+b\ge\nu_1\ge\nu_2\ge\dots\ge\nu_n\ge a-b} ^{}}
\so^{even}_{\nu}(x_1,x_2,\dots,x_n),
\end{multline}
with the understanding that the sum on the left-hand side ranges
over half-integral $p$ if $b$ is a half-integer, 
and that the sum on the right-hand side ranges
over integers $\nu_t$ if $a+b$ is an integer, and over half-integers
$\nu_t$ if $a+b$ is a half-integer.
The notation $(a^n)$ is short for the
vector consisting of $n$ components all of which equal to $a$, while
the notation $(b^{n-1},p)$ means the vector in which the first $n-1$
components are $b$, followed by a component $p$. If we use
\eqref{eq:soeven}, then \eqref{eq:sum-symplD} becomes
\begin{multline*} 
\frac {1} {2^{n^2-n+1}}
\det_{1\le h,t\le n}\big((\ep q^{h-1})^{n-t}+(\ep q^{h-1})^{-(n-t)}\big)
\so^{even}_{\big(((c-n+1)/2)^n\big)}(\ep q^{n-1},\ep q^{n-2},\dots,\ep q,\ep)\\
\times
\sum _{p=-(c-n+1)/2}
^{(c-n+1)/2}\so^{even}_{\big(((c-n+1)/2)^{n-1},p\big)}(\ep q^{n-1},\ep q^{n-2},\dots,\ep
q,\ep)\\
=\ep^{(c-n+1)n}\frac {1} {2^{n^2-n+1}}
\det_{1\le h,t\le n}\big((\ep q^{h-1})^{n-t}+(\ep
q^{h-1})^{-(n-t)}\big)\kern5cm\\
\times
\so^{even}_{\big(((c-n+1)/2)^n\big)}( q^{n-1}, q^{n-2},\dots, q,1)\\
\times
\sum _{p=-(c-n+1)/2} ^{(c-n+1)/2}\ep^{(c-n+1-2p)/2}
\so^{even}_{\big(((c-n+1)/2)^{n-1},p\big)}( q^{n-1}, q^{n-2},\dots,
q,1).
\end{multline*}
Clearly, the first determinant can be evaluated by means of
\eqref{eq:ortho3}. The even orthogonal character of shape $\big(((c-n+1)/2)^n\big)$
can be evaluated by using the definition \eqref{e25} and
\eqref{eq:ortho3}. (It should be observed that the second
determinant in the numerator on the right-hand side of \eqref{e25}
vanishes if $x_n=\pm1$.) Finally, the sum over $p$ of even orthogonal characters is
evaluated in Lemma~\ref{lem:spin} if $\ep=1$, respectively in
Lemma~\ref{lem:spin-1} if $\ep=-1$.

In order to evaluate $W_2(c)$ and $W_4(c)$, we proceed in a completely
analogous fashion. In particular, the two sums $W_2(c)$ and $W_4(c)$
are equal to
\begin{equation} \label{eq:SummeD2}
\sum _{c-n+1\ge\la_1\ge\dots\ge\la_n>0} ^{}
\frac {1} {2^{n^2-n+1}}
\det_{1\le h,t\le n}\big((\ep e^{\pi i (h-1)/m})^{\la_t+n-t}+(\ep e^{\pi i
(h-1)/m})^{-(\la_t+n-t)}\big),
\end{equation}
with $\ep=1$ or $\ep=-1$ depending on whether we want to express
$W_2(c)$ or $W_4(c)$. By the use of \eqref{e31}, 
this sum becomes
\begin{multline} \label{eq:sum-symplD2}
\frac {1} {2^{n^2-n+1}}
\det_{1\le h,t\le n}\big((\ep q^{h-1})^{n-t}+(\ep q^{h-1})^{-(n-t)}\big)
\sum _{c-n+1\ge\la_1\ge\dots\ge\la_n\ge1} ^{}
\so^{even}_\la(\ep ,\ep q,\dots,\ep q^{n-1}).
\end{multline}
The only difference to \eqref{eq:sum-symplD} is that in the 
sum $\la_n$ must be {\it positive}, instead of
just {\it non-negative}. Hence, the remaining steps are the same,
with the small modification that we use \eqref{eq:soeven} with 
$a=(c-n+2)/2$ and $b=(c-n)/2$ to see that
\eqref{eq:sum-symplD2} equals
\begin{multline*} 
\ep^{(c-n+1)n}\frac {1} {2^{n^2-n+1}}
\det_{1\le h,t\le n}\big((\ep q^{h-1})^{n-t}+(\ep
q^{h-1})^{-(n-t)}\big)\kern5cm\\
\times
\so^{even}_{\big(((c-n+2)/2)^n\big)}( q^{n-1}, q^{n-2},\dots, q,1)\\
\times
\sum _{p=-(c-n)/2} ^{(c-n)/2}\ep^{(c-n-2p)/2}
\so^{even}_{\big(((c-n)/2)^{n-1},p\big)}( q^{n-1}, q^{n-2},\dots,
q,1).
\end{multline*}
Again, the first determinant can be evaluated by means of
\eqref{eq:ortho3}, the even orthogonal character of shape $\big(((c-n+2)/2)^n\big)$
can be evaluated by using the definition \eqref{e25} and
\eqref{eq:ortho3}, while the sum over $p$ of even orthogonal characters is
evaluated in Lemma~\ref{lem:spin}, respectively in
Lemma~\ref{lem:spin-1}.

If everything is combined and simplified, the
claimed formula \eqref{eq:D+-e-asy2a} is eventually obtained.

In the case that $m$ is a half-integer, an adaption of the above
arguments shows that \eqref{eq:SD1} equals
\begin{multline*} 
{{\sum _{\fl{m}\ge\la_1>\dots>\la_n\ge0} ^{}}}
 \prod _{1\le h<t\le n} ^{}\bigg(
   \sin\frac {\pi (\la_h-\la_t)} {2m}\cdot
   \sin\frac {\pi (\la_h+\la_t)} {2m}\bigg)\\
+{{\sum _{\fl{m}\ge\la_1>\dots>\la_n>0}
^{}}}
 \prod _{1\le h<t\le n} ^{}\bigg(
   \sin\frac {\pi (\la_h-\la_t)} {2m}\cdot
   \sin\frac {\pi (\la_h+\la_t)} {2m}\bigg),
\end{multline*}
regardless of the parity of $k+\vert\et\vert$.
Hence, it equals $W_1(\fl{m})+W_2(\fl{m})$.
On the other hand, we have seen above how to evaluate $W_1(\fl{m})$ 
and $W_2(\fl{m})$. If this is substituted and the resulting
expressions are simplified, we eventually obtain 
\eqref{eq:D+-e-asy2b}.
\end{proof}

We conclude the section with the analogous results for the case of
walks with diagonal steps.

\begin{Theorem} \label{thm:11-asy}
Let $m$ be a positive integer or half-integer. Furthermore,
let $\et=(\et_1,\et_2,\break\dots,\et_n)$ and 
$\la=(\la_1,\la_2,\dots,\la_n)$ be vectors of integers or of
half-integers in the
alcove $\mathcal A^{\tilde D_{n}}_m$ of type $\tilde D_n$ {\em(}defined in \eqref{eq:alcD}{\em)}.
Then, as $k$ tends to infinity such that $k\equiv 2\et_j+2\la_j$ mod\/~$2$,
the number of random walks 
from $\et$ to $\la$ with exactly $k$ diagonal steps,
which stay in the alcove $\mathcal A^{\tilde D_{n}}_m$, is asymptotically
\begin{multline} \label{eq:Dd-asy}
 \frac {4^{n^2}} {2(8m)^n} \bigg(2^n \prod _{j=0}
^{n-1}\cos\frac {j\pi} {2m}\bigg)^k\\
\times
 \prod _{1\le h<t\le n} ^{}\bigg(\sin\frac {\pi (\et_h-\et_t)} {2m} \cdot
   \sin\frac {\pi (\la_h-\la_t)} {2m}
\sin\frac {\pi (\et_h+\et_t)} {2m} \cdot
   \sin\frac {\pi (\la_h+\la_t)} {2m}\bigg).
\end{multline}
\end{Theorem}

\begin{proof}We proceed as in the proof of
Theorem~\ref{thm:9-asy}. We shall again be brief here.

We have to find the asymptotic behaviour of \eqref{eq:Dd} (instead
of the rather similar \eqref{eq:Bd}, with which we dealt in
the proof of Theorem~\ref{thm:9-asy}). By applying arguments very
similar to those in the proof of Theorem~\ref{thm:9-asy}, we infer
that the dominating terms in the expansion of \eqref{eq:Dd} come
from the expansion of the third determinant. To be precise, these
dominating terms add up to $2^{n-1}$ times
\begin{equation} \label{eq:6-4}
\frac {1} {4m^n}
\(2^n\prod _{j=0} ^{n-1}\cos\frac {\pi j} {2m}\)^k
\det_{1\le h,t\le n}\(
\cos\frac {\pi (h-1)\la_t} {m}\)
\det_{1\le h,t\le n}\(
\cos\frac {\pi (h-1)\et_t} {m}\).
\end{equation}
We have already evaluated these determinants in \eqref{eq:cosdet}.
If we use these evaluations in \eqref{eq:6-4}, and
if the result is multiplied by $2^{n-1}$, then we obtain exactly
\eqref{eq:Dd-asy}.
\end{proof}

\begin{Theorem} \label{thm:11-asy2}
Let $m$ be a positive integer or half-integer. Furthermore,
let $\et=(\et_1,\et_2,\dots,\et_n)$ be a vector of integers or of
half-integers in the
alcove $\mathcal A^{\tilde D_{n}}_m$ of type $\tilde D_n$ 
{\em(}defined in \eqref{eq:alcD}{\em)}.
Then, as $k$ tends to infinity,
the number of random walks which start at $\et$ and proceed for exactly $k$
diagonal steps,
which stay in the alcove $\mathcal A^{\tilde D_{n}}_m$, is asymptotically
\begin{multline} \label{eq:Dd-asy2a}
 \frac {2\cdot4^{n^2-n}} {(8m)^n} \bigg(2^n \prod _{j=0}
^{n-1}\cos\frac {j\pi} {2m}\bigg)^k\\
\times
 \prod _{1\le h<t\le n} ^{}\bigg(
\sin\frac {\pi (\et_h-\et_t)} {2m} \cdot
   \sin\frac {\pi (t-h)} {2m}
\sin\frac {\pi (\et_h+\et_t)} {2m} \cdot
   \sin\frac {\pi (t+h-2)} {2m}\bigg)\\
\times
\left(
     \frac{\prod_{h=1}^{n}\sin\frac {\pi (m-n+2h)} {2m}}
     {\prod_{h=1}^{n-1}\sin\frac {\pi h} {2m}}
 \prod _{1\le h<t\le n} ^{}
\frac {\sin\frac {\pi (m-n+h+t-1)} {2m}\cdot \sin\frac {\pi (m-n+t+h)} {2m}} 
{\sin^2\frac {\pi (t+h-2)} {2m}}
\right.
\\
\cdot
     \sum_{k=1}^{n}\frac {(-1)^{n-k}(m-n+2k)}
        {\prod_{h=1}^{n} 
\sin\frac {\pi (m-n+k+h)} {2m}
         \prod_{h=1}^{k-1}\sin\frac {\pi h} {2m}
         \prod_{h=1}^{n-k}\sin\frac {\pi h} {2m}}\\
+
     \frac{\prod_{h=1}^{n}\sin\frac {\pi (m-n+2h-1)} {2m}}
     {\prod_{h=1}^{n-1}\sin\frac {\pi h} {2m}}
 \prod _{1\le h<t\le n} ^{}
\frac {\sin\frac {\pi (m-n+h+t-2)} {2m}\cdot \sin\frac {\pi (m-n+t+h-1)} {2m}} 
{\sin^2\frac {\pi (t+h-2)} {2m}}
\\
\cdot
     \sum_{k=1}^{n}\frac {(-1)^{n-k}(m-n+2k-2)}
        {\prod_{h=1}^{n} 
\sin\frac {\pi (m-n+k+h-1)} {2m}
         \prod_{h=1}^{k-1}\sin\frac {\pi h} {2m}
         \prod_{h=1}^{n-k}\sin\frac {\pi h} {2m}}\\
+
     \frac{\prod_{h=1}^{n}\sin\frac {\pi (m-n+2h-1)} {2m}}
     {\prod_{h=1}^{n-1}\sin\frac {\pi h} {2m}}
 \prod _{1\le h<t\le n} ^{}
\frac {\sin\frac {\pi (m-n+h+t)} {2m}\cdot \sin\frac {\pi (m-n+t+h-1)} {2m}} 
{\sin^2\frac {\pi (t+h-2)} {2m}}
\\
\cdot
     \sum_{k=1}^{n}\frac {(-1)^{n-k}(m-n+2k-1)}
        {\prod_{h=1}^{n} 
\sin\frac {\pi (m-n+k+h-1)} {2m}
         \prod_{h=1}^{k-1}\sin\frac {\pi h} {2m}
         \prod_{h=1}^{n-k}\sin\frac {\pi h} {2m}}\\
+
     \frac{\prod_{h=1}^{n}\sin\frac {\pi (m-n+2h-2)} {2m}}
     {\prod_{h=1}^{n-1}\sin\frac {\pi h} {2m}}
 \prod _{1\le h<t\le n} ^{}
\frac {\sin\frac {\pi (m-n+h+t-1)} {2m}\cdot \sin\frac {\pi (m-n+t+h-2)} {2m}} 
{\sin^2\frac {\pi (t+h-2)} {2m}}
\\
\left.
\cdot
     \sum_{k=1}^{n}\frac {(-1)^{n-k}(m-n+2k-3)}
        {\prod_{h=1}^{n} 
\sin\frac {\pi (m-n+k+h-2)} {2m}
         \prod_{h=1}^{k-1}\sin\frac {\pi h} {2m}
         \prod_{h=1}^{n-k}\sin\frac {\pi h} {2m}}
\right)
\end{multline}
if $m$ is an integer and $k+2\et_j$ is even, it is asymptotically
\begin{multline} \label{eq:Dd-asy2b}
 \frac {4^{n^2-n+1}} {(8m)^n} \bigg(2^n \prod _{j=0}
^{n-1}\cos\frac {j\pi} {2m}\bigg)^k\\
\times
 \prod _{1\le h<t\le n} ^{}\bigg(
\sin\frac {\pi (\et_h-\et_t)} {2m} \cdot
   \sin\frac {\pi (t-h)} {2m}
\sin\frac {\pi (\et_h+\et_t)} {2m} \cdot
   \sin\frac {\pi (t+h-2)} {2m}\bigg)\\
\times
\left(
     \frac{\prod_{h=1}^{n}\sin\frac {\pi (\fl{m}-n+2h)} {2m}}
     {\prod_{h=1}^{n-1}\sin\frac {\pi h} {2m}}
 \prod _{1\le h<t\le n} ^{}
\frac {\sin\frac {\pi (\fl{m}-n+h+t-1)} {2m}\cdot \sin\frac {\pi
(\fl{m}-n+t+h)} {2m}} 
{\sin^2\frac {\pi (t+h-2)} {2m}}
\right.
\\
\cdot
     \sum_{k=1}^{n}\frac {(-1)^{n-k}(\fl{m}-n+2k)}
        {\prod_{h=1}^{n} 
\sin\frac {\pi (\fl{m}-n+k+h)} {2m}
         \prod_{h=1}^{k-1}\sin\frac {\pi h} {2m}
         \prod_{h=1}^{n-k}\sin\frac {\pi h} {2m}}\\
+
     \frac{\prod_{h=1}^{n}\sin\frac {\pi (\fl{m}-n+2h-1)} {2m}}
     {\prod_{h=1}^{n-1}\sin\frac {\pi h} {2m}}
 \prod _{1\le h<t\le n} ^{}
\frac {\sin\frac {\pi (\fl{m}-n+h+t)} {2m}\cdot \sin\frac {\pi
(\fl{m}-n+t+h-1)} {2m}} 
{\sin^2\frac {\pi (t+h-2)} {2m}}
\\
\left.\cdot
     \sum_{k=1}^{n}\frac {(-1)^{n-k}(\fl{m}-n+2k-1)}
        {\prod_{h=1}^{n} 
\sin\frac {\pi (\fl{m}-n+k+h-1)} {2m}
         \prod_{h=1}^{k-1}\sin\frac {\pi h} {2m}
         \prod_{h=1}^{n-k}\sin\frac {\pi h} {2m}}
\right)
\end{multline}
if $m$ is a half-integer and $k+2\et_j$ is even, it is asymptotically
\begin{multline} \label{eq:Dd-asy2c}
 \frac {2\cdot4^{n^2-n+1}} {(8m)^n} \bigg(2^n \prod _{j=0}
^{n-1}\cos\frac {j\pi} {2m}\bigg)^k\\
\times
 \prod _{1\le h<t\le n} ^{}\bigg(
\sin\frac {\pi (\et_h-\et_t)} {2m} \cdot
   \sin\frac {\pi (t-h)} {2m}
\sin\frac {\pi (\et_h+\et_t)} {2m} \cdot
   \sin\frac {\pi (t+h-2)} {2m}\bigg)\\
\times
     \frac{\prod_{h=1}^{n}\sin\frac {\pi (m-n+2h-1)} {2m}}
     {\prod_{h=1}^{n-1}\sin\frac {\pi h} {2m}}
 \prod _{1\le h<t\le n} ^{}
\frac {\sin^2\frac {\pi (m-n+h+t-1)} {2m}}
{\sin^2\frac {\pi (t+h-2)} {2m}}\\
\times
     \sum_{k=1}^{n}\frac {(-1)^{n-k}(m-n+2k-1)}
        {\prod_{h=1}^{n} 
\sin\frac {\pi (m-n+k+h-1)} {2m}
         \prod_{h=1}^{k-1}\sin\frac {\pi h} {2m}
         \prod_{h=1}^{n-k}\sin\frac {\pi h} {2m}}
\end{multline}
if $m$ is an integer and $k+2\et_j$ is odd, and it is asymptotically
\begin{multline} \label{eq:Dd-asy2d}
 \frac {4^{n^2-n+1}} {(8m)^n} \bigg(2^n \prod _{j=0}
^{n-1}\cos\frac {j\pi} {2m}\bigg)^k\\
\times
 \prod _{1\le h<t\le n} ^{}\bigg(
\sin\frac {\pi (\et_h-\et_t)} {2m} \cdot
   \sin\frac {\pi (t-h)} {2m}
\sin\frac {\pi (\et_h+\et_t)} {2m} \cdot
   \sin\frac {\pi (t+h-2)} {2m}\bigg)\\
\times
\left(
     \frac{\prod_{h=1}^{n}\sin\frac {\pi (m-n+2h-1/2)} {2m}}
     {\prod_{h=1}^{n-1}\sin\frac {\pi h} {2m}}
 \prod _{1\le h<t\le n} ^{}
\frac {\sin^2\frac {\pi (m-n+h+t-1/2)} {2m}}
{\sin^2\frac {\pi (t+h-2)} {2m}}
\right.\kern2cm
\\
\times
     \sum_{k=1}^{n}\frac {(-1)^{n-k}(m-n+2k-1/2)}
        {\prod_{h=1}^{n} 
\sin\frac {\pi (m-n+k+h-1/2)} {2m}
         \prod_{h=1}^{k-1}\sin\frac {\pi h} {2m}
         \prod_{h=1}^{n-k}\sin\frac {\pi h} {2m}}
\\
+
     \frac{\prod_{h=1}^{n}\sin\frac {\pi (m-n+2h-3/2)} {2m}}
     {\prod_{h=1}^{n-1}\sin\frac {\pi h} {2m}}
 \prod _{1\le h<t\le n} ^{}
\frac {\sin^2\frac {\pi (m-n+h+t-3/2)} {2m}}
{\sin^2\frac {\pi (t+h-2)} {2m}}\\
\left.
\times
     \sum_{k=1}^{n}\frac {(-1)^{n-k}(m-n+2k-3/2)}
        {\prod_{h=1}^{n} 
\sin\frac {\pi (m-n+k+h-3/2)} {2m}
         \prod_{h=1}^{k-1}\sin\frac {\pi h} {2m}
         \prod_{h=1}^{n-k}\sin\frac {\pi h} {2m}}
\right)
\end{multline}
if $m$ is a half-integer and $k+2\et_j$ is odd.
\end{Theorem}

\begin{proof}
Here, in view of Theorem~\ref{thm:11-asy}, we have to carry out the
sum of \eqref{eq:Dd-asy} over all possible $\la$, i.e., over all
$\la_1>\la_2>\dots>\la_{n-1}>\vert\la_n\vert$ with $\la_1+\la_2<2m$, where 
$k\equiv 2\et_j+2\la_j$ mod\/~$2$. 
Leaving away factors which are independent of
$\la$, the problem is to compute the sum
\begin{equation} \label{eq:SD2}
{\underset{\la_1+\la_2<2m}
{\sum _{\la_1>\la_2>\dots>\la_{n-1}>\vert\la_n\vert} 
^{}}\kern-1.2cm}{}^{\displaystyle\prime}
\kern1.3cm
 \prod _{1\le h<t\le n} ^{}\bigg(
   \sin\frac {\pi (\la_h-\la_t)} {2m}\cdot
   \sin\frac {\pi (\la_h+\la_t)} {2m}\bigg),
\end{equation}
where the sum $\sum{}^{\textstyle\prime}$ 
is restricted to integral, respectively to half-integral
$\la_1,\la_2,\dots,\la_n$, depending on whether $k+2\et_j$ is even or odd. 

As in previous proofs, we have to distinguish between several cases.
Let us first restrict the sum \eqref{eq:SD2} to integral
$\la_1,\la_2,\dots,\la_n$. Then, by repeating the argument from the
proof of Theorem~\ref{thm:10-asy2} of replacement of $\la_1$ by
$2m-\la_1$ and of replacement of $\la_n$ by $-\la_n$, 
and using notation from that proof, 
we obtain that the sum \eqref{eq:SD2} equals
$W_1(m)+W_1(m-1)+W_2(m)+W_2(m-1)$ if $m$ is an integer, and it equals
$2W_1(\fl m)+2W_2(\fl m)$ if $m$ is a half-integer. 
In the proof of Theorem~\ref{thm:10-asy2}
it was shown how to evaluate the sum $W_1(c)$ and $W_2(c)$ by means of
Lemmas~\ref{lem:spin} and \ref{lem:spin-1}. 
If the results are substituted, we arrive at
\eqref{eq:Dd-asy2a} and \eqref{eq:Dd-asy2b}.

If we restrict the sum \eqref{eq:SD2} to half-integral
$\la_1,\la_2,\dots,\la_n$, then we are faced with the problem of
evaluating the sum
\begin{equation*} 
X(c)=\sum _{c\ge\la_1>\dots>\la_n>0} ^{}
\kern-.8cm{}^{\displaystyle\prime}
\kern.8cm
 \prod _{1\le h<t\le n} ^{}\bigg( \sin\frac {\pi (\la_h-\la_t)} {2m}
\cdot \sin\frac {\pi (\la_h+\la_t)} {2m}\bigg)
\end{equation*}
where the sum $\sum{}^{\textstyle\prime}$ is
over all {\it half-integral\/} $\la_1,\la_2,\dots,\la_n$.
More specifically, using the replacement of $\la_1$ by $2m-\la_1$
and the replacement of $\la_n$ by $-\la_n$
another time, one sees readily that
the sum \eqref{eq:SD2} equals $4X(m-\frac {1} {2})$
if $m$ is an integer, and it equals $2X(m)+2X(m-1)$ if $m$ is a
half-integer.

In order to evaluate $X(c)$, where $c$ is a half-integer, 
we use \eqref{eq:cosdet} to rewrite it as
\begin{align*} 
\sum _{c\ge\la_1>\dots>\la_n>0} ^{}
\kern-.8cm{}^{\displaystyle\prime}
\kern.8cm
&
2^{-n^2+2n-1}
\det_{1\le h,t\le n}\big(\cos(\pi (h-1)\la_t/m)\big)\\
&=\sum _{c\ge\la_1>\dots>\la_n>0} ^{}
\kern-.8cm{}^{\displaystyle\prime}
\kern.8cm
\frac {1} {2^{n^2-n+1}}
\det_{1\le h,t\le n}\big((e^{\pi i (h-1)/m})^{\la_t}+(e^{\pi i
(h-1)/m})^{-\la_t}\big).
\end{align*}
Replacing $\la_t$ by $\la_t+n-t$, $t=1,2,\dots,n$, we obtain the
expression
\begin{equation} \label{eq:SummeD3}
\sum _{c-n+1\ge\la_1\ge\dots\ge\la_n\ge\frac {1} {2}} 
^{}\kern-1.2cm{}^{\displaystyle\prime}
\kern1.3cm
\frac {1} {2^{n^2-n+1}}
\det_{1\le h,t\le n}\big((e^{\pi i (h-1)/m})^{\la_t+n-t}+(e^{\pi i
(h-1)/m})^{-(\la_t+n-t)}\big).
\end{equation}
Again, this determinant can be expressed in terms of an
{\it even orthogonal character}. (See \eqref{e25} for the definition.)
Writing, as before, $q$ for $e^{\pi i/2m}$, the sum in
\eqref{eq:SummeD3}
equals
\begin{equation} \label{eq:sum-ortho2}
\frac {1} {2^{n^2-n+1}}
\det_{1\le h,t\le n}\big((q^{h-1})^{n-t}+(q^{h-1})^{-(n-t)}\big)\\
\sum _{c-n+1\ge\la_1\ge\dots\ge\la_n\ge\frac {1} {2}} 
^{}\kern-1.2cm{}^{\displaystyle\prime}
\kern1.3cm
\so^{even}_\la(q^{n-1},q^{n-2},\dots,1),
\end{equation}
where, still, the sum $\sum{}^{\textstyle\prime}$ is
over all half-integral $\la_1,\la_2,\dots,\la_n$.
(The reader should note that the second determinant in \eqref{e25} vanishes
if one of the variables $x_t$, $t=1,2,\dots,n$, is $1$, 
which is the case in our situation.)
The sum in \eqref{eq:sum-ortho2} can also be simplified by means of
\eqref{eq:soeven}.
Use of this formula in \eqref{eq:sum-ortho2} gives
\begin{multline*} 
\frac {1} {2^{n^2-n+1}}
\det_{1\le h,t\le n}\big((q^{h-1})^{n-t}+(q^{h-1})^{-(n-t)}\big)\\
\times
\so^{even}_{\big(((2c-2n+3)/4)^n\big)}(q^{n-1},q^{n-2},\dots,1)\\
\times
\sum _{p=-(2c-2n+1)/4} ^{(2c-2n+1)/4}
\so^{even}_{\big(((2c-2n+3)/4)^{n-1},p\big)}(q^{n-1},q^{n-2},\dots,1).
\end{multline*}
Again, the first determinant can be evaluated by means of
\eqref{eq:ortho3}, the even orthogonal character of shape $\big(((2c-2n+3)/4)^n\big)$
can be evaluated by using the definition \eqref{e25} and
\eqref{eq:ortho3}, while the sum over $p$ of even orthogonal characters is
evaluated in Lemma~\ref{lem:spin}.
Substitution of the result and simplification eventually
leads to the expressions \eqref{eq:Dd-asy2c} and \eqref{eq:Dd-asy2d}.
\end{proof}

\section*{Acknowledgements}
I am indebted to the referee for an extremely careful reading of the
manuscript, for the many suggestions which improved the presentation, 
and for pointing out a flaw in the original statement and 
proof of Lemma~\ref{lem:1}.

\bigskip

\appendix

\global\def\theTheorem{\mbox{A
}}
\setcounter{Theorem}{0}

\section{A saddle point approximation}\label{app:A}

The subject of this appendix is the saddle point approximation which
is needed in the proof of Theorem~\ref{thm:3-asy}.

\begin{Lemma} \label{lem:1}
Let $m$ be a positive integer and $d$ be an integer or half-integer.
Let $\om=e^{2\pi i/m}$, and let $\{r_1,r_2,\dots,r_n\}$ be a subset
of $\{0,1,\dots,m-1\}$, where the $r_j$'s are pairwise distinct. 
Furthermore, let 
\begin{equation} \label{eq:Cmax} 
C=C(r_1,\dots,r_n;m):=\max_\th
\prod _{j=1} ^{n}2\left\vert \cos\(\frac {\pi(\th+r_j)} {m}\)\right\vert,
\end{equation}
and let $\th_1,\th_2,\dots,\th_s$ be the values of $\th$ where the
maximum in \eqref{eq:Cmax} is attained.
Then, as $k$ tends to infinity such that 
$d+\frac {nk} {2}$ is an integer, 
the coefficient of $z^{d+\frac {nk} {2}}$ in 
\begin{equation} \label{eq:omprod}
\prod _{j=1} ^{n}(1+\om^{r_j}z)^k
\end{equation}
is asymptotically
\begin{equation} \label{eq:Sattel}
\frac {1} {\sqrt{2\pi k}}C^k
\sum _{\ell=1} ^{s}\ep(r_1,\dots,r_n;m)
\frac {\om^{\frac {k} {2}\sum _{j=1} ^{n}r_j- d\th_\ell}} { \sqrt{
c_0(\th_\ell) }},
\end{equation}
where 
$$\ep(r_1,\dots,r_n;m)=\sgn 
\prod _{j=1} ^{n}\cos\(\frac {\pi(\th+r_j)} {m}\)$$
and
$$c_0(\th)=\sum _{j=1} ^{n}\(2\cos\frac {\pi(\th+r_j)} {m}\)^{-2},
$$
except if the sum in \eqref{eq:Sattel} vanishes, in which case the
asymptotic order of the coefficient of $z^{d+\frac {nk} {2}}$ in
\eqref{eq:omprod} is strictly less than $C^k$.
\end{Lemma}

\begin{Remark}
The proof of the lemma below shows that, in fact, each $\th_\ell$ is a
solution of the equation
\begin{equation} \label{eq:th0}
\sum _{j=1} ^{n}{\tan\frac {\pi(\th+r_j)} {m}}=0.
\end{equation}
\end{Remark}

\begin{proof}[Proof of Lemma~\ref{lem:1}]
We apply the saddle point method (see \cite{FlSeAA}). (We would in
fact like to directly 
apply a general ``law of large powers," such as for example
Theorem~6.5 in \cite{FlSeAA}. However, I was not able to find an
applicable theorem in the literature. In particular, Theorem~6.5 from
\cite{FlSeAA} does not apply in our situation since it requires that
the coefficients of the series of which the powers are formed are
positive.)

We begin by writing the coefficient 
of $z^{d+\frac {nk} {2}}$ in the product \eqref{eq:omprod} as an
integral,
\begin{equation} \label{eq:intC}
\frac {1} {2\pi i}\int _{\mathcal C} ^{}
\frac {\prod _{j=1} ^{n}(1+\om^{r_j}z)^k}
{z^{d+\frac {nk} {2}}}\frac {dz} {z},
\end{equation}
where $\mathcal C$ is some contour encircling the origin
counter-clockwise. 
The task is to choose a contour $\mathcal C$ such that the contribution of the
integrand is concentrated very close to a finite number of points on the
contour (the saddle
points), whereas otherwise the contributions are negligible. The
equation for the saddle points is (see \cite[Sec.~6.3]{FlSeAA})
$$\frac {d} {dz}\(\frac {\prod _{j=1} ^{n}(1+\om^{r_j}z)^k}
{z^{d+\frac {nk} {2}}}\)=0.$$
Equivalently, this is
\begin{equation} \label{eq:SP}
\frac {k} {2}\sum _{j=1} ^{n}\frac {\om^{r_j}z-1} {\om^{r_j}z+1}-d=0.
\end{equation}
For the saddle point method to work, it is often not necessary to
determine the solution(s) of the saddle point equation exactly. It
is often sufficient to find a suitable approximation for large $k$.
This is also the case here. Clearly, if $k$ is large, then the
coefficient of $k/2$ in \eqref{eq:SP} must be very small. Therefore,
as a first approximation, we consider solutions of the equation
\begin{equation} \label{eq:z0}
\sum _{j=1} ^{n}\frac {\om^{r_j}z-1} {\om^{r_j}z+1}=0.
\end{equation}
This is (after denominators have been cleared) a polynomial equation
of degree $n$. Therefore we must expect (up to) $n$ different
solutions to \eqref{eq:z0}. 

We claim that all $n$ solutions to \eqref{eq:z0} are distinct and
moreover have modulus $1$. To see this, we substitute $e^{2\pi i\th/m}$ 
for $z$ in
\eqref{eq:z0}. After a little manipulation, we obtain the equivalent
equation
\begin{equation} \label{eq:th}
\sum _{j=1} ^{n}{\tan\frac {\pi(\th+r_j)} {m}}=0.
\end{equation}
The left-hand side of \eqref{eq:th} is defined for $0\le\th\le m$
except at the values (which have to be taken modulo $m$)
$$\th=-r_1+\frac m2,-r_2+\frac m2,\dots,-r_n+\frac m2,$$
all of which are distinct by assumption. These values
are simple poles of the function on the left-hand side of
\eqref{eq:th}. Between two successive poles the function
is monotone increasing and continuous
(here, we regard the interval $[0,m]$ as a circular interval,
identifying $0$ and $m$), ranging
there from $-\infty$ to $+\infty$. Hence, in each open interval
bounded by two successive poles there lies exactly 
one solution to \eqref{eq:th}. Since
the number of poles, and, thus, of such intervals, is $n$,
there are $n$ values of $\th$ in the range $[0,m)$ satisfying
equation~\eqref{eq:th}. In turn, each solution to \eqref{eq:th}
produces a solution to \eqref{eq:z0} with modulus $1$. This
establishes that, indeed, all $n$ solutions to \eqref{eq:z0} have
modulus $1$ and are pairwise distinct.

Let $\th_1,\th_2,\dots,\th_n$ be the above solutions to \eqref{eq:th},
and, for a given $\th_\ell$, 
let $z_\ell=e^{2\pi i\th_\ell/m}$ be the corresponding
solution to \eqref{eq:z0}. Obviously, we have $1+\om^{r_j}z_\ell\ne0$
for all $\ell$ and $j$. 

We shall need a slightly better approximation of the saddle points,
which will approximate
\eqref{eq:SP} up to second order. Let $y_\ell:=z_\ell(1+\frac {c_\ell} {k})$.
Substitution of $y_\ell$ in \eqref{eq:SP} leads to the equation
\begin{align*} \label{}
\frac {k} {2}\sum _{j=1} ^{n}&\frac {\om^{r_j}y_\ell-1}
{\om^{r_j}y_\ell+1}-d\\
&=\frac {k} {2}\sum _{j=1} ^{n}\frac {\om^{r_j}z_\ell-1} {\om^{r_j}z_\ell+1}
\(1+\frac {c_\ell\om^{r_j}z_\ell} {k(\om^{r_j}z_\ell-1)}-
\frac {c_\ell\om^{r_j}z_\ell} {k(\om^{r_j}z_\ell+1)}+O\(k^{-2}\)\)-d\\
&=
c_\ell\sum _{j=1} ^{n}
\frac {\om^{r_j}z_\ell} {(\om^{r_j}z_\ell+1)^2}-d+O\(k^{-1}\)\\
&=
c_\ell\sum _{j=1} ^{n}
\frac {1} {\(2\cos\frac {\pi (r_j+\th_\ell)} {m}\)^2}-d+O\(k^{-1}\)\\
&=0,
\end{align*}
where we used \eqref{eq:z0} in the third line. Thus, $y_\ell$ will
approximate a solution to \eqref{eq:SP} up to second order, if we
choose
$$c_\ell=\frac {d} {\displaystyle\sum _{j=1} ^{n}
\(2\cos\frac {\pi (r_j+\th_\ell)} {m}\)^{-2}}.$$

Returning to the integral \eqref{eq:intC} that we want to
approximate, we choose as the contour $\mathcal C$ the circle
\begin{equation} \label{eq:contour} 
\textstyle\left\{\(1+\frac {c_\ell} {k}\)e^{it}:
0\le t\le 2\pi\right\},
\end{equation}
which contains the (approximated) 
saddle points $y_\ell$. Thus, writing $\bar\th$ for
$2\pi \th_\ell/m$ for short, the portion 
of the integral \eqref{eq:intC} contributed by a neighbourhood
$[\bar\th-\ep,\bar\th+\ep]$ of
$\bar\th$ becomes
\begin{equation} \label{eq:int}
\frac {1} {2\pi }\int _{-\ep} ^{\ep}
\frac {\prod _{j=1} ^{n}\(1+\om^{r_j}(1+\frac {c_\ell }
{k})e^{i(\bar\th+t)}\)^k}
{\((1+\frac {c_\ell } {k})e^{i(\bar\th+t)}\)^{d+\frac {nk} {2}}} {dt}.
\end{equation}
In order to estimate this integral, we compute a Taylor expansion of 
$$\log\(1+\om^{r_j}\(1+\frac {c_\ell }
{k}\)e^{i(\bar\th+t)}\),$$ 
thereby also neglecting terms which are of order
$O(k^{-2})$:
\begin{multline*} \label{}
\log\(1+\om^{r_j}\(1+\frac {c_\ell }
{k}\)e^{i(\bar\th+t)}\)\\=
\log (1 + {e^{{ {i  {\bar\th}} }}} \om^{r_j}) 
 + 
\frac {c_\ell } {k}  {\frac { {e^{{ {i  {\bar\th}} }}} \om^{r_j} 
 } 
    { {{\left( 1 + {e^{{ {i  {\bar\th}} }}} \om^{r_j} \right) }
        } }}
 + it\( {\frac { {e^{{ {i  {\bar\th}} }}}  \om^{r_j}} 
      {1 + {e^{{ {i  {\bar\th}} }}} \om^{r_j}}}
+ \frac {c_\ell } {k}{\frac {  {e^{{ {i  {\bar\th}} }}} \om^{r_j} 
 } 
    { {{\left( 1 + {e^{{ {i  {\bar\th}} }}} \om^{r_j} \right) }^
        2} }}\)\\
- 
  \frac {t^2} {2}\(  {\frac {{e^{{ {i  {\bar\th}} }}}  \om^{r_j}} 
      { {{\left( 1 + {e^{{ {i  {\bar\th}} }}} \om^{r_j} \right) }^
          2}}}
+  \frac {c_\ell } {k}{\frac { {e^{{ {i  {\bar\th}} }}} \om^{r_j} 
      \left(    1 - 
        {e^{{ {i  {\bar\th}} }}} {} \om^{r_j} \right) } 
    { {{\left( 1 + {e^{{ {i  {\bar\th}} }}} \om^{r_j} \right) }^
        3} }}\)
+ O(t^3)+O\(k^{-2}\).
\end{multline*}
This expansion is only valid for 
\begin{equation} \label{eq:tkBed} 
\left\vert\frac {\(1+\frac {c_\ell }
{k}\)e^{i(\bar\th+t)}-e^{i\bar\th}\om^{r_j}} {
1+e^{i\bar\th}\om^{r_j}}\right\vert< 1,
\end{equation}
whence, in general, it will only
be valid for $t$ in a neighbourhood of $0$ and $k$ large enough.
More precisely, choose $k_0$ so that the left-hand side of
\eqref{eq:tkBed} with $t=0$ is less than $1$, i.e.,
$$
\left\vert\frac { {c_\ell }
} {k_0(
1+e^{i\bar\th}\om^{r_j})}\right\vert< 1.
$$
Then the above expansion of the logarithm will be valid for $-\ep\le
t\le \ep$ and $k\ge k_0$ for some $\ep$ which is independent of $k$.

If this is substituted in \eqref{eq:int}, then, after some
simplification based on \eqref{eq:th} and the definition of $c_\ell $, we
obtain
\begin{multline} \label{eq:intexp}
\frac {1} {2\pi e^{ i\bar\th(d+\frac {nk} {2}) }}
\prod _{j=1} ^{n} (1 + {e^{{ {i  {\bar\th}} }}} \om^{r_j})^k\\
\times
\int _{-\ep} ^{\ep}
\exp\bigg( {
-   \frac {kt^2} {2}\Big(  \textstyle
\sum\limits _{j=1} ^{n}{\frac {{e^{{ {i  {\bar\th}} }}}  \om^{r_j}} 
      { {{\left( 1 + {e^{{ {i  {\bar\th}} }}} \om^{r_j} \right) }^
          2}}}
+ \frac {c_\ell } k\sum\limits _{j=1} ^{n} {\frac { {e^{{ {i  {\bar\th}} }}} \om^{r_j} 
      \left(    1 - 
        {e^{{ {i  {\bar\th}} }}} {} \om^{r_j} \right) } 
    { {{\left( 1 + {e^{{ {i  {\bar\th}} }}} \om^{r_j} \right) }^
        3} }}\Big)
+ O(kt^3)+O\(k^{-1}\)
}\bigg).
\end{multline}

Next we estimate the integral in \eqref{eq:intexp}. One argues
as is usual in such a situation: the essential contribution to the
integral comes from the range $\vert t\vert<(\log k)/\sqrt k$, the rest
being of order $O\(k^{-1}\)$. In this range, the error terms
$ O(kt^3)+O\(k^{-1}\)$ are of the order $O\big((\log^3k)/\sqrt k\big)$.
Finally, one extends the integral to the complete range
$-\infty<t<\infty$, upon making another error of the order
$O\(k^{-1}\)$. Leaving these details to the reader, up to an error of
$O\big((\log^3k)/\sqrt k\big)$ the integral in \eqref{eq:intexp} is
asymptotically
\begin{equation} \label{eq:intasy}
\int _{-\infty} ^{\infty}\exp\(-S\frac {kt^2} {2}\)=\sqrt{\frac
{2\pi} {kS}},
\end{equation}
where 
$$S=\sum _{j=1} ^{n}{\frac {{e^{{ {i  {\bar\th}} }}}  \om^{r_j}} 
      { {{\left( 1 + {e^{{ {i  {\bar\th}} }}} \om^{r_j} \right) }^
          2}}
+ \frac {c_\ell } k\sum _{j=1} ^{n} {\frac { {e^{{ {i  {\bar\th}} }}} \om^{r_j} 
      \left(    1 - 
        {e^{{ {i  {\bar\th}} }}} {} \om^{r_j} \right) } 
    { {{\left( 1 + {e^{{ {i  {\bar\th}} }}} \om^{r_j} \right) }^
        3} }}}
=
c_0(\th_\ell)+O\(\frac {1} {k}\).
$$
with $c_0(\th_\ell)$ as given in the statement of the theorem.
The evaluation \eqref{eq:intasy} is valid since $c_0>0$, and hence
$\Re S>0$ for $k$ large enough.

Clearly, the second term in the definition of $S$ makes only a
negligible contribution to the asymptotic behaviour of
\eqref{eq:intasy}. If the result is
substituted in \eqref{eq:intexp}, we obtain that the portion of the
integral \eqref{eq:intC} near the saddle point $z_\ell$, given
explicitly in \eqref{eq:int}, contributes
$$
\frac {\prod _{j=1} ^{n}(1+\om^{r_j+\th_\ell})^k}
{\om^{\th_\ell(d+\frac {nk} {2})}}
\frac {1} { \sqrt{2\pi
c_0(\th_\ell) k}},
$$
as $k$ tends to infinity,
up to an error of $O\big((\log^3k)/\sqrt k\big)$.
Summation of these contributions, neglecting those which are
asymptotically smaller, yields indeed the expression
\eqref{eq:Sattel} after some manipulation.

It remains to show that the contribution of the remaining portions of
the integral \eqref{eq:intC} is negligible. This is a routine matter,
the details of which we leave to the reader. We content ourselves to
indicate that one may start with the observation that the modulus of
the integrand of \eqref{eq:intC} along the contour \eqref{eq:contour}
is bounded by
\begin{equation} \label{eq:bound} 
\left\vert \frac {\prod _{j=1} ^{n}\(1+\om^{r_j}(1+\frac {c_\ell }
{k})e^{it}\)^k}
{\((1+\frac {c_\ell } {k})e^{it}\)^{d+\frac {nk}
{2}}}\right\vert<
\text{constant}^{\sqrt k}
\left\vert  {\prod _{j=1} ^{n}\(1+\om^{r_j}
e^{it}\)^k}
\right\vert,
\end{equation} 
the constant depending on $r_1,r_2,\dots,r_n$ but not on $k$, as long
as $t$ stays away from tiny neighbourhoods of $\pi-\frac {2\pi r_j}
{m}$, $j=1,2,\dots,n$, determined by the condition that 
$\vert 1+\om^{r_j}e^{it}\vert\le 1/\sqrt k$,
and that the right-hand side of \eqref{eq:bound} 
is exponentially smaller than $C^k$ for
any $t$ in these remaining portions, while the contributions in these
tiny neighbourhoods is even super-exponentially small.. 
\end{proof}

\global\def\theTheorem{\mbox{B
}}
\setcounter{Theorem}{0}

\section{Some determinants}

In our computations we need frequently the following determinant evaluations.
All of them are readily proved by
the standard argument that proves Vandermonde-type determinant
evaluations.

\begin{Lemma} \label{lem:dets}
Let $n$ by a non-negative integer. Then
\begin{align} \label{eq:ortho1}
\det\limits_{1\le h,t\le n}\(x_h^{t-1/2}+x_h^{1/2-t}\)&=
(x_1x_2\cdots x_n)^{-n+1/2}
\prod _{1\le h<t\le n} ^{}(x_h-x_t)(1-x_hx_t)\prod _{h=1}
^{n}(x_h+1),\\
\label{eq:ortho2}
\det\limits_{1\le h,t\le n}\(x_h^{t-1/2}-x_h^{-t+1/2}\)&=
(x_1x_2\cdots x_n)^{-n+1/2}
\prod _{1\le h<t\le n} ^{}(x_h-x_t)(1-x_hx_t)\prod _{h=1}
^{n}(x_h-1),\\
\label{eq:ortho3}
\det_{1\le h,t\le n}(x_h^{t-1}+x_h^{-(t-1)})&=
2\cdot(x_1\cdots x_n)^{-n+1}\prod
_{1\le h<t\le n} ^{}(x_h-x_t)(1-x_hx_t),\\
\label{eq:sympl}
\det_{1\le h,t\le n}(x_h^t-x_h^{-t})&=(x_1\cdots x_n)^{-n}\prod
_{1\le h<t\le n} ^{}(x_h-x_t)(1-x_hx_t)
\prod _{h=1} ^{n}(x_h^2-1).
\end{align}
\end{Lemma}

\section{Odd and even orthogonal characters, and Schur functions 
at special values of the arguments}

\global\def\theTheorem{\mbox{C\arabic{Theorem}}}
\setcounter{Theorem}{0}

In this appendix we provide the evaluations of odd
orthogonal characters, Schur functions 
of rectangular shape and nearly rectangular shape, and of certain
sums of even orthogonal characters, where the variables are
specialized in peculiar ways. The evaluations of the odd
orthogonal characters are needed for
the proof of Theorem~\ref{thm:5-asy2} (on which, in turn, hinges also the
proof of Theorem~\ref{thm:1-asy2}) and the proof of 
Theorem~\ref{thm:4-asy2}, 
the evaluations of the special Schur functions are needed for the
proofs of Theorems~\ref{thm:6-asy2}, \ref{thm:8-asy2}, and
\ref{thm:9-asy2}, while the evaluations of the sums of specialized
even orthogonal characters are needed for the proofs of
Theorems~\ref{thm:10-asy2} and \ref{thm:11-asy2}.

Recall that, given a partition
$\la=(\la_1,\la_2,\dots,\la_n)$ 
(i.e., a non-increasing sequence of non-negative integers) or 
half-partition (which, by definition, is 
a non-increasing sequence of positive odd integers
divided by 2), the {\it odd orthogonal character}\break $\so^{odd}_\la(x_1,x_2,\dots,x_n)$
is defined by (see \cite[(24.28)]{FuHaAA})
\begin{equation} \label{e12}
\so^{odd}_\la(x_1,x_2,\dots,x_n)=\frac {\det\limits_{1\le h,t\le
n}(x_h^{\la_t+n-t+1/2}-x_h^{-(\la_t+n-t+1/2)})}
{\det\limits_{1\le h,t\le
n}(x_h^{n-t+1/2}-x_h^{-(n-t+1/2)})}.
\end{equation}
It is not difficult to see that the denominator in \eqref{e12} does in
fact cancel out, so that any odd orthogonal character 
$\so^{odd}_\la(x_1,x_2,\dots,x_n)$ is in fact a {\it Laurent polynomial\/} in
$x_1,x_2,\dots,x_n$ (i.e., a polynomial in 
$x_1,x_1^{-1},x_2,x_2^{-1},\dots,x_n,x_n^{-1}$),
and is thus well-defined for any choice of the
variables $x_1,x_2,\dots,x_n$ such that all of them are non-zero. 

As earlier in Section~\ref{sec:4},
we shall use the notation $\big((p/2)^n\big)$ for the
vector of $n$ components, all of them equal to $p/2$.

\begin{Lemma} \label{lem:so1}
Let $m$ and $n$ be positive integers with $m\ge n$, and
let $q=e^{\pi i/m}$. Then we have
\begin{equation} \label{eq:so1a}
\so^{odd}_{\big(((m-n)/2)^n\big)}(q^{n-1},q^{n-3},\dots,q^{-n+3},q^{-n+1})=
\frac {\displaystyle m^{n/2}
\prod _{h=1} ^{n/2}\cot\frac {(2h-1)\pi} {2m}}
{\displaystyle 2^{\binom n2}\prod _{1\le h<t\le n}
^{}\sin\frac {(t-h)\pi} {m}
}
\end{equation}
if $n$ is even, and 
\begin{equation} \label{eq:so1b}
\so^{odd}_{\big(((m-n)/2)^n\big)}(q^{n-1},q^{n-3},\dots,q^{-n+3},q^{-n+1})=
\frac {\displaystyle m^{(n+1)/2}
\prod _{h=1} ^{(n-1)/2}\cot\frac {h\pi} {m}}
{\displaystyle 2^{\binom n2}\prod _{1\le h<t\le n}
^{}\sin\frac {(t-h)\pi} {m}}
\end{equation}
if $n$ is odd.
\end{Lemma}

\begin{proof}
In principle, we would like to specialize in the
definition \eqref{e12} of the odd orthogonal character. 
However, we face the difficulty that, because
of the peculiar choice of the variables $x_1,x_2,\dots,x_n$ in
\eqref{eq:so1a} and \eqref{eq:so1b}, 
direct substitution gives an indeterminate
expression $0/0$. More precisely, the problem is the pairs of
reciprocal variables (such as $q^{n-1}$ and $q^{-n+1}$, $q^{n-3}$ and
$q^{-n+3}$, etc.). In the case
that $n$ is odd, we have the additional problem that one of the
variables is $q^0=1$. 

To overcome this problem, we have recourse to de l'H\^opital's rule.
We have to distinguish between two cases,
depending on whether $n$ is even or odd.

Let first $n$ be even. In that case we must compute the limit as
\begin{equation} \label{eq:lim1}
x_1\to q^{n-1},\ 
x_2\to q^{n-3},\ \dots,\
x_{n/2}\to q
\end{equation}
of the right-hand side of \eqref{e12} with $x_{n/2+1}=q^{-1}$, \dots,
$x_{n-1}=q^{-n+3}$, $x_n=q^{-n+1}$. Upon a simultaneous 
rearrangement and change of signs 
of the rows of the two determinants in \eqref{e12},
one sees that this is equivalent to computing the limit
\eqref{eq:lim1} of $(\det A_1)/(\det B_1)$, where
\begin{equation*} \label{eq:detA1}
A_1=\begin{pmatrix} 
x_h^{\frac {m+n+1} {2}-t}-x_h^{-(\frac {m+n+1} {2}-t)}&
1\le h\le n/2\\
q^{(2n-2h+1)(\frac {m+n+1} {2}-t)}-q^{-(2n-2h+1)(\frac {m+n+1} {2}-t)}&
n/2<h\le n
\end{pmatrix}_{1\le h,t\le n}
\end{equation*}
and
\begin{equation*} \label{eq:detB1a}
B_1=\begin{pmatrix} 
x_h^{n-t+\frac {1} {2}}-x_h^{-(n-t+\frac {1} {2})}&
1\le h\le n/2\\
q^{(2n-2h+1)(n-t+\frac {1} {2})}-q^{-(2n-2h+1)(n-t+\frac {1} {2})}&
n/2<h\le n
\end{pmatrix}_{1\le h,t\le n}.
\end{equation*}
Now applying de l'H\^opital's rule, this limit is equal to the limit
\eqref{eq:lim1} of 
\begin{equation} \label{eq:det/det}
\frac {\det A_2} {
\frac {\partial} {\partial x_1}
\frac {\partial} {\partial x_2}\cdots
\frac {\partial} {\partial x_{n/2}}
\det B_1},
\end{equation}
where
\begin{align} 
\notag
A_2&=\begin{pmatrix} 
({\frac {m+n+1} {2}-t})
(q^{(n-2h+1)(\frac {m+n+1} {2}-t-1)}+q^{-(n-2h+1)(\frac {m+n+1}
{2}-t+1)})&
1\le h\le n/2\\
q^{(2n-2h+1)(\frac {m+n+1} {2}-t)}-q^{-(2n-2h+1)(\frac {m+n+1} {2}-t)}&
n/2<h\le n
\end{pmatrix}\\
\label{eq:detA2}
&=\begin{pmatrix} 
 (-1)^{\frac {n} {2}-h}i\,q^{-(n-2h+1)}({\frac {m+n+1} {2}-t})\hfill\\
\times
( q^{(n-2h+1)(\frac {n+1} {2}-t)}- q^{-(n-2h+1)(\frac {n+1}
{2}-t)})&
1\le h\le n/2\\
 (-1)^{{n}-h}i\,
(q^{(2n-2h+1)(\frac {n+1} {2}-t)}+q^{-(2n-2h+1)(\frac {n+1} {2}-t)})&
n/2<h\le n
\end{pmatrix}.
\end{align}
The limit of the denominator of \eqref{eq:det/det} is readily
obtained, because the determinant of $B_1$ can actually be evaluated
by means of \eqref{eq:ortho2}. To be precise, with
\begin{equation} \label{eq:xn}
x_{n/2+1}=q^{n-1},\ \dots,\ x_{n-1}=q^{3},\ x_n=q,
\end{equation}
we have
\begin{equation} \label{eq:detB1}
\det B_1=(x_1x_2\cdots x_n)^{-n+1/2}
\prod _{1\le h<t\le n} ^{}(x_t-x_h)(1-x_hx_t)\prod _{h=1}
^{n}(x_h-1).
\end{equation}
Thus, if we differentiate the product on the right-hand side of
\eqref{eq:detB1} with
respect to $x_1$, say, (using the product rule, of course),
and subsequently set $x_1=q^{n-1}$
(which is exactly what we want to do; see \eqref{eq:lim1}), 
then it is just one term in
the derivative which contributes, namely
\begin{equation*} 
\frac {1} {x_1-x_{n/2+1}}(x_1x_2\cdots x_n)^{-n+1/2}
\prod _{1\le h<t\le n} ^{}(x_t-x_h)(1-x_hx_t)\prod _{h=1}
^{n}(x_h-1),
\end{equation*}
in which the factor $(x_1-x_{n/2+1})$ cancels;
all other terms vanish because of the occurrence of the factor
$(x_1-x_{n/2+1})$. (Recall that $x_{n/2+1}=q^{n-1}$; see
\eqref{eq:xn}.) An analogous argument applies to the other pairs
of variables. Hence, the denominator of \eqref{eq:det/det} is equal to
\begin{multline} \label{eq:Nenn1}
(-1)^{\binom {n/2}2}(q^{n^2/2})^{-n+1/2}
\prod _{1\le h<t\le n/2}
^{}(q^{n-2t+1}-q^{n-2h+1})^4(1-q^{2n-2h-2t+2})^4\\
\times
\prod _{h=1} ^{n/2}(q^{n-2h+1}-1)^2(1-q^{2n-4h+2}).
\end{multline}

Next we devote ourselves to the evaluation of the determinant of
$A_2$, with $A_2$ given by \eqref{eq:detA2}. For $t=1,2,\dots,\frac
{n} {2}$ we subtract column $n-t+1$ from column $t$. As a result we
obtain that $\det A_2$ is equal to $\det A_3$, where $A_3$ is the
$n\times n$ block matrix
\begin{equation} \label{eq:detA3}
A_3=\begin{pmatrix} A_3^{(1)}&*\\
0&A_3^{(2)}\end{pmatrix},
\end{equation}
with $A_3^{(1)}$ the $\frac {n} {2}\times \frac {n} {2}$ matrix
$$A_3^{(1)}=\Big(
 (-1)^{\frac {n} {2}-h}i\,q^{-(n-2h+1)}m
( q^{(n-2h+1)(\frac {n+1} {2}-t)}- q^{-(n-2h+1)(\frac {n+1}
{2}-t)})\Big)_{1\le h,t\le n/2},$$
and $A_3^{(2)}$ the $\frac {n} {2}\times \frac {n} {2}$ matrix
$$A_3^{(2)}=\Big(
 (-1)^{\frac {n} {2}-h}i\,
(q^{(n-2h+1)(\frac {1} {2}-t)}+q^{-(n-2h+1)(\frac {1} {2}-t)})
\Big)_{1\le h,t\le n/2}.$$
Clearly, the determinant $\det A_3$ is equal to the product $(\det
A_3^{(1)})\cdot(\det A_3^{(2)})$. Using again \eqref{eq:ortho2}, we have
\begin{multline} \label{eq:detA31}
\det A_3^{(1)}=i^{n/2}q^{-n^2/4}m^{n/2}
(q^{n^2/4})^{-n/2+1/2}\\
\times
\prod _{1\le h<t\le n/2} ^{}(q^{n-2h+1}-q^{n-2t+1})
(1-q^{2n-2h-2t+2})\prod _{h=1}
^{n/2}(q^{n-2h+1}-1),
\end{multline}
whereas by means of \eqref{eq:ortho1} we have
\begin{multline} \label{eq:detA32}
\det A_3^{(2)}=(-1)^{\binom {n/2}2}i^{n/2}
(q^{n^2/4})^{(1-n)/2}\\
\times
\prod _{1\le h<t\le n/2} ^{}(q^{n-2h+1}-q^{n-2t+1})(1-q^{2n-2h-2t+2})
\prod _{h=1} ^{n/2}(1+q^{n-2h+1}).
\end{multline}
If we now combine \eqref{eq:detA3}, \eqref{eq:detA31},
\eqref{eq:detA32}, and \eqref{eq:Nenn1}, use the fact that $\det
A_2=\det A_3$, and substitute all this in
\eqref{eq:det/det}, then, after some simplification, we obtain
\begin{align*} \label{eq:so1}
\so&^{odd}_{\big(((m-n)/2)^n\big)}(q^{n-1},q^{n-3},\dots,q^{-n+3},q^{-n+1})\\
&=\frac {\displaystyle (-1)^{n/2}m^{n/2}
\prod _{h=1} ^{n/2}(1+q^{n-2h+1})}
{\displaystyle \prod _{1\le h<t\le n/2}
^{}(q^{h-t}-q^{t-h})^2
\prod _{1\le h,t\le n/2}
(q^{-n+h+t-1}-q^{n-h-t+1})
\prod _{h=1} ^{n/2}(1-q^{n-2h+1})
},
\end{align*}
which can be rewritten as \eqref{eq:so1a}.

\medskip
Now let $n$ be odd. We proceed in a completely analogous manner. Again, the
task is to compute the specialized odd orthogonal character, 
by means of a limit of its definition \eqref{e12}. Here, the
difficulty is not only the pairs of reciprocal variables, but also
that one of the variables is $q^0=1$.

What we must compute is the limit as
\begin{equation} \label{eq:lim2}
x_1\to q^{n-1},\ 
x_2\to q^{n-3},\ \dots,\
x_{(n-1)/2}\to q^2,\
x_{(n+1)/2}\to 1,
\end{equation}
of the right-hand side of \eqref{e12} with $x_{(n+3)/2}=q^{-2}$, \dots,
$x_{n-1}=q^{-n+3}$, $x_n=q^{-n+1}$. Upon a simultaneous 
rearrangement and change of signs 
of the rows of the two determinants in \eqref{e12},
one sees that this is equivalent to computing the limit
\eqref{eq:lim2} of $(\det A_4)/(\det B_4)$, where
\begin{equation*} \label{eq:detA4}
A_4=\begin{pmatrix} 
x_h^{\frac {m+n+1} {2}-t}-x_h^{-(\frac {m+n+1} {2}-t)}&
1\le h\le (n+1)/2\\
q^{(2n-2h+2)(\frac {m+n+1} {2}-t)}-q^{-(2n-2h+2)(\frac {m+n+1} {2}-t)}&
(n+1)/2<h\le n
\end{pmatrix}_{1\le h,t\le n}
\end{equation*}
and
\begin{equation*} \label{eq:detB4a}
B_4=\begin{pmatrix} 
x_h^{n-t+\frac {1} {2}}-x_h^{-(n-t+\frac {1} {2})}&
1\le h\le (n+1)/2\\
q^{(2n-2h+2)(n-t+\frac {1} {2})}-q^{-(2n-2h+2)(n-t+\frac {1} {2})}&
(n+1)/2<h\le n
\end{pmatrix}_{1\le h,t\le n}.
\end{equation*}
Applying de l'H\^opital's rule, this limit is equal to the limit
\eqref{eq:lim2} of 
\begin{equation} \label{eq:det/det2}
\frac {\det A_5} {
\frac {\partial} {\partial x_1}
\frac {\partial} {\partial x_2}\cdots
\frac {\partial} {\partial x_{(n+1)/2}}
\det B_4},
\end{equation}
where
\begin{align} 
\notag
A_5&=\begin{pmatrix} 
({\frac {m+n+1} {2}-t})
(q^{(n-2h+1)(\frac {m+n+1} {2}-t-1)}+q^{-(n-2h+1)(\frac {m+n+1}
{2}-t+1)})&
1\le h\le (n+1)/2\\
q^{(2n-2h+2)(\frac {m+n+1} {2}-t)}-q^{-(2n-2h+2)(\frac {m+n+1} {2}-t)}&
(n+1)/2<h\le n
\end{pmatrix}\\
\label{eq:detA5}
&=\begin{pmatrix} 
 (-1)^{\frac {n+1} {2}-h}q^{-(n-2h+1)}({\frac {m+n+1} {2}-t})\hfill\\
\times
( q^{(n-2h+1)(\frac {n+1} {2}-t)}+ q^{-(n-2h+1)(\frac {n+1}
{2}-t)})&
1\le h\le (n+1)/2\\
 (-1)^{{n}+1-h}
(q^{(2n-2h+2)(\frac {n+1} {2}-t)}-q^{-(2n-2h+2)(\frac {n+1} {2}-t)})&
(n+1)/2<h\le n
\end{pmatrix}.
\end{align}
The limit of the denominator of \eqref{eq:det/det2} is 
obtained in the same way as before, by using the complete
factorization of $\det B_4$ by means of \eqref{eq:ortho2}. 
(It should be noted that the 
difference between the denominator of \eqref{eq:det/det2} and that of
\eqref{eq:det/det} is the number of differentiations, and the number
of variables with respect to which the limit is taken.)
After a small calculation,
it turns out that the denominator of \eqref{eq:det/det2} is equal to
\begin{multline} \label{eq:Nenn2}
(-1)^{\binom {(n+1)/2}2}
(q^{(n^2-1)/2})^{-n+1/2}
\prod _{1\le h<t\le (n-1)/2} ^{}(q^{n-2t+1}-q^{n-2h+1})^4
(1-q^{2n-2h-2t+2})^4\\
\times
\prod _{h=1} ^{(n-1)/2}(q^{n-2h+1}-1)^6(1-q^{2n-4h+2}).
\end{multline}

Next we devote ourselves to the evaluation of the determinant of
$A_5$, with $A_5$ given by \eqref{eq:detA5}. For $t=1,2,\dots,\frac
{n-1} {2}$ we add column $n-t+1$ to column $t$. As a result we
obtain that $\det A_5$ is equal to $\det A_6$, where $A_6$ is the
$n\times n$ block matrix
\begin{equation} \label{eq:detA6}
A_6=\begin{pmatrix} A_6^{(1)}&*\\
0&A_6^{(2)}\end{pmatrix},
\end{equation}
with $A_6^{(1)}$ the $\frac {n+1} {2}\times \frac {n+1} {2}$ matrix
\begin{multline*} 
A_6^{(1)}
=\Bigg(
 (-1)^{\frac {n+1} {2}-h}q^{-(n-2h+1)}m\\
\times
\begin{cases} 
q^{(n-2h+1)(\frac {n+1} {2}-t)}+ q^{-(n-2h+1)(\frac {n+1}
{2}-t)}&1\le t\le \frac {n-1} {2}\\
1&t=\frac {n+1} {2}
\end{cases}
\Bigg)_{1\le h,t\le (n+1)/ {2}},
\end{multline*}
and $A_6^{(2)}$ the $\frac {n-1} {2}\times \frac {n-1} {2}$ matrix
$$A_6^{(2)}=\Big(
 (-1)^{\frac {n+1} {2}-h}
(q^{(n-2h+1)(-t)}-q^{-(n-2h+1)(-t)})
\Big)_{1\le h,t\le (n-1)/2}.$$
(It should be noted that the entries in column $(n+1)/2$ of
$A_6^{(1)}$ are exactly a half of what would result from 
substituting $t=(n+1)/2$
in the definition of the entries of the other columns.)
Clearly, the determinant $\det A_6$ is equal to the product $(\det
A_6^{(1)})\cdot(\det A_6^{(2)})$. Using \eqref{eq:ortho3}, we have
\begin{multline} \label{eq:detA61}
\det A_6^{(1)}=(-1)^{(n^2-1)/8}q^{-(n^2-1)/4}m^{(n+1)/2}
(q^{(n^2-1)/4})^{-n/2+1/2}\\
\times
\prod _{1\le h<t\le (n+1)/2}
^{}(q^{n-2h+1}-q^{n-2t+1})(1-q^{2n-2h-2t+2}),
\end{multline}
whereas by means of \eqref{eq:sympl} we have
\begin{multline} \label{eq:detA62}
\det A_6^{(2)}=q^{-(n^2-1)/4}
(q^{(n^2-1)/4})^{-(n-1)/2}\\
\times
\prod _{1\le h<t\le (n-1)/2} ^{}(q^{n-2h+1}-q^{n-2t+1})(1-q^{2n-2h-2t+2})
\prod _{h=1} ^{(n-1)/2}(q^{2n-4h+2}-1).
\end{multline}
If we now combine \eqref{eq:detA6}, \eqref{eq:detA61},
\eqref{eq:detA62}, and \eqref{eq:Nenn2}, use the fact that $\det
A_5=\det A_6$, and substitute all this in
\eqref{eq:det/det2}, then, after some simplification, we obtain
\begin{align*} 
\so&^{odd}_{\big(((m-n)/2)^n\big)}(q^{n-1},q^{n-3},\dots,q^{-n+3},q^{-n+1})\\
&\kern1cm
=\frac {\displaystyle m^{(n+1)/2}}
{\displaystyle 
\prod _{1\le h<t\le (n-1)/2} ^{}(q^{h-t}-q^{t-h})^2
\prod _{1\le h,t\le (n-1)/2} ^{}(q^{n-h-t+1}-q^{-n+h+t-1})
}\\
&\kern8.5cm
\times
\prod _{h=1} ^{(n-1)/2}
\frac 
{(q^{\frac {n+1} {2}-h}+q^{-\frac {n+1} {2}+h})} 
{(q^{\frac {n+1} {2}-h}-q^{-\frac {n+1} {2}+h})^3},
\end{align*}
 which can be rewritten as \eqref{eq:so1b}.
\end{proof}

The next lemma provides a similar evaluation of a rectangularly shaped
odd orthogonal character, the special values of the arguments at which
the character is evaluated being exactly the negative values of
those in Lemma~\ref{lem:so1}. This evaluation is needed in the proof
of Theorem~\ref{thm:5-asy2}, however only in the case that $m$ is
even. For the sake of completeness, we state also the corresponding
result for odd $m$ without proof.

\begin{Lemma} \label{lem:so2}
Let $m$ and $n$ be positive integers with $m\ge n$, and
let $q=e^{\pi i/m}$. Then we have
\begin{equation} \label{eq:so2a}
\so^{odd}_{\big(((m-n)/2)^n\big)}(-q^{n-1},-q^{n-3},\dots,-q^{-n+3},-q^{-n+1})=
\frac {\displaystyle m^{n/2}
\prod _{h=1} ^{n/2}\tan\frac {(2h-1)\pi} {2m}}
{\displaystyle 2^{\binom n2}\prod _{1\le h<t\le n}
^{}\sin\frac {(t-h)\pi} {m}
}
\end{equation}
if $n$ is even (regardless of $m$), 
\begin{equation} \label{eq:so2c}
\so^{odd}_{\big(((m-n)/2)^n\big)}(-q^{n-1},-q^{n-3},\dots,-q^{-n+3},-q^{-n+1})\\
=\frac {\displaystyle (-1)^{(m-n)/2} m^{(n-1)/2}
\prod _{h=1} ^{(n-1)/2}\tan\frac {h\pi} {m}}
{\displaystyle 2^{\binom n2}\prod _{1\le h<t\le n}
^{}\sin\frac {(t-h)\pi} {m}}
\end{equation}
if both $n$ and $m$ are odd, and
\begin{equation} \label{eq:so2b}
\so^{odd}_{\big(((m-n)/2)^n\big)}(-q^{n-1},-q^{n-3},\dots,-q^{-n+3},-q^{-n+1})
=0
\end{equation}
if $n$ is odd and $m$ is even.
\end{Lemma}

\begin{proof}
Again, when we would directly specialize in the
definition \eqref{e12} of the odd orthogonal character,
then we face the difficulty that we obtain an indeterminate
expression $0/0$. As we mentioned before the statement of the theorem,
we are only going to discuss the case that $m$ is even. The arguments
are however completely analogous if $m$ is odd.

Let first $n$ be even. In that case we must compute the limit as
\begin{equation} \label{eq:lim3}
x_1\to -q^{n-1},\ 
x_2\to -q^{n-3},\ \dots,\
x_{n/2}\to -q
\end{equation}
of the right-hand side of \eqref{e12} with $x_{n/2+1}=-q^{-1}$, \dots,
$x_{n-1}=-q^{-n+3}$, $x_n=-q^{-n+1}$. Upon a simultaneous
rearrangement and change of signs 
of the rows of the two determinants in \eqref{e12}, and using the
equality $-1=q^m$, one sees that this is equivalent to computing the limit
\eqref{eq:lim3} of $(\det A_7)/(\det B_7)$, where
\begin{equation*} \label{eq:detA7}
A_7=\begin{pmatrix} 
x_h^{\frac {m+n+1} {2}-t}-x_h^{-(\frac {m+n+1} {2}-t)}&
1\le h\le n/2\\
q^{(2n+m-2h+1)(\frac {m+n+1} {2}-t)}-q^{-(2n+m-2h+1)(\frac {m+n+1} {2}-t)}&
n/2<h\le n
\end{pmatrix}_{1\le h,t\le n}
\end{equation*}
and
\begin{equation*} \label{eq:detB6a}
B_7=\begin{pmatrix} 
x_h^{n-t+\frac {1} {2}}-x_h^{-(n-t+\frac {1} {2})}&
1\le h\le n/2\\
q^{(2n+m-2h+1)(n-t+\frac {1} {2})}-q^{-(2n+m-2h+1)(n-t+\frac {1} {2})}&
n/2<h\le n
\end{pmatrix}_{1\le h,t\le n}.
\end{equation*}
Now applying de l'H\^opital's rule, this limit is equal to the limit
\eqref{eq:lim3} of 
\begin{equation} \label{eq:det/det6}
\frac {\det A_8} {
\frac {\partial} {\partial x_1}
\frac {\partial} {\partial x_2}\cdots
\frac {\partial} {\partial x_{n/2}}
\det B_7},
\end{equation}
where
\begin{align} 
\notag
A_8&=\begin{pmatrix} 
({\frac {m+n+1} {2}-t})
(q^{(n+m-2h+1)(\frac {m+n+1} {2}-t-1)}+q^{-(n+m-2h+1)(\frac {m+n+1}
{2}-t+1)})&
1\le h\le n/2\\
q^{(2n+m-2h+1)(\frac {m+n+1} {2}-t)}-q^{-(2n+m-2h+1)(\frac {m+n+1} {2}-t)}&
n/2<h\le n
\end{pmatrix}\\
\label{eq:detA8}
&=\begin{pmatrix} 
 (-1)^{\frac {m} {2}-h-t}q^{-(n-2h+1)}({\frac {m+n+1} {2}-t})\hfill\\
\times
( q^{(n-2h+1)(\frac {n+1} {2}-t)}+ q^{-(n-2h+1)(\frac {n+1}
{2}-t)})&
1\le h\le n/2\\
 (-1)^{\frac {n+m} {2}-h-t+1}
(q^{(2n-2h+1)(\frac {n+1} {2}-t)}-q^{-(2n-2h+1)(\frac {n+1} {2}-t)})&
n/2<h\le n
\end{pmatrix}.
\end{align}
The limit of the denominator of \eqref{eq:det/det6} is again readily
obtained, because the determinant of $B_7$ can actually be evaluated
by means of \eqref{eq:ortho2}. 
The result is that the denominator of \eqref{eq:det/det6} is equal to
\begin{multline} \label{eq:Nenn3}
(q^{n^2/2})^{-n+1/2}
\prod _{1\le h<t\le n/2}
^{}(q^{n-2t+1}-q^{n-2h+1})^4(1-q^{2n-2h-2t+2})^4\\
\times
\prod _{h=1} ^{n/2}(q^{n-2h+1}+1)^2(1-q^{2n-4h+2}).
\end{multline}

Next we devote ourselves to the evaluation of the determinant of
$A_8$, with $A_8$ given by \eqref{eq:detA8}. For $t=1,2,\dots,\frac
{n} {2}$ we subtract column $n-t+1$ from column $t$. As a result we
obtain that $\det A_8$ is equal to $\det A_9$, where $A_9$ is the
$n\times n$ block matrix
\begin{equation} \label{eq:detA9}
A_9=\begin{pmatrix} A_9^{(1)}&*\\
0&A_9^{(2)}\end{pmatrix},
\end{equation}
with $A_9^{(1)}$ the $\frac {n} {2}\times \frac {n} {2}$ matrix
$$A_9^{(1)}=\Big(
 (-1)^{\frac {m} {2}-h-t} q^{-(n-2h+1)}m
( q^{(n-2h+1)(\frac {n+1} {2}-t)}+ q^{-(n-2h+1)(\frac {n+1}
{2}-t)})\Big)_{1\le h,t\le n/2},$$
and $A_9^{(2)}$ the $\frac {n} {2}\times \frac {n} {2}$ matrix
$$A_9^{(2)}=\Big(
 (-1)^{\frac {m} {2}-h-t+1}
(q^{(n-2h+1)(\frac {1} {2}-t)}-q^{-(n-2h+1)(\frac {1} {2}-t)})
\Big)_{1\le h,t\le n/2}.$$
Clearly, the determinant $\det A_9$ is equal to the product $(\det
A_9^{(1)})\cdot(\det A_9^{(2)})$. Using \eqref{eq:ortho1}, we have
\begin{multline} \label{eq:detA91}
\det A_9^{(1)}=(-1)^{\binom {n/2+1}2}q^{-n^2/4}m^{n/2}
(q^{n^2/4})^{-n/2+1/2}\\
\times
\prod _{1\le h<t\le n/2} ^{}(q^{n-2h+1}-q^{n-2t+1})
(1-q^{2n-2h-2t+2})\prod _{h=1}
^{n/2}(q^{n-2h+1}+1),
\end{multline}
whereas by means of \eqref{eq:ortho2} we have
\begin{multline} \label{eq:detA92}
\det A_9^{(2)}=
(q^{n^2/4})^{1/2-n/2}\\
\times
\prod _{1\le h<t\le n/2} ^{}(q^{n-2h+1}-q^{n-2t+1})(1-q^{2n-2h-2t+2})
\prod _{h=1} ^{n/2}(q^{n-2h+1}-1).
\end{multline}
If we now combine \eqref{eq:detA9}, \eqref{eq:detA91},
\eqref{eq:detA92}, and \eqref{eq:Nenn3}, use the fact that $\det
A_8=\det A_9$, and substitute all this in
\eqref{eq:det/det6}, then, after some simplification, we obtain
\begin{align*} \label{eq:so2}
\so&^{odd}_{\big(((m-n)/2)^n\big)}(q^{n-1},q^{n-3},\dots,q^{-n+3},q^{-n+1})\\
&=\frac {\displaystyle m^{n/2}
\prod _{h=1} ^{n/2}(1-q^{n-2h+1})}
{\displaystyle \prod _{1\le h<t\le n/2}
^{}(q^{h-t}-q^{t-h})^2
\prod _{1\le h,t\le n/2}
(q^{-n+h+t-1}-q^{n-h-t+1})
\prod _{h=1} ^{n/2}(1+q^{n-2h+1})
},
\end{align*}
which can be rewritten as \eqref{eq:so2a}.

\medskip
Now let $n$ be odd. We proceed in a completely analogous manner. Again, the
task is to compute the specialized odd orthogonal character in
\eqref{eq:orth}, by means of a limit of its definition \eqref{e12}. 

What we must compute is the limit as
\begin{equation} \label{eq:lim4}
x_1\to -q^{n-1},\ 
x_2\to -q^{n-3},\ \dots,\
x_{(n-1)/2}\to -q^2,
\end{equation}
of the right-hand side of \eqref{e12} with $x_{(n+1)/2}=-1,
x_{(n+3)/2}=-q^{-2}$, \dots,
$x_{n-1}=-q^{-n+3}$, $x_n=-q^{-n+1}$. Upon a simultaneous 
rearrangement and change of signs 
of the rows of the two determinants in \eqref{e12},
one sees that this is equivalent to computing the limit
\eqref{eq:lim4} of $(\det A_{10})/(\det B_{10})$, where
\begin{equation*} \label{eq:detA10}
A_{10}=\begin{pmatrix} 
x_h^{\frac {m+n+1} {2}-t}-x_h^{-(\frac {m+n+1} {2}-t)}&
1\le h\le (n-1)/2\\
q^{(2n+m-2h)(\frac {m+n+1} {2}-t)}-q^{-(2n+m-2h)(\frac {m+n+1} {2}-t)}&
(n-1)/2<h\le n
\end{pmatrix}_{1\le h,t\le n}
\end{equation*}
and
\begin{equation*} \label{eq:detB10a}
B_{10}=\begin{pmatrix} 
x_h^{n-t+\frac {1} {2}}-x_h^{-(n-t+\frac {1} {2})}&
1\le h\le (n-1)/2\\
q^{(2n+m-2h)(n-t+\frac {1} {2})}-q^{-(2n+m-2h)(n-t+\frac {1} {2})}&
(n-1)/2<h\le n
\end{pmatrix}_{1\le h,t\le n}.
\end{equation*}
The determinant $\det B_{10}$ can again be evaluated explicitly by
means of \eqref{eq:ortho2}, and it turns out to be non-zero. However,
the determinant $\det A_{10}$ vanishes because all the entries in
its $n$-th row are 0 (because of $q^m=-1$). 
Therefore the quotient $(\det A_{10})/(\det B_{10})$
is zero, 
whence also its limit \eqref{eq:lim4}. Thus, the claim 
\eqref{eq:so2b} is established.
\end{proof}

Our next lemma provides special evaluations of so-called {\it Schur
functions}, which are needed in the proof of
Theorem~\ref{thm:6-asy2}. Recall that
for any partition $\la=(\la_1,\la_2,\dots,\la_N)$ 
(i.e., a non-increasing sequence of non-negative integers) 
the Schur function\break  $s_\la(x_1,x_2,\dots,x_N)$ is defined
by (see \cite[p.~403, (A.4)]{FuHaAA}, \cite[I, (3.1)]{MacdAC},
or \cite[Prop.~1.4.4]{LascAZ})
\begin{equation} \label{e4}
s_\lambda (x_1, x_2,\dots ,x_N)=\frac {
 \det\limits _{1\le h,t\le N}(x_h^{\lambda _t+N-t})} {
\det\limits_{1\le h,t\le N}(x_h^{N-t})}.
\end{equation}
Again, it is not difficult to see that the denominator in \eqref{e4} does in
fact cancel out, so that any Schur function
$s_\la(x_1,x_2,\dots,x_N)$ is in fact a {\it polynomial\/} in
$x_1,x_2,\dots,x_N$, and is thus well-defined for any choice of the
variables $x_1,x_2,\dots,x_N$.

In the lemma below, the notation $\big((m-n-1)^n\big)$ is a short
notation for the partition in which the first $n$ parts are $m-n-1$,
followed by $n+1$ parts all of which are 0.

\begin{Lemma} \label{lem:Schur}
Let $q=e^{\pi i/m}$. Then we have
\begin{equation} \label{eq:Schur1a}
s_{\big((m-n-1)^n\big)}(q^{n},q^{n-1},\dots,q,-1,q^{-1},\dots,q^{-n+1},q^{-n})
=
2^{-n^2}\frac {\prodl _{h=1} ^{n/2}\tan^2\frac {(2h-1)\pi} {2m}}
{\prodl _{h=1} ^{n+1}\prodl _{t=1} ^{n}\left\vert\sin\frac {(2t-2h+1)\pi}
{2m}\right\vert}
\end{equation}
if both $m$ and $n$ are even, 
\begin{equation} \label{eq:Schur1b}
s_{\big((m-n-1)^n\big)}(q^{n},q^{n-1},\dots,q,-1,q^{-1},\dots,q^{-n+1},q^{-n})
=
2^{-n^2}\frac {\prodl _{h=1} ^{(n+1)/2}\tan^2\frac {(2h-1)\pi} {2m}}
{\prodl _{h=1} ^{n+1}\prodl _{t=1} ^{n}\left\vert\sin\frac {(2t-2h+1)\pi}
{2m}\right\vert}
\end{equation}
if $m$ is even and $n$ is odd,
\begin{multline} \label{eq:Schur2a}
s_{\big((m-n-1)^n\big)}(q^{n},q^{n-1},\dots,q,-1,q^{-1},\dots,q^{-n+1},q^{-n})\\
=
2^{-n^2}\prodl _{h=1} ^{n/2}\frac
{\sin^2\frac {(2h-1)\pi} {2m}}
{\cos^2\frac {h\pi} {m}}
\frac {1}
{\prodl _{h=1} ^{n+1}\prodl _{t=1} ^{n}\left\vert\sin\frac {(2t-2h+1)\pi}
{2m}\right\vert}
\end{multline}
if $m$ is odd and $n$ is even, and 
\begin{equation} \label{eq:Schur2b}
s_{\big((m-n-1)^n\big)}(q^{n},q^{n-1},\dots,q,-1,q^{-1},\dots,q^{-n+1},q^{-n})
=0
\end{equation}
if both $m$ and $n$ are odd.
\end{Lemma}

\begin{proof}
We have to evaluate \eqref{e4} with $N=2n+1$,
$\la=(m-n-1,m-n-1,\dots,m-n-1,0,\dots,0)$ (where $m-n-1$ is repeated
$n$ times), $x_1=q^n$, $x_2=q^{n-1}$, \dots, $x_n=q$, $x_{n+1}=-1$,
$x_{n+2}=q^{-1}$, \dots, $x_{2n}=q^{-n+1}$, $x_{2n+1}=q^{-n}$. The
denominator of \eqref{e4} can be easily evaluated as it is just a
Vandermonde determinant. If we use that $q^m=-1$,
under these specializations the numerator becomes
\begin{equation} \label{eq:num}
\det\begin{pmatrix} D_{11}&D_{12}\\
l_{1}&l_2\\
D_{21}&D_{22}\end{pmatrix},
\end{equation}
where $D_{11}$ is the $n\times n$ matrix
\begin{equation*} \label{eq:D11}
D_{11}=\begin{pmatrix} (-1)^{n+1-h}(q^{n+1-h})^{n-t}\end{pmatrix}_{1\le h,t\le
n},
\end{equation*}
$D_{12}$ is the $n\times (n+1)$ matrix
\begin{equation*} \label{eq:D12}
D_{12}=\begin{pmatrix} (q^{n+1-h})^{n+1-t}\end{pmatrix}_{1\le h\le n,\,1\le t\le
n+1},
\end{equation*}
$l_1$ is the (row) vector
\begin{equation*} \label{eq:l1}
l_1=\begin{pmatrix} (-1)^{m+n-t}\end{pmatrix}_{1\le t\le
n},
\end{equation*}
$l_2$ is the (row) vector
\begin{equation*} \label{eq:l2}
l_2=\begin{pmatrix} (-1)^{n+1-t}\end{pmatrix}_{1\le t\le
n+1},
\end{equation*}
$D_{21}$ is the $n\times n$ matrix
\begin{equation*} \label{eq:D21}
D_{21}=\begin{pmatrix} (-1)^{h}(q^{-h})^{n-t}\end{pmatrix}_{1\le h,t\le
n},
\end{equation*}
and $D_{22}$ is the $n\times (n+1)$ matrix
\begin{equation*} \label{eq:D22}
D_{22}=\begin{pmatrix} (q^{-h})^{n+1-t}\end{pmatrix}_{1\le h\le n,\,1\le t\le
n+1}.
\end{equation*}
For the evaluation of the determinant \eqref{eq:num}, we have to
distinguish between four cases, depending on the parities of $m$ and
$n$. 

If both $m$ and $n$ are even, we subtract column $n+1+t$ from column
$t$, $t=1,2,\dots,n$. The result is that the first $n$ entries in
the odd numbered rows become zero. If we rearrange the rows so that
these zeroes are moved to the bottom (i.e., row $2$ is moved up to
first position, row $4$ is moved up to second position, etc.), then
we obtain that the determinant \eqref{eq:num} is equal to
\begin{equation} \label{eq:num1}
(-1)^{\binom n2}\det\begin{pmatrix} E_{1}&*\\
0&E_{2}\end{pmatrix},
\end{equation}
where $E_{1}$ is the $n\times n$ matrix
$$
E_{1}=\begin{pmatrix} 2(-1)^{n+1-2h}(q^{n+1-2h})^{n-t}\end{pmatrix}_{1\le h,t\le
n},
$$
and $E_{2}$ is the $(n+1)\times (n+1)$ matrix
$$
E_{2}=\begin{pmatrix} (q^{n+2-2h})^{n+1-t}&h\ne \frac {n} {2}+1\\
(-1)^{n+1-t}&h=\frac {n} {2}+1\end{pmatrix}_{1\le h,t\le
n+1}.
$$
Clearly, the determinant in \eqref{eq:num1} (and thus the determinant
in \eqref{eq:num} that we want to evaluate) is equal to the product
$(\det E_1)\cdot (\det E_2)$. Both of the latter determinants are
Vandermonde determinants, and are therefore easily evaluated. 
If the result is substituted in \eqref{e4}, together with the
denominator evaluation, the claimed expression \eqref{eq:Schur1a} is
obtained after some simplification.

In the other three cases we proceed in a similar manner. There,
in the determinant \eqref{eq:num} we {\it add\/} column $n+1+t$ to column
$t$, $t=1,2,\dots,n$. If the parities of $m$ and $n$ are different,
then, again, the rows can be rearranged so that a
block form is obtained (this rearrangement is different in the two
cases), in which both the upper-left and the
lower-right blocks are Vandermonde matrices. Finally, if both 
$m$ and $n$ are odd, then, after addition of the columns as described
above, the first $n$ entries in $n+2$ (!) rows become zero (to be
precise, these are the odd numbered rows {\it and\/} the $(n+1)$-st
row). Thus, the determinant of this matrix vanishes.
\end{proof}

Next, we turn to two further special evaluations of Schur functions,
which are needed in the proofs of Theorems~\ref{thm:8-asy2} and
\ref{thm:9-asy2}. These special evaluations are given in 
Lemmas~\ref{lem:SchurA} and \ref{lem:SchurB} below). 
In the proofs, we make use of the standard basic hypergeometric
notation introduced earlier in \eqref{eq:qhyp}, and
in particular of the following summation formula.

\begin{Lemma} \label{lem:Hyp}
For any non-negative integer $N$ and any indeterminate $b$, we have
\begin{equation} \label{eq:Hyp}
{} _{2} \phi _{1} \! \left [   \begin{matrix} \let \over / \def\frac#1#2{#1
  / #2} {q^{-2 N}}, b\\ \let \over / \def\frac#1#2{#1 / #2} {{{q^{4 -
  2 N}}}\over b}\end{matrix} ;{q^2}, {\displaystyle {\frac{{q^3}} b}} \right ]
={\frac {({\let \over / \def\frac#1#2{#1 / #2} {b\over q}}; q) _{N} \,
     ({\let \over / \def\frac#1#2{#1 / #2} q^2}; q^2) _{N} }
   {{q^N}\,({\let \over / \def\frac#1#2{#1 / #2} {b\over {{q^2}}}}; {q^2})
      _{N}\, (q;q)_N }}.
\end{equation}
\end{Lemma} 
\begin{proof}We write the hypergeometric series on the left-hand side
of \eqref{eq:Hyp} as a sum over $k$,
say. The reader should observe that, because of the upper parameter 
$q^{-2N}$, the sum is in fact a finite sum, with $k$ running from
$0$ up to $N$.
Now we reverse the order of summation in the sum,
i.e., we replace $k$ by $N-k$. If we rewrite the resulting sum in
basic hypergeometric notation, then we obtain
$$
{\frac {\left( 1 - b {q^{2N-2}} \right)  
}
   {{q^N}\left( 1 - {b {{q^{-2}}}} \right)  }}
     {} _{2} \phi _{1} \! \left [             \begin{matrix} \let \over /
      \def\frac#1#2{#1 / #2} {q^{-2 N}}, {b\over {{q^2}}}\\ \let \over /
      \def\frac#1#2{#1 / #2} {{{q^{2 - 2 N}}}\over b}\end{matrix} ;{q^2},
      {\displaystyle {\frac {{q^3}} b}} \right ].
$$
This series can be summed by means of the summation formula
(see (\cite[Ex.~1.8]{GaRaAA})
$${} _{2} \phi _{1} \! \left [ \begin{matrix} \let \over / a^2, {{a^2}\over b}\\
   \let \over / b\end{matrix} ;q^2, {\displaystyle {\frac {b q} {a^2}}}
   \right ] =  \frac {1} {2}
   \left( {\frac {(\let \over / -{b\over a} ;q)
_\infty\,(a;q)_\infty\,  (q ;q^2) _\infty} 
       {\let \over / (b;q^2)_\infty\, ({{b q}\over {a^2}} ;q^2) _\infty} } + 
        {\frac {(\let \over / {b\over a} ;q) _\infty\,
(-a;q)_\infty\, (q ;q^2) _\infty}  
       {\let \over / (b;q^2)_\infty\, ({{b q}\over {a^2}} ;q^2) _\infty} } \right) ,
$$
upon letting $a$ tend to $q^{-N}$ and replacing $b$ by $q^{2-2N}/b$.
After some simplification, one arrives at the right-hand side of
\eqref{eq:Hyp}.
\end{proof}

Similar to earlier notational conventions,
in the lemma below, the notation $(c^{n-p},(c-1)^p)$ is a short
notation for the partition in which the first $n-p$ parts are $c$,
the next $p$ parts are $c-1$, followed by $n+1$ parts all of which
are $0$.

\begin{Lemma} \label{lem:SchurA}
Let $n$ and $c$ be positive integers, let $p$ be a non-negative integer
with $0\le p\le n$, and let $q$ be an
indeterminate. Then we have
\begin{multline} \label{eq:SchurA1}
s_{\big(c^{n-p},(c-1)^p\big)}(q^{2n-1},q^{2n-3},\dots,q^3,q,1,q^{-1},q^{-3},
\dots,q^{-2n+3},q^{-2n+1})\\
=\prod _{h=1} ^{2n}\frac {\(q^{\frac {c+h} {2}}-q^{-\frac {c+h} {2}}\)} 
{\(q^{\frac {h} {2}}-q^{-\frac {h} {2}}\)}
\prod _{h=1} ^{n}\prod _{t=1} ^{n}
\frac {\(q^{c+n+t-h}-q^{-c-n-t+h}\)} 
{\(q^{n+t-h}-q^{-n-t+h}\)}
\prod _{h=1} ^{n}
\frac {\(q^{h}-q^{-h}\)^2} 
{\(q^{c+p+h}-q^{-c-p-h}\)}\\
\times
\frac {1} {\displaystyle\prod _{h=1} ^{p}\(q^h-q^{-h}\)
\prod _{h=1} ^{n-p}\(q^h-q^{-h}\)}
\frac {(q^{\frac {c} {2}}-q^{-\frac {c} {2}})
(q^{\frac {c} {2}+p}+q^{-\frac {c} {2}-p})} 
{(q^{c+p}-q^{-c-p})}.
\end{multline}
\end{Lemma}

\begin{proof}
Using the symmetry of the Schur function, we have to evaluate \eqref{e4} with $N=2n+1$,
$\la=(c,\dots,c,c-1,\dots,c-1,0,\dots,0)$ (where $c$ is repeated
$n-p$ times and $c-1$ is repeated $p$ times), $x_1=q^{2n-1}$, $x_2=q^{2n-3}$, \dots,$x_{n-1}=q^3$, $x_n=q$, 
$x_{n+1}=q^{-1}$, $x_{n+2}=q^{-3}$, \dots, $x_{2n-1}=q^{-2n+3}$,
$x_{2n}=q^{-2n+1}$,
$x_{2n+1}=1$. The
denominator of \eqref{e4} can be easily evaluated as it is just a
Vandermonde determinant. On the other hand,
under these specializations the numerator becomes
\begin{equation*} 
\det\begin{pmatrix} F_{1}&F_{2}&F_{3}\\
l_{3}&l_4&l_{5}
\end{pmatrix},
\end{equation*}
where $F_{1}$ is the $(2n)\times (n-p)$ matrix
\begin{equation*} 
F_{1}=\begin{pmatrix} (q^{2n+1-2h})^{c+2n+1-t}\end{pmatrix}_{1\le h\le
2n,\,1\le t\le n-p},
\end{equation*}
$F_{2}$ is the $(2n)\times p$ matrix
\begin{equation*} 
F_{2}=\begin{pmatrix} (q^{2n+1-2h})^{c+n+p-t}\end{pmatrix}_{1\le h\le
2n,\,1\le t\le p},
\end{equation*}
$F_{3}$ is the $(2n)\times (n+1)$ matrix
\begin{equation*} 
F_{3}=\begin{pmatrix} (q^{2n+1-2h})^{n+1-t}\end{pmatrix}_{1\le h\le 2n,\,1\le t\le
n+1},
\end{equation*}
$l_3$ is the (row) vector of length $n-p$ consisting entirely of $1$'s,
$l_4$ is the (row) vector of length $p$ consisting entirely of $1$'s,
and $l_5$ is the (row) vector of length $n+1$ consisting entirely of
$1$'s.

We consider first the case where $p$ is strictly between $0$ and $n$,
that is, where $0<p<n$.
We subtract the $(t+1)$-st column from the $t$-th column,
$t=1,2,\dots,2n$. Clearly, this makes the last row become
$(0,0,\dots,0,1)$. So, if we subsequently expand the determinant with
respect to the last row and factor $q^{2n+1-2h}-1$ out of the
$h$-th row, $h=1,2,\dots,2n$, we obtain the expression 
\begin{equation} \label{eq:numA}
\bigg(\prod _{h=1} ^{2n}(q^{2n+1-2h}-1)\bigg)
\det\begin{pmatrix} G_{1}&G_2&G_3&G_4&G_{5}
\end{pmatrix},
\end{equation}
where $G_{1}$ is the $(2n)\times (n-p-1)$ matrix
\begin{equation*} 
G_{1}=\begin{pmatrix} (q^{2n+1-2h})^{c+2n-t}\end{pmatrix}_{1\le h\le
2n,\,1\le t\le n-p-1},
\end{equation*}
$G_{2}$ is the $(2n)\times 1$ matrix (i.e., column of length $2n$)
\begin{equation*} 
G_{2}=\begin{pmatrix} (q^{2n+1-2h})^{c+n+p-1}\(1+q^{2n+1-2h}\)\end{pmatrix}_{1\le h\le
2n},
\end{equation*}
$G_{3}$ is the $(2n)\times (p-1)$ matrix
\begin{equation*} 
G_{3}=\begin{pmatrix} (q^{2n+1-2h})^{c+n+p-t-1}\end{pmatrix}_{1\le h\le
2n,\,1\le t\le p-1},
\end{equation*}
$G_{4}$ is the $(2n)\times 1$ matrix (i.e., column of length $2n$)
\begin{equation*} 
G_{4}=\begin{pmatrix} (q^{2n+1-2h})^n
\displaystyle\sum _{k=0} ^{c-1}q^{(2n+1-2h)k}\end{pmatrix}_{1\le h\le 2n},
\end{equation*}
and $G_{5}$ is the $(2n)\times n$ matrix
\begin{equation*} 
G_{5}=\begin{pmatrix} (q^{2n+1-2h})^{n-t}\end{pmatrix}_{1\le h\le 2n,\,1\le t\le
n}.
\end{equation*}
We now use linearity in the $(n-p)$-th and $n$-th columns
to convert \eqref{eq:numA} into
\begin{multline} \label{eq:numA1}
\bigg(\prod _{h=1} ^{2n}(q^{2n+1-2h}-1)\bigg)
\Bigg(\sum _{k=0} ^{c-1}
\det\begin{pmatrix} G_{1}&G_2^{(1)}&G_3&G_4^{(k)}&G_{5}
\end{pmatrix}\\+
\sum _{k=0} ^{c-1}
\det\begin{pmatrix} G_{1}&G_2^{(2)}&G_3&G_4^{(k)}&G_{5}
\end{pmatrix}\Bigg),
\end{multline}
where $G_1$, $G_3$ and $G_5$ are as above, 
where $G_2^{(1)}$ is the 
$(2n)\times 1$ matrix 
\begin{equation*} 
G_{2}^{(1)}=\begin{pmatrix} (q^{2n+1-2h})^{c+n+p}
\end{pmatrix}_{1\le h\le 2n},
\end{equation*}
where $G_2^{(2)}$ is the 
$(2n)\times 1$ matrix 
\begin{equation*} 
G_{2}^{(2)}=\begin{pmatrix} (q^{2n+1-2h})^{c+n+p-1}
\end{pmatrix}_{1\le h\le 2n},
\end{equation*}
and where $G_4^{(k)}$ is the 
$(2n)\times 1$ matrix 
\begin{equation*} 
G_{4}^{(k)}=\begin{pmatrix} (q^{2n+1-2h})^{n+k}
\end{pmatrix}_{1\le h\le 2n}.
\end{equation*}
The determinants in \eqref{eq:numA1} are determinants of
the form 
\begin{equation} \label{eq:Vand}
\det\(X_t^{2n+1-2h}\)_{1\le h,t\le 2n}=
\bigg(\prod _{t=1} ^{2n}x_t^{2n-1}\bigg)\det\((X_t^{-2})^{h-1}\)_{1\le h,t\le 2n}.
\end{equation}
Since the last determinant is a Vandermonde determinant, we can
evaluate the determinants in \eqref{eq:numA1}.
Thus, under our specializations, the numerator in \eqref{e4} becomes 
\begin{multline} \label{eq:numA2}
\bigg(\prod _{h=1} ^{2n}(q^{2n+1-2h}-1)\bigg)
\bigg(\prod _{1\le h<t\le n} ^{}(q^{t-h}-q^{-t+h})^2\bigg)
\bigg(\prod _{h=1} ^{n}\prod _{t=1}
^{n}(q^{c+n+t-h}-q^{-c-n-t+h})\bigg)\\
\times
\left(
\sum _{k=0} ^{c-1}
\frac{\displaystyle
\prod _{h=1} ^{n}(q^{c+n-h-k}-q^{-c-n+h+k})
\prod _{h=1} ^{n}(q^{h+k}-q^{-h-k})}
{\displaystyle
 (q^{c+p-k-1}-q^{-c-p+k+1})
 \prod _{h=1} ^{p-1}(q^h-q^{-h}) \prod _{h=1} ^{n-p}(q^h-q^{-h})
 \prod _{h=1} ^{n}(q^{c+p+h-1}-q^{-c-p-h+1})
} \right.\\\left.
+
\sum _{k=0} ^{c-1}
\frac{\displaystyle
\prod _{h=1} ^{n}(q^{c+n-h-k}-q^{-c-n+h+k})
\prod _{h=1} ^{n}(q^{h+k}-q^{-h-k})}
{\displaystyle
 (q^{c+p-k}-q^{-c-p+k})
 \prod _{h=1} ^{p}(q^h-q^{-h}) \prod _{h=1} ^{n-p-1}(q^h-q^{-h})
 \prod _{h=1} ^{n}(q^{c+p+h}-q^{-c-p-h})
}\right).
\end{multline}
The task now is to simplify the two sums. In the first sum in
\eqref{eq:numA2} we replace $k$ by $k-1$. After this replacement,
we put the two sums together, that is, for any fixed $k$ 
we add the $k$-th summands of the two sums. 
Thus, the expression \eqref{eq:numA2} simplifies to
\begin{multline} \label{eq:numA3}
\frac{\displaystyle
\bigg(\prod _{h=1} ^{2n}(q^{2n+1-2h}-1)\bigg)
\bigg(\prod _{1\le h<t\le n} ^{}(q^{t-h}-q^{-t+h})^2\bigg)
\bigg(\prod _{h=1} ^{n}\prod _{t=1}
^{n}(q^{c+n+t-h}-q^{-c-n-t+h})\bigg)}
{\displaystyle
 \bigg(\prod _{h=1} ^{p}(q^h-q^{-h}) \bigg)
 \bigg(\prod _{h=1} ^{n-p}(q^h-q^{-h})\bigg)
 \bigg(\prod _{h=0} ^{n}(q^{c+p+h}-q^{-c-p-h})\bigg)
}\\
\times
\Bigg(
\sum _{k=0} ^{c}
\Bigg(\prod _{h=1} ^{n-1}(q^{c+n-h-k}-q^{-c-n+h+k})
(q^{h+k}-q^{-h-k})\Bigg)
(q^n-q^{-n})\kern4cm\\
\cdot\Big(
q^{-k}(q^{-c-n-p}+q^{c+n-p}-q^{c+p-n}-q^{c+n+p})
-q^{k}(q^{-c-n-p}+q^{-c+n-p}-q^{-c+p-n}-q^{c+n+p})
\Big)
\Bigg).
\end{multline}
In what follows, we concentrate on the sum in \eqref{eq:numA3}.
We split the sum in two sums according to the additive decomposition of the
last factor of the summand and rewrite the two sums in
basic hypergeometric notation,
\begin{multline} \label{eq:numA4}
(-1)^nq^{cn-c-n}(q^{-2n-2c+2};q^2)_{n-1}\,(q^{2};q^2)_n\\
\times\Bigg(
(q^{-c-n-p}+q^{c+n-p}-q^{c+p-n}-q^{c+n+p})\
{} _{2} \phi _{1} \! \left [   \begin{matrix} 
q^{-2c},q^{2n}\\q^{2-2n-2c}
\end{matrix} ;{q^2}, {\displaystyle q^{1-2n}}\right ]\\
-(q^{-c-n-p}+q^{-c+n-p}-q^{-c+p-n}-q^{c+n+p})\
{} _{2} \phi _{1} \! \left [   \begin{matrix} 
q^{-2c},q^{2n}\\q^{2-2n-2c}
\end{matrix} ;{q^2}, {\displaystyle q^{3-2n}}\right ]
\Bigg).
\end{multline}
We write the second $_2\phi_1$-series as a sum over $k$, say, 
reverse the order of summation, that is, we replace $k$ by $c-k$, and
then we write the resulting sum again in basic hypergeometric
notation. We obtain a $_2\phi_1$-series which turns out to be
identical with the first $_2\phi_1$-series in \eqref{eq:numA4}.
If we put everything together and simplify, then \eqref{eq:numA4}
becomes
\begin{multline} \label{eq:numA5}
(-1)^nq^{cn-2c-2n-p}(q^{-2n-2c+2};q^2)_{n-1}\,(q^{2};q^2)_n
(1-q^c)(1-q^{c+2n})(1+q^{c+2p})\\
\times
{} _{2} \phi _{1} \! \left [   \begin{matrix} 
q^{-2c},q^{2n}\\q^{2-2n-2c}
\end{matrix} ;{q^2}, {\displaystyle q^{1-2n}}\right ].
\end{multline}
Next we apply the contiguous relation
$$
{} _{2} \phi _{1} \! \left [ \begin{matrix} \let \over / a, b\\ \let \over /
   c\end{matrix} ;q, {\displaystyle z} \right ] =
  {} _{2} \phi _{1} \! \left [ \begin{matrix} \let \over / a, b q\\ \let \over /
    c\end{matrix} ;q, {\displaystyle z} \right ] - 
 b z 
   {\frac {        (1 - a )  }        {(1 - c)}}
{} _{2} \phi _{1} \! \left [ \begin{matrix} \let \over / a q, b q\\ \let
        \over / c q\end{matrix} ;q, {\displaystyle z} \right ]$$
to the $_2\phi_1$-series in \eqref{eq:numA5}.
This transforms the expression \eqref{eq:numA5} into
\begin{multline} \label{eq:numA6}
(-1)^nq^{cn-2c-2n-p}(q^{-2n-2c+2};q^2)_{n-1}\,(q^{2};q^2)_n
(1-q^c)(1-q^{c+2n})(1+q^{c+2p})\\
\times
\left(
{} _{2} \phi _{1} \! \left [   \begin{matrix} 
q^{-2c},q^{2n+2}\\q^{2-2n-2c}
\end{matrix} ;{q^2}, {\displaystyle q^{1-2n}}\right ]
-q\frac {1-q^{-2c}} {1-q^{2-2n-2c}}
{} _{2} \phi _{1} \! \left [   \begin{matrix} 
q^{-2c+2},q^{2n+2}\\q^{4-2n-2c}
\end{matrix} ;{q^2}, {\displaystyle q^{1-2n}}\right ]
\right).
\end{multline}
Both $_2\phi_1$-series in the last line can be evaluated by means of
Lemma~\ref{lem:Hyp}. If we substitute the
result in \eqref{eq:numA6}, put this back in \eqref{eq:numA3},
and divide the result by the numerator in \eqref{e4}
subject to our specializations, which we evaluated by means of the
Vandermonde determinant evaluation, we arrive at the right-hand side
of \eqref{eq:SchurA1} after some simplification.

If $p=0$, then we can proceed in a completely analogous manner. In
fact, the computations are somewhat simpler in this case, so that we
leave the details to the reader. That the formula works also for
$p=n$ can then be checked by verifying that \eqref{eq:SchurA1} for $p=n$
agrees with \eqref{eq:SchurA1} for $p=0$ and $c$ replaced by $c-1$.

This completes the proof of the lemma.
\end{proof}

\begin{Lemma} \label{lem:SchurB}
Let $n$ and $c$ be positive integers, and let $q$ be an
indeterminate. Then we have
\begin{multline} \label{eq:SchurB1}
s_{(c^n)}(q^{2n-1},q^{2n-3},\dots,q^3,q,-1,q^{-1},q^{-3},
\dots,q^{-2n+3},q^{-2n+1})\\
=\prod _{h=1} ^{2n}
\frac {\(q^{\frac {c+h} {2}}-(-1)^{c+h}q^{-\frac {c+h} {2}}\)} 
{\(q^{\frac {h} {2}}-(-1)^hq^{-\frac {h} {2}}\)}
\prod _{h=1} ^{n-1}\prod _{t=1} ^{n}
\frac {\(q^{c+n+t-h}-q^{-c-n-t+h}\)} 
{\(q^{n+t-h}-q^{-n-t+h}\)}.
\end{multline}
\end{Lemma}

\begin{proof}
We replace $q$ by $-q$ in the $p=0$ case of Lemma~\ref{lem:SchurA}.
Then, the left-hand side of \eqref{eq:SchurA1} agrees exactly with
the left-hand side of \eqref{eq:SchurB1}. The products on the
right-hand sides are not in the same form, but they are equivalent
because the extra terms in Lemma~\ref{lem:SchurA} cancel out if
$p=0$.
\end{proof}

The final two lemmas concern sums of {\it even orthogonal characters}, which
are needed in the proofs of Theorems~\ref{thm:10-asy2} and
\ref{thm:11-asy2}.
Given a non-increasing sequence 
$\la=(\la_1,\la_2,\dots,\la_n)$ of integers or half-integers with
$\la_{n-1}\ge\vert\la_n\vert$, the even orthogonal character\break
$\so^{even}_\la(x_1,x_2,\dots,x_n)$ is
defined by (see \cite[(24.40)]{FuHaAA})
\begin{equation} \label{e25}
\so^{even}_\la(x_1,x_2,\dots,x_n)
=\frac {\det\limits_{1\le h,t\le
n}(x_h^{\la_t+n-t}+x_h^{-(\la_t+n-t)})+
\det\limits_{1\le h,t\le
n}(x_h^{\la_t+n-t}-x_h^{-(\la_t+n-t)})} 
{\det\limits_{1\le h,t\le
n}(x_h^{n-t}+x_h^{-(n-t)})}.
\end{equation}

\begin{Lemma} \label{lem:spin}
Let $n$ be a positive integer and $c$ be a non-negative
half-integer. Then we have
\begin{multline} \label{eq:spinSumme}
\sum _{p=-c} ^{c}\so^{even}_{\big((c^{n-1},p)\big)}(q^{n-1},q^{n-2},\dots,1)
\\=\prod_{1\le h<t\le n}\frac{q^{(2c+t+h-1)/2}-q^{-(2c+t+h-1)/2}}
                     {q^{(t+h-2)/2}-q^{-(t+h-2)/2}}
     \frac{\displaystyle\prod_{h=1}^{n}\(q^{c+h-1/2}-q^{-(c+h-1/2)}\)}
     {\displaystyle\prod_{h=1}^{n-1}\(q^{h/2}-q^{-h/2}\)}\\
\times
     \sum_{k=1}^{n}\frac {(-1)^{k-1}(2c+2k-1)}
        {\displaystyle \prod_{h=1}^{n} \(q^{(2c+k+h-1)/2}-q^{-(2c+k+h-1)/2}\)
         \prod_{h=1}^{k-1}\(q^{h/2}-q^{-h/2}\)
         \prod_{h=1}^{n-k}\(q^{h/2}-q^{-h/2}\)}.
\end{multline}
\end{Lemma}

\begin{proof}For convenience, let us write
$\soe_\la(x_1,x_2,\dots,x_n)$ for just ``one half" in the definition
\eqref{e25} of even orthogonal characters, that is,
\begin{equation} \label{e26}
\soe_\la(x_1,x_2,\dots,x_n)
=\frac {\det\limits_{1\le h,t\le
n}(x_h^{\la_t+n-t}+x_h^{-(\la_t+n-t)})} 
{\det\limits_{1\le h,t\le
n}(x_h^{n-t}+x_h^{-(n-t)})}.
\end{equation}
It should be observed that we have
\begin{equation} \label{eq:spinspin}
\soe_\la(x_1,\dots,x_{n-1},1)=
\so^{even}_\la(x_1,\dots,x_{n-1},1),
\end{equation}
since the second determinant in \eqref{e25} vanishes if $x_n=1$.
Now, by \eqref{e26} and
the determinant evaluation \eqref{eq:ortho3}, we have
\begin{align} \notag 
\sum _{p=-c} ^{c}&\soe_{\big((c^{n-1},p)\big)}(x_1,x_2,\dots,x_n)\\
\notag
&=\frac {1} {2}
\Bigg(\prod_{1\le h<t\le n} ^{}\dfrac {1} {(x_h^{\frac {1} {2}}x_t^{-\frac {1} {2}}
-x_h^{-\frac {1} {2}}x_t^{\frac {1} {2}})(x_h^{\frac {1} {2}}x_t^{\frac
{1} {2}}-x_h^{-\frac {1} {2}}x_t^{-\frac
{1} {2}})}\Bigg)\\
\notag
&\kern1cm
\times
\sum _{p=-c} ^{c}\det_{1\le h,t\le n}\begin{pmatrix} 
x_h^{c+n-t}+x_h^{-(c+n-t)}&1\le t<n\\
x_h^p+x_h^{-p}&t=n\end{pmatrix}\\
\notag
&=\frac {1} {2}
\Bigg(\prod_{1\le h<t\le n} ^{}\dfrac {1} {(x_h^{\frac {1} {2}}x_t^{-\frac {1} {2}}
-x_h^{-\frac {1} {2}}x_t^{\frac {1} {2}})(x_h^{\frac {1} {2}}x_t^{\frac
{1} {2}}-x_h^{-\frac {1} {2}}x_t^{-\frac
{1} {2}})}\Bigg)\\
\notag
&\kern1cm
\times\det_{1\le h,t\le n}\begin{pmatrix} 
x_h^{c+n-t}+x_h^{-(c+n-t)}&1\le t<n\\
2\dfrac {x_h^{c+\frac {1} {2}}-x_h^{-(c+\frac {1} {2})}}
{x_h^{\frac {1} {2}}-x_h^{-\frac {1} {2}}}&t=n\end{pmatrix}\\
\notag
&=
\Bigg(\prod_{1\le h<t\le n} ^{}\dfrac {1} {(x_h^{\frac {1} {2}}x_t^{-\frac {1} {2}}
-x_h^{-\frac {1} {2}}x_t^{\frac {1} {2}})(x_h^{\frac {1} {2}}x_t^{\frac
{1} {2}}-x_h^{-\frac {1} {2}}x_t^{-\frac
{1} {2}})}\Bigg)
\Bigg(\prod _{h=1} ^{n}\frac {1} {x_h^{\frac {1} {2}}-x_h^{-\frac {1}
{2}}}\Bigg)\\
\notag
&\kern1cm
\times
\det_{1\le h,t\le n}\begin{pmatrix} 
x_h^{c+n-t+\frac {1} {2}}-x_h^{-(c+n-t+\frac {1} {2})}
+x_h^{c+n-t-\frac {1} {2}}-x_h^{-(c+n-t-\frac {1} {2})}&1\le t<n\\
{x_h^{c+\frac {1} {2}}-x_h^{-(c+\frac {1} {2})}}
&t=n\end{pmatrix}\\
\notag
&=
\Bigg(\prod_{1\le h<t\le n} ^{}\dfrac {1} {(x_h^{\frac {1} {2}}x_t^{-\frac {1} {2}}
-x_h^{-\frac {1} {2}}x_t^{\frac {1} {2}})(x_h^{\frac {1} {2}}x_t^{\frac
{1} {2}}-x_h^{-\frac {1} {2}}x_t^{-\frac
{1} {2}})}\Bigg)
\Bigg(\prod _{h=1} ^{n}\frac {1} {x_h^{\frac {1} {2}}-x_h^{-\frac {1}
{2}}}\Bigg)\\
\notag
&\kern1cm
\times
\det_{1\le h,t\le n}\(
x_h^{c+n-t+\frac {1} {2}}-x_h^{-(c+n-t+\frac {1} {2})}\),
\end{align}
where the last line arises by obvious elementary column operations
from the next-to-last line. In this identity, we want to put
$x_1=q^{n-1}$, $x_2=q^{n-2}$, \dots, $x_n=1$. However, we cannot
directly put $x_n=1$, because this would lead to an expression $0/0$.
Therefore instead, we have to take the limit $x_n\to1$. Thereby,
using \eqref{eq:spinspin}, we get
\begin{align} \notag
\sum _{p=-c}
^{c}&\so^{even}_{\big((c^{n-1},p)\big)}(q^{n-1},q^{n-2},\dots,1)\\
\notag
&=
\Bigg(\prod_{1\le h<t\le n} ^{}\dfrac {1} {
(q^{(t-h)/2}-q^{-(t-h)/2})(q^{(h+t-2)/2}-q^{-(h+t-2)/2}}\Bigg)
\Bigg(\prod _{h=1} ^{n-1}\frac {1} {q^{h/2}-q^{-h/2}}\Bigg)\\
\label{eq:spin1}
&\kern1cm
\times
\det_{1\le h,t\le n}\begin{pmatrix}
q^{(n-h)(c+n-t+\frac {1} {2})}-q^{-(n-h)(c+n-t+\frac {1} {2})}
&1\le h<n\\
2c+2n-2t+1
&h=n\end{pmatrix}.
\end{align}
Now we expand the determinant along the last row. Then each of the
appearing minors can be evaluated by means of \eqref{eq:sympl}. 
If we substitute the result in \eqref{eq:spin1}, we obtain
\begin{align*} 
\sum _{p=-c}
^{c}&\so^{even}_{\big((c^{n-1},p)\big)}(q^{n-1},q^{n-2},\dots,1)\\
&=
\Bigg(\prod_{1\le h<t\le n} ^{}\dfrac {1} {
(q^{(t-h)/2}-q^{-(t-h)/2})(q^{(h+t-2)/2}-q^{-(h+t-2)/2})}\Bigg)
\Bigg(\prod _{h=1} ^{n-1}\frac {1} {q^{h/2}-q^{-h/2}}\Bigg)\\
&\kern1cm
\times
\sum _{k=1} ^{n}(-1)^{n+k}(2c+2n-2k+1)
\underset{t\ne k}{\prod _{t=1} ^{n}}
(q^{c+n-t+\frac {1} {2}}-q^{-(c+n-t+\frac {1} {2})})
\\
&\kern2cm\cdot
\underset{h,t\ne k}{\prod
_{1\le h<t\le n} ^{}}
(q^{(t-h)/2}-q^{-(t-h)/2})(q^{(2c+2n-h-t+1)/2}-q^{-(2c+2n-h-t+1)/2}).
\end{align*}
After replacing $k$ by $n+1-k$ and performing
some simplification, we arrive at the claimed expression.
\end{proof}

\begin{Lemma} \label{lem:spin-1}
Let $n$ be a positive integer and $c$ be a non-negative
half-integer. Then we have
\begin{multline} \label{eq:spinSumme-1}
\sum _{p=-c} ^{c}(-1)^{c-p}
\so^{even}_{\big((c^{n-1},p)\big)}(q^{n-1},q^{n-2},\dots,1)
\\=\prod_{1\le h<t\le n}\frac{q^{(2c+t+h-1)/2}-q^{-(2c+t+h-1)/2}}
                     {q^{(t+h-2)/2}-q^{-(t+h-2)/2}}
\prod_{h=1}^{n-1}     \frac1
     {\displaystyle\(q^{h/2}+q^{-h/2}\)}.
\end{multline}
\end{Lemma}

\begin{proof}According to the definition \eqref{e25} of even orthogonal
characters and the determinant evaluation \eqref{eq:ortho3}, we have
\begin{align*} \notag 
\sum _{p=-c} ^{c}&(-1)^{c-p}\soe_{\big((c^{n-1},p)\big)}(x_1,x_2,\dots,x_n)\\
\notag
&=\frac {1} {2}
\Bigg(\prod_{1\le h<t\le n} ^{}\dfrac {1} {(x_h^{\frac {1} {2}}x_t^{-\frac {1} {2}}
-x_h^{-\frac {1} {2}}x_t^{\frac {1} {2}})(x_h^{\frac {1} {2}}x_t^{\frac
{1} {2}}-x_h^{-\frac {1} {2}}x_t^{-\frac
{1} {2}})}\Bigg)\\
\notag
&\kern1cm
\times
\sum _{p=-c} ^{c}(-1)^{c-p}\det_{1\le h,t\le n}\begin{pmatrix} 
x_h^{c+n-t}+x_h^{-(c+n-t)}&1\le t<n\\
x_h^p+x_h^{-p}&t=n\end{pmatrix}\\
\notag
&=\frac {1} {2}
\Bigg(\prod_{1\le h<t\le n} ^{}\dfrac {1} {(x_h^{\frac {1} {2}}x_t^{-\frac {1} {2}}
-x_h^{-\frac {1} {2}}x_t^{\frac {1} {2}})(x_h^{\frac {1} {2}}x_t^{\frac
{1} {2}}-x_h^{-\frac {1} {2}}x_t^{-\frac
{1} {2}})}\Bigg)\\
\notag
&\kern1cm
\times\det_{1\le h,t\le n}\begin{pmatrix} 
x_h^{c+n-t}+x_h^{-(c+n-t)}&1\le t<n\\
2\dfrac {x_h^{c+\frac {1} {2}}+x_h^{-(c+\frac {1} {2})}}
{x_h^{\frac {1} {2}}+x_h^{-\frac {1} {2}}}&t=n\end{pmatrix}\\
\notag
&=
\Bigg(\prod_{1\le h<t\le n} ^{}\dfrac {1} {(x_h^{\frac {1} {2}}x_t^{-\frac {1} {2}}
-x_h^{-\frac {1} {2}}x_t^{\frac {1} {2}})(x_h^{\frac {1} {2}}x_t^{\frac
{1} {2}}-x_h^{-\frac {1} {2}}x_t^{-\frac
{1} {2}})}\Bigg)
\Bigg(\prod _{h=1} ^{n}\frac {1} {x_h^{\frac {1} {2}}+x_h^{-\frac {1}
{2}}}\Bigg)\\
\notag
&\kern1cm
\times
\det_{1\le h,t\le n}\begin{pmatrix} 
x_h^{c+n-t+\frac {1} {2}}+x_h^{-(c+n-t+\frac {1} {2})}
+x_h^{c+n-t-\frac {1} {2}}+x_h^{-(c+n-t-\frac {1} {2})}&1\le t<n\\
{x_h^{c+\frac {1} {2}}+x_h^{-(c+\frac {1} {2})}}
&t=n\end{pmatrix}\\
\notag
&=
\Bigg(\prod_{1\le h<t\le n} ^{}\dfrac {1} {(x_h^{\frac {1} {2}}x_t^{-\frac {1} {2}}
-x_h^{-\frac {1} {2}}x_t^{\frac {1} {2}})(x_h^{\frac {1} {2}}x_t^{\frac
{1} {2}}-x_h^{-\frac {1} {2}}x_t^{-\frac
{1} {2}})}\Bigg)
\Bigg(\prod _{h=1} ^{n}\frac {1} {x_h^{\frac {1} {2}}+x_h^{-\frac {1}
{2}}}\Bigg)\\
&\kern1cm
\times
\det_{1\le h,t\le n}\(
x_h^{c+n-t+\frac {1} {2}}+x_h^{-(c+n-t+\frac {1} {2})}\),
\end{align*}
where $\soe_\la(x_1,x_2,\dots,x_n)$ is defined in \eqref{e26}, and
where the last line arises by obvious elementary column operations
from the next-to-last line. Now we put 
$x_1=q^{n-1}$, $x_2=q^{n-2}$, \dots, $x_n=1$ in this identity,
in which case, due to \eqref{eq:spinspin}, the ``incomplete" even orthogonal character
$\soe_{\big((c^{n-1},p)\big)}(x_1,x_2,\dots,x_n)$ becomes the
specialized even orthogonal character 
$\so^{even}_{\big((c^{n-1},p)\big)}(q^{n-1},q^{n-2},\dots,1)$,
and we use
\eqref{eq:ortho3} to evaluate the determinant in the last line.
After some simplification we arrive at \eqref{eq:spinSumme-1}.
\end{proof}

\end{document}